\documentclass[10pt,reqno]{amsart}
\usepackage{enumerate}
\usepackage{amsfonts}
\usepackage{amsmath}
\usepackage{amsthm}
\usepackage{latexsym}
\usepackage{multicol}
\usepackage{verbatim}
\usepackage{amscd}
\usepackage{amssymb}
\usepackage{amsmath, amscd,wasysym,xcolor}
\usepackage{graphicx}
\usepackage{booktabs,caption,fixltx2e}
\usepackage{setspace}
\graphicspath{ {images/} }
\usepackage{multicol}
\usepackage[colorlinks=true, linkcolor=blue, citecolor=red, urlcolor=blue]{hyperref}
\usepackage{hyperref}
\usepackage{hyperref}
\advance\textwidth by 1.2in \advance\oddsidemargin by -.65in \advance\evensidemargin by -.65in
\newtheorem*{cor}{Corollary}
\newtheorem*{lem}{Lemma}
\newtheorem*{prop}{Proposition}

\newtheorem*{saturation}{Saturation Theorem of Knutson and Tao}

\theoremstyle{definition}
\newtheorem*{defn}{Definition}

\theoremstyle{definition}
\newtheorem*{thm}{Theorem}

\newtheorem*{rem}{Remark}

\newcounter{cnt}
 \makeatletter
		\def\mydggeometry{\makeatletter\dg@YGRID=1\dg@XGRID=20\unitlength=0.003pt\makeatother}
		\makeatother \theoremstyle{remark}
		
		\numberwithin{equation}{section}

		\let\bwdg\bigwedge
		\def\bigwedge{{\textstyle\bwdg}}

		\newcommand{\thmref}[1]{Theorem~\ref{#1}}
		
		\newcommand{\secref}[1]{Section~\ref{#1}}
		\newcommand{\lemref}[1]{Lemma~\ref{#1}}
		\newcommand{\propref}[1]{Proposition~\ref{#1}}
		\newcommand{\propnref}[1]{Proposition$'$~\ref{#1}}
		
		\newcommand{\remref}[1]{Remark~\ref{#1}}

		\newcommand{\Ker}{\operatorname{Ker}}

		\newcommand{\nc}{\newcommand}
		\newcommand{\rnc}{\renewcommand}

		\nc{\cal}{\mathcal} \nc{\goth}{\mathfrak} \rnc{\bold}{\mathbf}
		\nc{\fk}{\mathfrak}
		\nc{\germ}{\mathfrak}

		\renewcommand{\Bbb}{\mathbb}
		\nc\bomega{{\mbox{\boldmath $\omega$}}}
		\nc\bxi{{\mbox{\boldmath $\xi$}}}
		\nc\blambda{{\mbox{\boldmath $\lambda$}}}
		\nc\bOmega{{\mbox{\boldmath $\Omega$}}}
		\nc\bchi{{\mbox{\boldmath $\chi$}}}\nc\bnu{{\mbox{\boldmath $\nu$}}}
		\nc\bpsi{{\mbox{\boldmath $\Psi$}}} \nc\bast{{\mbox{\boldmath $\ast$}}}
		\nc\balpha{{\mbox{\boldmath $\alpha$}}}
		\nc\bpi{{\mbox{\boldmath $\pi$}}}
		\nc\bsigma{{\mbox{\boldmath $\sigma$}}} \nc\bcN{{\mbox{\boldmath $\cal{N}$}}} \nc\bcm{{\mbox{\boldmath $\cal{M}$}}} 
		\nc\bLambda{{\mbox{\boldmath$\Lambda$}}} 
		\nc\bll{{\mbox{\boldmath$\ell$}}} \nc\bgamma{{\mbox{\boldmath$\gamma$}}}
		
		\nc\p{{\mbox{\boldmath$\rho$}}} 
		\nc\bmu{{\mbox{\boldmath$\mu$}}}
		
		\newcommand{\mfg }{\mathfrak{g}}
		
		\newcommand{\lie}[1]{\mathfrak{#1}}
		\makeatletter
\def\subsection{\def\@secnumfont{\bfseries}\@startsection{subsection}{2}%
{\parindent}{.5\linespacing\@plus.5\linespacing}{-.5em}%
{\normalfont\bfseries}}
\makeatother

\nc{\Hom}{\operatorname{Hom}}
\nc{\td}{\operatorname{\tilde{d}}}\nc{\D}{\operatorname{d}}
\nc{\mode}{\operatorname{mod}}
\nc{\End}{\operatorname{End}} \nc{\wh}[1]{\widehat{#1}} \nc{\Ext}{\operatorname{Ext}} \nc{\ch}{\text{ch}} \nc{\ev}{\operatorname{ev}}
\nc{\Ob}{\operatorname{Ob}} \nc{\soc}{\operatorname{soc}} \nc{\rad}{\operatorname{rad}} \nc{\head}{\operatorname{head}}
\nc{\A}{\operatorname{\bold{b}}}

\def\ker{\operatorname{Ker}}

\nc{\Cal}{\cal} \nc{\Xp}[1]{X^+(#1)} \nc{\Xm}[1]{X^-(#1)}
\nc{\on}{\operatorname} \nc{\Z}{{\mathbb{Z}}} \nc{\J}{{\cal J}} \nc{\C}{{\mathbb{C}}} \nc{\Q}{{\bold Q}} \nc{\R}{{\mathbb{R}}}
\nc{\K}{\bold{\kappa}}

\nc{\F}{\operatorname{\boF}} 

\nc{\N}{{\Bbb N}} \nc\boa{\bold a} \nc\bob{\bold b} \nc\boc{\bold c} \nc\bod{\bold d} \nc\boe{\bold e} \nc\bof{\bold f} \nc\bog{\bold g}
\nc\boh{\bold h} \nc\boi{\bold i} \nc\boj{\bold j} \nc\bok{\bold k} \nc
\bol{\bold l} \nc\bom{\bold m} \nc\bon{\bold n} \nc\boo{\bold o}
\nc\bop{\bold p} \nc\boq{\bold q} \nc\bor{\bold r} \nc\bos{\bold s} \nc\boT{\bold t} \nc\boF{\bold F} \nc\bou{\bold u} \nc\bov{\bold v}
\nc\bow{\bold w} \nc\boz{\bold z} \nc\boy{\bold y} \nc\ba{\bold A} \nc\bb{\bold B} \nc\bc{\bold C} \nc\bd{\bold D} \nc\be{\bold E} \nc\bg{\bold
G} \nc\bh{\bold H} \nc\bi{\bold I} \nc\bj{\bold J} \nc\bk{\bold K} \nc\bl{\bold L} \nc\bm{\bold M} \nc\bn{\bold N} \nc\bo{\bold O} \nc\bp{\bold
P} \nc\bq{\bold Q} \nc\br{\bold R} \nc\bs{\bold S} \nc\bt{\bold T} \nc\bu{\bold U} \nc\bv{\bold V} \nc\bw{\bold W} \nc\bz{\bold Z} \nc\bx{\bold
x} \nc\KR{\bold{KR}} \nc\rk{\bold{rk}} \nc\het{\text{ht }}

\nc\fn{{fin}}  \nc\af{{aff}}  \nc\tr{{tor}} \nc\btilde{\bold{\tilde{\bold{H}}}}

\nc{\mpp}{\rotatebox[origin=c]{180}{\pm}}
\nc{\iA}{\rotatebox[origin=c]{180}{A}}
\nc{\iB}{\rotatebox[origin=c]{180}{B}}
\nc{\iC}{\rotatebox[origin=c]{180}{C}}

\nc\eps{\epsilon}
\nc\toa{\tilde a} \nc\tob{\tilde b} \nc\toc{\tilde c} \nc\tod{\tilde d} \nc\toe{\tilde e} \nc\tof{\tilde f} \nc\tog{\tilde g} \nc\toh{\tilde h}
\nc\toi{\tilde i} \nc\toj{\tilde j} \nc\tok{\tilde k} \nc\tol{\tilde l} \nc\tom{\tilde m}  \nc\ton{\tilde n} \nc\too{\tilde o} \nc\toq{\tilde q}
\nc\tor{\tilde r} \nc\tos{\tilde s} \nc\toT{\tilde t} \nc\tou{\tilde u} \nc\tov{\tilde v} \nc\tow{\tilde w} \nc\toz{\tilde z}
\begin{document}

\title[Graded Character Formula for Fusion Products and LR-Coefficients in Type $A_2$]{Graded Character Formula for Fusion Products of Irreducible Modules and Littlewood-Richardson Coefficients in Type $A_2$}
\author{Tanusree Khandai}
\address{ Indian Institute of Science Education and Research, Mohali, Punjab, India}
\email{tanusree@iisermohali.ac.in }
\author{Shushma Rani}
\address{Indian Institute of Science Education and Research, Mohali, Punjab, India}
\email{shushmarani95@gmail.com, ph16067@iisermohali.ac.in.}



\begin{abstract}

We investigate a specific class of CV modules for $\mathfrak{sl}_3$ and establish an exact sequence for these modules. Utilizing dimension arguments, we demonstrate that this module is isomorphic to the fusion product of irreducible modules, thereby offering a new proof of the conjecture regarding the independence of fusion products from parameters. By analyzing the filtration of the kernel within the exact sequence, we derive the graded character formula for fusion product modules. Moreover, we leverage the graded character to deduce the algebraic characterization of the Littlewood-Richardson (LR) coefficients and present an alternative proof of the saturation theorem in type $A_2$.
\end{abstract}

\maketitle
\section{Introduction}

Let $\mathfrak g$ be a finite-dimensional simple Lie algebra over the complex field and $\mathfrak g[t]$ the corresponding current algebra. 
In this paper, we study a class of graded finite-dimensional representations of $\mathfrak{sl}_3[t]$ modules called fusion product modules which are known to generalize several important classes of representations of the current Lie algebras like the local Weyl modules, Demazure modules and Kirillov-Reshitikhin modules.

For any $\boa =(a_1, a_2,\cdots,a_k)\in\mathbb C^k$, the tensor product of the irreducible $\mathfrak g$-modules $V(\lambda_1)$, $V(\lambda_2), \cdots, $ $V(\lambda_k)$ with the highest weight $\lambda_1, \cdots, \lambda_k$, acquires  the structure of a $\mathfrak g[t]$-module called an  {\it{evaluation module}} as follows:
$$x\otimes P(t).v_1\otimes v_2\otimes \cdots \otimes v_k =\sum_{i=1}^k P(a_i) v_1\otimes \cdots v_{i-1}\otimes x.v_i\otimes v_{i+1}\otimes\cdots v_k. $$  We denote this $\mathfrak g[t]$-module by $V(\blambda,\boa)$, where $\blambda:=(\lambda_1, \lambda_2, \cdots, \lambda_k)$. When the $a_i$'s are all distinct, $V(\blambda,\boa)$ is an irreducible $\mathfrak g[t]$-module. Constructing a $\mathfrak g$-equivariant filtration on an irreducible $\lie g[t]$-module $V(\blambda,\boa)$, the associated fusion product $ V(\lambda_1)^{a_1}\ast \cdots \ast V(\lambda_k)^{a_k}$, was introduced in \cite{FL}. These modules  had been proposed as a way of constructing the generalized versions of the Kostka polynomials, and it was conjectured that as a $\lie g[t]$-module, they are independent of the point of evaluation. While the conjecture has not been completely resolved, it has been shown to hold in several special cases \cite{CL,fourlitt2007, Naoi2012, F.homogeneous,fmm17,kl, cv15, Naoi2012, BK}.

Given a dominant integral weight $\lambda=\sum\limits_{i=1}^n r_i\omega_i$, a $\mathbb{Z}_+$-graded cyclic $\mathfrak g[t]$-module $W_{loc}(\lambda)$, with the highest weight $\lambda$ was introduced in \cite{CPWeyl}. This module is called the local Weyl module with the highest weight $\lambda$. It is defined via generators and relations, and it has the property that every finite-dimensional cyclic $\mathfrak g[t]$ module generated by the highest weight vector of weight $\lambda$ is a quotient of the local Weyl module $W_{loc}(\lambda)$. By construction, fusion product of evaluation modules are highest weight cyclic representations of current Lie algebras and  hence quotients of local Weyl modules. Further the underlying vector space of the fusion product modules and consequently their dimensions are known. 
In \cite{CL}, by constructing explicit bases for local Weyl modules over $\mathfrak{sl}_{n+1}[t]$, it was shown that the dimension of 
$W_{loc}(\lambda)$ is equal to the dimension of the fusion product $V(\omega_1)^{z_{11}}\ast \cdots \ast V(\omega_1)^{z_{1r_1}}\ast \cdots \ast V(\omega_n)^{z_{nr_n}}$, thereby proving that the fusion product of fundamental representations of $\mathfrak{sl}_{n+1}$ with highest weight $\lambda$,   is in fact isomorphic to $W_{loc}(\lambda)$. This established the conjecture in \cite{CL} for the fusion product of fundamental representations of $\mathfrak{sl}_{n+1}$. In \cite{cv15},  suitable quotients of local Weyl modules called the Chari-Venkatesh (CV) modules were introduced. Once again by determining an explicit basis for this new family of modules, it was proved that every fusion product module over $\mathfrak{sl}_2[t]$ is isomorphic to a Chari-Venkatesh module and this helped obtain a  reproof of the conjecture in \cite{CL} for $\lie g$ of type $A_1$.

The Chari-Venkatesh modules are a family of finite-dimensional quotients of the local Weyl modules were introduced in \cite{cv15}. These modules are associated with a family of partitions $\bxi$ that are indexed by the set of positive roots of $\mathfrak g$, and they satisfy certain natural compatibility conditions.
They have garnered significant interest in recent years for their remarkable connections to fusion products of irreducible modules, as well as their role in investigating the independence of fusion products from parameters. In this paper, we focus on a specific class of Chari-Venkatesh modules over $\mathfrak{sl}_3[t]$.

Our main objective is to explore the properties of CV-modules over $\mathfrak{sl}_3[t]$ and establish a deeper understanding of their structure. To this end, we present an exact sequence for these modules, offering a novel perspective on their interrelationships. By employing dimension arguments, we demonstrate the isomorphism between CV modules and certain fusion products, thereby providing a compelling reproof of the conjecture regarding the independence of fusion products from parameters. This conjecture for two irreducible modules of type $A_2$ has already been resolved in \cite{BK}. However our approach of study is different.
In fact, in \secref{example}, via an example, we  show that our methods can even be extended to understand the structure of fusion products of more than two irreducible $\mathfrak{sl}_3[t]$-modules. 

Further, we investigate the filtration of the kernel within the exact sequence, leading to the derivation of the graded character formula for fusion product modules. This formula enables us to obtain an algebraic characterization of the Littlewood-Richardson (LR) coefficients - a fundamental topic in the representation theory of symmetric groups. The insights gained from the graded character contribute to an alternate proof of the saturation theorem in type $A_2$. In summary, this manuscript presents a comprehensive study of CV-modules over $\mathfrak{sl}_3[t]$, unveiling their intricate relationships with fusion products and LR coefficients.

The paper is structured as follows. In Section \ref{prelim}, we establish the notations and recall the definitions of the local Weyl module and the fusion product module. In Section \ref{CV-mod}, we review the definition and properties of the Chari-Venkatesh modules introduced in \cite{cv15, k18}. 
In Section \ref{fusion-sec}, starting with the statement of the main \thmref{filtration} we study the Chari-Venkatesh modules related to a pair of dominant integral weights of $\mathfrak{sl}_3$. We show that a series of short exact sequences of $\mathfrak{sl}_3[t]$-modules can be associated with such a module. In Section \ref{sec S_lambda,mu} and \ref{proof}, we analyze these short sequences and use dimension arguments on the corresponding modules to complete the proof of the main result. In Section \ref{graded char}, we present the graded character of fusion product modules and, using it derive another definition of the Littlewood-Richardson coefficients \thmref{littlewood coeff}. In Section \ref{section 8}, we give an alternate proof of the saturated tensor product theorem for modules of type $A_2$. This research contributes to the understanding of CV-modules and their connections with fusion products and certain combinatorial identities.
\subsection*{\em{Acknowledgements}}{\em
The authors express their gratitude to Professor Vyjayanthi Chari for raising pertinent questions regarding this work, to Niranjan Nehra for the initial discussions, and to Mrigendra Singh Kushwava for the discussion on the saturated tensor product theorem. Additionally, the second author acknowledges the support of the CSIR-SRF grant with file number 09/947(0082)/2017-EMR-I.}
\section[Preliminaries]{Preliminaries}\label{prelim} Throughout this paper, $\C$ will denote the field of complex numbers and $\Z$ the set of integers, $\Z_+$ (resp. $\N$)  the set of non-negative integers (resp. positive integers).
\subsection{} Let $\mathfrak{g}= \mathfrak{sl}_{3}$ be the Lie algebra of $3 \times 3$  matrices with trace zero, $\mathfrak{h}$ be the Cartan subalgebra of $\mathfrak{g}$ consisting of diagonal matrices. Let $\Delta=\{\alpha_1, \alpha_2\}$ be the simple roots of $\mathfrak g$ with respect to $\mathfrak h$, $\theta=\alpha_1+\alpha_2$ be the longest root of $\lie g$ and  $R^+=\{\alpha_1,\alpha_2, \theta\}$ (respectively, $R^-=-R^+$) be the set of positive (negative) roots of $\mathfrak g$ with respect to $\Delta$. Let $\{\omega_i:  i=1,2\}$ be the set of fundamental weights of $\mathfrak{g}$ which is dual to $\{\alpha_i: i=1,2\}$, $P= \Z\omega_1+\Z\omega_2$ (resp. $Q=\Z\alpha_1+\Z\alpha_2$)  be the weight lattice (resp. root lattice) of $\mathfrak{g}$ and $P^+=\Z_+\omega_1+\Z_+\omega_2$ (resp. $Q^+=\Z_+\alpha_1+\Z_+\alpha_2$). For $\alpha\in R^+$, let $\mathfrak g_{\alpha}$ be the corresponding root space. For $\alpha_i \in \Delta$, fix non-zero elements $x_{\alpha_i}^\pm\in \mathfrak g_{\pm\alpha_i}$ and $h_{\alpha_i}\in \mathfrak h$ such that 
$$[h_{\alpha_i},x_{\alpha_i}^\pm]=\pm2 x^\pm_{\alpha_i} \hspace{.25cm}[h_{\alpha_i},x^\pm_{\alpha_i}]=\pm 2x^\pm_{\alpha_i}.$$
Then $\mathfrak g_{\pm{\theta}}$ is spanned by $x^\pm_{\theta}=\pm [x_{\alpha_1}^{\pm},x_{\alpha_2}^\pm]$ and $h_{\theta}=[x_{\theta}^+,x_{\theta}^-]=h_{\alpha_1}+h_{\alpha_2}$. Set
$\mathfrak n^\pm = \underset{\alpha\in R^+}{\oplus} \mathfrak g_{\pm \alpha}.$

For $\lambda \in P^+$, let $V(\lambda)$ be the irreducible finite-dimensional $\mathfrak g$-module with the highest weight vector $v_\lambda$ of weight  $\lambda$, i.e, the cyclic module generated by $v_\lambda$ with defining relations:
\begin{equation}\label{irr.rel}
\mathfrak n^+.v_\lambda =0, \hspace{.25cm} h.v_\lambda=\lambda(h)v_\lambda, \hspace{.25cm} (x_{\alpha_i}^-)^{\lambda(h_{\alpha_i})+1}v_\lambda =0, \quad \quad \forall \,\,  h\in \mathfrak h \, \text{and} \, i=1,2. \end{equation}
Let $W$ be the Weyl group of $\mathfrak g$ and $w_0$ the longest element of $W$. The for all $\lambda\in P^+$, the module $V(-w_0\lambda)$ is the $\mathfrak g$-dual of $V (\lambda)$.

Let $\mathbb {Z}[P]$ be the group ring of $P$ with integer coefficients and basis $e(\mu)$ for $  \mu \in P$. The character of the $\lie g$-module $V(\lambda)$ is the element of $\mathbb {Z}[P]$ given by:
$$\ch_{\mathfrak{g}}V(\lambda) = \sum_{\mu \in P}\dim V(\lambda)_{\mu} e(\mu),\quad \quad \text{where} \, \, V(\lambda)_\mu = \{ v\in V(\lambda) : hv=\mu(h)v\, \forall h\in \lie h\}.$$
\subsection{} Let $\mathfrak g[t]$ be the current algebra of the Lie algebra $\mathfrak g$. As a vector space $\fk g[t]= \mathfrak{g} \otimes \mathbb {C}[t]$, where $\mathbb {C}[t]$ denotes the polynomial ring in indeterminate $t$. Clearly, $\mathfrak g[t]$ is a Lie algebra with a Lie  bracket given by 
$$ \quad  [a\otimes f, b\otimes g]=[a,b]\otimes fg, \quad \quad \forall \, a, b \in \mathfrak{g}, \, \, f,g \in \C[t].$$ $\C[t]$ has a natural $\mathbb Z_+$ grading, determined by the degree of the polynomial. Let $\bold U(\mathfrak g[t])$ be the universal enveloping algebra of $\mathfrak{g}[t]$. With the $\mathbb Z_+$-grading inherited from $ \mathfrak{g}[t]$, $\bold U(\mathfrak g[t])$ is $\mathbb Z_+$-graded. We say an element $X \in \bold U(\mathfrak{g}[t])$ has grade $r_1+r_2+\cdots+r_p$ if $X$ is of the form $(x_{\alpha_1}\otimes t^{r_1})(x_{\alpha_2}\otimes t^{r_2})\cdots(x_p\otimes t^{r_p}).$ For a positive integer $s$, we denote by $\bold{U}(\mathfrak{g}[t])[s]$ the subspace of $\bold{U}(\mathfrak{g}[t])$ spanned by the $s$-graded elements.

For every $\alpha\in R^+$, define the power series in the indeterminate $u$, 
\begin{eqnarray*}
H_\alpha(u) = \exp \left(- \sum_{r=1}^\infty \frac{h_\alpha\otimes t^r}{r} u^r\right), \end{eqnarray*} 
and for $r,s \in \mathbb N$, set
\begin{equation*}
S(r,s) = \{(b_p)_{p\geq 0}: b_p\in \mathbb Z_+, {\sum\limits_{p\geq 0}} b_p = r, \sum\limits_{p\geq 0} pb_p = s \}, \quad \quad \quad 
x(r,s) = {\sum\limits_{(b_p)_{p\geq 0}\in S(r,s)}} \prod\limits_{i=0}^s (x\otimes t^i)^{(b_i)}. \end{equation*}
The following result was proved in \cite[Lemma 7.5]{Garland} and formulated in \cite[Lemma 1.3]{CPWeyl} as follows.
\begin{lem}\label{Garland} Given $s\in \mathbb N$, $r\in \mathbb Z_+$ and $\alpha\in R^+$, we have
\begin{equation}\label{G.eqn}
	(x_\alpha^+\otimes t)^{(s)}(x_\alpha^-\otimes 1)^{(r+s)} - (-1)^s (\sum\limits_{k\geq 0} x^-_\alpha(r,r+s-k){P_\alpha(u)}_k) \in \bu(\mathfrak g[t])\mathfrak n^+[t].
\end{equation}  where
$P_\alpha(u)_k$ denotes the coefficient of $u^k$ in the power series $H_\alpha(u)$ and $(y)^{(s)} := \frac{y^s}{s!}$ for $y \in \mathfrak g[t]$.
\end{lem} 
\subsection{} A graded representation of $\mathfrak g[t]$ is a $\mathbb Z_+$-graded vector space such that $$V= \underset{r\in \mathbb Z_+}{\bigoplus} V[r], \hspace{1.5cm} \bu(\mathfrak g[t])[s].V[r]\subseteq V[r+s], \, \, \forall\, r,s\in \mathbb Z_+.$$ If $U, V$ are two graded representations of $\mathfrak g[t]$, then we say $\psi: U\rightarrow V$ is a morphism of graded $\mathfrak g[t]$-modules if  $\psi(U[r])\subseteq V[r]$ for all $r\in \mathbb Z_+$. For $s\in \mathbb  Z$, let $\tau^\ast_s$ be the grade-shifting operator given by 
$$\tau^\ast_s(V)[k] = V[k-s] \quad \quad \forall\, \, k\in \mathbb Z_+,$$ and  graded representations $V$ of $\mathfrak g[t]$. Given $z\in \mathbb C$ and a $\mathfrak g$-module $U$, let  $\ev_z(U)$ denote the corresponding evaluation module for $\mathfrak g[t]$. Clearly, for $z=0$ $\ev_0(U)$ is a graded representation and 
$\ev_0(U)[0]=U.$
\subsection{} For $\lambda \in P^+$, the local Weyl module $W_{loc}(\lambda)$, is the  $\mathfrak{g}[t]$-module generated by an element $w_{\lambda}$ with the following defining relations: 
\begin{equation}
(\mathfrak n^+ \otimes \C[t])w_\lambda= 0,\quad \,(h_{\alpha_i}\otimes t^s)w_\lambda= \lambda(h_{\alpha_i})\delta_{s,0}w_\lambda, \quad \,
(x_{\alpha_i}^-\otimes 1)^{\lambda(h_{\alpha_i})+1}w_\lambda =0,
\end{equation} 
for all $i=1,2$ and $s \in \Z_+$. Clearly,  $W_{loc}(\lambda)$ is a graded $\mathfrak{g}[t]$ module and every  finite-dimensional cyclic module with highest weight $\lambda$ is a quotient of $W_{loc}(\lambda)$. These modules were introduced in \cite{CPWeyl}. 
\label{W.loc}
\subsection{} \label{fusion}
Given a $s$-tuple of dominant integral weights $\blambda=(\lambda_1,\cdots,\lambda_s)\in (P^+)^s$ and $\bold z= (z_1,\cdots, z_s),$ a $s$-tuple of distinct complex numbers, it was proved in \cite{Cinvent} that the evaluation module
$V(\blambda,\boz) $ 
is an irreducible $\mfg[t]$ module . The $\mathbb N$-grading in $\bold{U}(\fk g[t])$ induces a $\mathfrak g$-equivariant grading on $V(\blambda,\boz)$ given 
by 
$$V(\blambda, \boz){[k]} = \underset{0\leq r\leq k}{\bigoplus} \bold{U}(\fk g[t])[r]. v_1\otimes\cdots\otimes v_s,$$ where $v_i$ is the generator of $V(\lambda_i)$ for $1\leq i\leq s$. The associated graded  $\mfg[t]$-module 
$$ V(\lambda_1)^{z_1}\ast \cdots \ast V(\lambda_s)^{z_s} := V(\blambda,\boz){[0]} \oplus \underset{k\in \mathbb N}{\bigoplus} \dfrac{\bold{V}(\blambda,\boz){[k]}}{\bold{V}(\blambda,\boz){[k-1]}} $$ is called the {\it{fusion product}} of $V(\lambda_1),\cdots, V(\lambda_s)$ at $\bold{z}$. These modules were introduced in \cite{FL}.

In this context, the following lemma is useful:
\begin{lem} Suppose that $V(\blambda,\boz)$ is an irreducible $\mathfrak g[t]$-module. 
If for $v\in \bold {V}(\blambda,\boz)$,  $\bar{v}$ denotes its image in $ V(\lambda_1)^{z_1}\ast \cdots \ast V(\lambda_s)^{z_s}$  then
$$ x\otimes t^p.\bar{v} = x\otimes (t-a_1)\cdots(t-a_p). {\bar{v}}, \hspace{.25cm} \forall\ x\in \mathfrak g, \ \text{and}\ 
a_1,\cdots, a_p\in \mathbb C.$$ 
\end{lem}
\begin{rem}\label{sec2.1}
Given a $k$-tuple of dominant integral weights $\lambda_1, \ldots, \lambda_k$, and a tuple $\bold z= (z_1,\cdots, z_k),$ of pairwise distinct complex numbers, the fusion product module $V(\lambda_1)^{z_1} \ast \ldots \ast V(\lambda_k)^{z_k}$ is a cyclic, finite dimensional, highest weight module with the highest weight $\sum\limits_{i=1}^k \lambda_i$. Hence, by the universal property of local Weyl modules, it is a  quotient of $W_{loc}(\sum\limits_{i=1}^{k} \lambda_i)$.
Since the underlying vector space of the given fusion product
is $\underset{i=1}{\overset{k}{\otimes}} V(\lambda_i)$,  the dimension of  $V(\lambda_1)^{z_1}\ast \ldots \ast V(\lambda_k)^{z_k}$ is equal to  $\prod\limits_{i=1}^{k} \dim V(\lambda_i) $.
\end{rem}
\vspace{.25cm}
\section{The Chari Venkatesh Modules}\label{CV-mod}
In this section, we recall the definition and properties of a family of quotients of local Weyl modules that were introduced by Chari and Venkatesh in \cite{cv15}. These modules are referred to as Chari Venkatesh modules ($CV$- modules). 
\begin{defn}
Given $\lambda\in P^+$, a  $|R^+|$-tuple of partitions $ {\pmb{\xi}} =({\xi}^\alpha)_{\alpha\in R^+}$ is said to be $\lambda$-compatible if,  $${\xi}^\alpha =( \xi^\alpha_1\geq \xi^\alpha_2\geq \cdots), \hspace{.35cm} \text{and} \hspace{.35cm} \sum_{i\geq 1}\xi_i^\alpha =\lambda(h_\alpha), \hspace{.25cm} \forall \ \alpha\in R^+.$$ Let $\Xi_\lambda$ be the set of all $\lambda$-compatible $R^+$-tuple of partitions. Given a $\bxi=(\xi_\alpha)_{\alpha\in R^+}\in\Xi_\lambda$, the Chari-Venkatesh module $V( {\pmb{\xi}})$ is defined as the cyclic $\mathfrak{g}[t]$-module generated by a non-zero vector  $v_{\pmb{\xi}}$ with the following defining relations:  
\begin{equation}
	(\mathfrak{n}^+ \otimes \mathbb C[t])v_{\pmb{\xi}}= 0,\quad \,(h_{\alpha_i}\otimes t^s)v_{\pmb{\xi}}= \lambda(h_{\alpha_i})\delta_{s,0}v_{\pmb{\xi}}, \quad \, 
	(x_{\alpha_i}^-\otimes 1)^{\lambda(h_{\alpha_i})+1}v_{\pmb{\xi}} =0, \quad \, \text{for}\, i=1,2, \end{equation} 
\begin{equation}
	(x^+_\alpha\otimes t)^{(s)}(x^-_\alpha\otimes1)^{(r+s)}v_{\pmb{\xi}} =0,
	\hspace{.5cm}  \alpha\in R^+, r,s\in \mathbb N, s+r\geq 1+rk + \sum_{j\geq k+1} \xi^{\alpha}_j,\ \text{for}\ k\in \mathbb N.
	\end{equation} \end{defn}
	\noindent By definition, a $CV$-module, $V({\pmb{\xi}})$, is the graded quotient of local Weyl module $W_{loc}(\lambda)$ by the submodule generated by graded elements$$\{(x^+_\alpha\otimes t)^{(s)}(x^-_\alpha\otimes1)^{(r+s)}w_{\lambda} :  \alpha\in R^+, r,s\in \mathbb N, s+r\geq 1+rk + \sum_{j\geq k+1} \xi^{\alpha}_j,\ \text{for}\ k\in \mathbb N\}.$$
	\subsection{}
For $k\in \Z_+$ and $r,s\in \mathbb N$, set, 
$$\begin{array}{c} {_k}S(r,s) = \{(b_p)_{p\geq 0}\in S(r,s): b_p=0, \text{for}\ 
p<k\},  \quad S(r,s)_{k} 
= \{(b_p)_{p\geq 0}\in S(r,s): b_p=0, \text{for}\ p\geq k\},\end{array}$$
$$\begin{array}{c} 
{_k}x(r,s) = {\sum\limits_{(b_p)_{p\geq 0}\in {_k}S(r,s)}} \prod_{i=0}^s (x\otimes t^i)^{(b_i)}, \quad
x(r,s)_{k}
=  {\sum\limits_{(b_p)_{p\geq 0}\in S(r,s)_{k}}} \prod_{i=0}^s (x\otimes t^i)^{(b_i)}\end{array}$$
The following lemma was proved in \cite[Lemma 2.5,\  Proposition 2.6]{cv15}. 
\begin{lem}\label{cv.2.5} Let $r,s,k\in \mathbb N$ and $L\in \Z_+$ such that $r+s\geq kr+L$. Then the following holds.
\begin{enumerate}
\item[(i)] For  $\alpha\in R^+$, 
$$x^-_\alpha(r,s) = {_k}x^-_{\alpha}(r,s) + \sum_{(r',s')} x^-_\alpha(r-r',s-s')_{k}\ {_k}x^-_\alpha(r',s'),$$ where sum is taken over all pairs $r',s'\in \mathbb N$ such that $r'<r, s'<s$ and $r'+s'\geq r'k+L.$
\item[(ii)] If $V$ is a $\mathfrak g[t]$-module and $v\in V$, then for $\alpha\in R^+$,
$$x^-_\alpha(r,s)v =0,\hspace{.25cm} \text{ if and only if}\hspace{.25cm} {_k}x^-_\alpha(r,s).v=0.$$ \end{enumerate}\end{lem}

The following lemma is deduced from Lemma \ref{Garland}, \cite[Proposition 2.6]{cv15}. 
\label{cv.lemma}
\begin{lem} Given a $R^+$-tuple of partitions ${\pmb{\xi}} = (\xi^{\alpha})_{\alpha\in R^+}$, let $V({\pmb{\xi}})$ be the associated $CV$-module with generator $v_{\pmb{\xi}}$. Then for all $\alpha\in R^+$, we have
\begin{align} &(x^+_\alpha\otimes t)^{(r)}(x^-_\alpha\otimes 1)^{(r+s)}.v_{\pmb{\xi}} = x_\alpha^-(r,s).v_{\pmb{\xi}},\\
&(x_{\alpha}^-(r,s)-{_1}x_{\alpha}^-(r,s))v_{\pmb{\xi}}\in \sum\limits_{r'<r} \bu(\mathfrak{sl}_3(\mathbb C))x_{\alpha}^-(r',s)v_{\pmb{\xi}},	
\end{align}
\end{lem}
\proof By part(2) of above Lemma \ref{cv.2.5}, for $r,s,k \in \mathbb{N}$ and $L \in \mathbb{Z}_+$, 
$$x^-_\alpha(r,s)v_{\pmb{\xi}} =0,\hspace{.25cm} \text{ if and only if}\hspace{.25cm} {_k}x^-_\alpha(r,s).v_{\pmb{\xi}}=0.$$
By \cite[Proposition 2.7]{cv15}, $V({\pmb{\xi}})$ be the $CV$-module is generated by $v_{\pmb{\xi}}$ with defining relations of $W_{loc}(\lambda)$ and $ {_k}x^-_\alpha(r,s).v_{\pmb{\xi}}=0.$
Therefore by definition of CV-module, $$(x^+_\alpha\otimes t)^{(r)}(x^-_\alpha\otimes 1)^{(r+s)}.v_{\pmb{\xi}} = x_\alpha^-(r,s).v_{\pmb{\xi}}. $$
Again using part(1) of Lemma \ref{cv.2.5} for $k=1$, we have $s=s'$ and
$$x^-_\alpha(r,s) = {_1}x^-_{\alpha}(r,s) + \sum_{r'} x^-_\alpha(r-r',0)_{1}\ {_1}x^-_\alpha(r',s),$$ where sum is taken over $r' \in \mathbb N$ such that $r'<r$. Since $x^-_\alpha(r-r',0)_{1} \in U(\mathfrak{sl}_3)$. Thus,
$$(x_{\alpha}^-(r,s)-{_1}x_{\alpha}^-(r,s))v_{\pmb{\xi}}\in \sum\limits_{r'<r} \bu(\mathfrak{sl}_3(\mathbb C))x_{\alpha}^-(r',s)v_{\pmb{\xi}}.$$ Hence the lemma.\endproof
\subsection{The set $P^+{(\lambda,k)} $ and module $\mathcal{F}_{\blambda}$ for $\blambda \in P^+(\lambda,k)$} 
For $\lambda \in P^+$ and a positive integer $k$, let 
$$P^+(\lambda,k)= \{\blambda=(\lambda_1,\cdots,\lambda_k)\in (P^+)^{\times k}: \sum\limits_{i=1}^{k} \lambda_i=\lambda \}.$$
Given $\blambda\in P^+(\lambda,k)$ and $\alpha \in R^+$, let $p(\lambda(h_\alpha)) = (\lambda_1(h_\alpha)^{\downarrow}\geq \ldots \geq \lambda_k(h_\alpha)^{\downarrow})$ be the $k$-tuple of integers obtained by rearranging the components of $(\lambda_1(h_\alpha), \ldots,\lambda_k(h_\alpha))$ in non-increasing order. Clearly, $p(\lambda(h_\alpha))$ is a partition of $\lambda(h_\alpha)$. We shall refer to such a partition as $\blambda_\alpha$-compatible  partition of $\lambda(h_\alpha)$ and denote it by $\blambda_{\alpha}^{\downarrow}$. 
Let ${\pmb{\xi}}_{\blambda} =(\blambda_\alpha^\downarrow)_{\alpha\in R^+}$. Clearly $\bxi_\blambda\in \Xi_\lambda$. Let $\cal F_\blambda$ be the CV module associated $\bxi_\blambda.$
\vspace{.025cm}

\noindent Define an ordering on $P^+(\lambda,k)$ as follows. Given $\blambda=(\lambda_{1},\cdots,\lambda_{k})\in P^+(\lambda, k)$ and $\bmu=(\mu_{1}, \cdots, \mu_{k})\in P^+(\lambda, k)$ we say $\blambda$ majorizes $\bmu$ and write $\blambda \succeq \bmu$ if for every $\alpha\in R^+$, 
$$\sum_{j=i}^k \lambda_{j}(h_\alpha)^\downarrow \geq \sum_{j=i}^k \mu_{j}(h_\alpha)^\downarrow, \hspace{.75cm} \text{for all}\, \ 1\leq i\leq k.$$ 
Denoting the image of $w_\lambda$ in $\mathcal F_{\blambda}$ by $v_\blambda$, we see that $\mathcal F_\blambda$ is a graded $\mathfrak g[t]$ module generated by $v_\blambda$ with defining relations:
\begin{equation} 
(\mathfrak{n}^- \otimes \mathbb C[t])v_\blambda= 0,\quad\,(h_{\alpha_i}\otimes t^s)v_\blambda= \lambda(h_{\alpha_i})\delta_{s,0}v_\blambda, \quad \, 
(x_{\alpha_i}^-\otimes 1)^{\lambda(h_{\alpha_i})+1}v_\blambda =0, \quad \, \text{for}\, i=1,2,
\end{equation}
\begin{equation}
(x_\alpha^+\otimes t)^s(x_\alpha^-\otimes 1)^{r+s} v_\blambda=0, \hspace{.25cm} 
s+r\geq 1+r\ell + \sum_{j\geq \ell+1} \lambda_j(h_\alpha)^\downarrow,  \ \forall \ \ell\in \mathbb N, \ \alpha \in R^+.
\end{equation}
\noindent It is easy to see that $\bxi: P^+(\lambda,k)\rightarrow \Xi_\lambda$ given by $\blambda\mapsto \bxi_\blambda$ is a well-defined map. Clearly, $\bxi$ is not injective for $k>2$. For example, $\blambda=(2\omega_2+\omega_1,\omega_2,2\omega_1)$ and $\bmu=(2\omega_2, \omega_1,2\omega_1+\omega_2)$ are elements of $P^+(3\omega_1+3\omega_2)$ such that $\bxi_\lambda=\bxi_\bmu$.

\subsection{} The following lemma lists some properties of the modules ${\mathcal F}_{\blambda}$ for $\blambda\in P^+(\lambda,k)$.
\begin{lem}\label{prop of F_lambda} Let $\blambda =(\lambda_1, \ldots, \lambda_k)$, $\bmu=(\mu_1, \ldots, \mu_k) \in P^+(\lambda,k)$. Then 
\begin{enumerate}
\item[(i)]  ${\cal F}_{\bmu}$ is a quotient of ${\cal F}_{\blambda}$ whenever $\blambda\succeq \bmu$. In particular, 
$\ev_0(V(\lambda))$ is the unique irreducible quotient of $\cal F_{\blambda}$, for all $\blambda\in P^+(\lambda,k)$.
\item[(ii)] $\mathfrak{g}\otimes t^k \mathbb{C}[t]. {\cal F}_{\blambda} =0 $. Hence, ${\cal F}_{\blambda}$ is a module for the Lie algebra $\mathfrak{g} \otimes \mathbb{C}[t] /(t^k)$.
\item[(iii)]\label{quotient.lem}  The fusion product module 
$ V(\lambda_1)^{z_1}\ast \cdots \ast V(\lambda_k)^{z_k}$  is a quotient of ${\mathcal F}_\blambda$ for all $k$-tuples of distinct complex numbers $(z_1,\cdots,z_k)$, . 
\end{enumerate}
\end{lem}
\proof Given $\lambda\in P^+$ and $\blambda, \bmu \in P^+(\lambda,k)$, let $v_\blambda$, $v_{\bmu}$ be the image of $w_\lambda$ in ${\cal F}_\blambda$ and ${\cal F}_\bmu$ respectively.

(i) Using definition, for all $\alpha\in R^+$,
$$(x_\alpha^+\otimes t)^{(s)}(x_\alpha^-\otimes 1)^{(r+s)}.v_{\bmu}=0, \, \text{whenever}\,   r,s,q\in \mathbb N, \, \text{satisfy}\, r+s\geq 1+qr+ \sum_{j=q+1}^k \mu_{j}(h_\alpha)^\downarrow.$$
Since $\blambda\succeq \bmu$, \quad  
$\sum_{j=q+1}^k \lambda_{j}(h_\alpha)^\downarrow \geq \sum_{j=q+1}^k \mu_{j}(h_\alpha)^\downarrow,$ for all  $0\leq q\leq k-1,\, \alpha\in R^+ $, we have 
for $r',s', q\in \mathbb N$ satisfying $r'+s'\geq 1+r'q+\sum_{j=q+1}^k \lambda_{j}(h_\alpha)^\downarrow,$
$$(x_\alpha^+\otimes t)^{(s')}(x_\alpha^-\otimes 1)^{(r'+s')}v_{\bmu}=0, \, \text{whenever}\,   r',s',q\in \mathbb N, \, \text{satisfy}\, r'+s'\geq 1+qr'+ \sum_{j=q+1}^k \lambda_{j}(h_\alpha)^\downarrow.$$ Hence  there exists a surjective $\mathfrak g[t]$-module homomorphism $\Psi^{\blambda}_{\bmu}: {\cal F}_{\blambda}\rightarrow {\cal F}_{\bmu}$,   such that $\Psi^{\blambda}_{\bmu}(v_{\blambda})= v_{\bmu}.$  Since $\blambda\succeq (\lambda,0,\cdots,0)$, for all $\blambda\in P^+(\lambda,k)$, and  by definition, $\cal F_{(\lambda,0,\cdots)}=\ev_0(V(\lambda))$, we conclude that $\ev_0(V(\lambda))$ is the unique irreducible quotient of ${\mathcal F}_{\blambda}$ for all $\blambda\in P^+(\lambda,k)$. 

(ii) Given $\blambda\in P^+(\lambda, k)$,  for all $\alpha\in R^+$, 
$\blambda_\alpha^\downarrow$ is a partition of $\lambda(h_\alpha)$ with atmost $k$ parts. Hence by  definition of ${\mathcal F}_\blambda$, we have
$$(x_\alpha^+\otimes t)^{(s)}(x_\alpha^-\otimes 1)^{(r+s)}.v_\blambda =0, \hspace{.35cm}\ r+s\geq 1+kr.$$ In particular for $r=1$, we have $x_\alpha^-\otimes t^k.v_\blambda=0 $ for all $\alpha\in R^+$, which implies $\mathfrak{g}\otimes t^k\mathbb C[t].v_\blambda =0.$

(iii) We know  $ev_{z_i}V(\lambda_i)$ is a quotient of the local Weyl module $W_{loc}(\lambda_i)$ for $1\le i\le k$. If $v_i$ is the image of $w_{\lambda_i}$ in $ev_{z_i}(V(\lambda_i))$, then $ V(\lambda_1)^{z_1}\ast \cdots \ast V(\lambda_k)^{z_k}$ is an integrable $\mathfrak g[t]$-module of highest weight $\sum\limits_{i=1}^k \lambda_i = \lambda$, generated by the vector $v_1\ast \cdots\ast v_k$ and hence is a quotient of $W_{loc}(\lambda)$. This implies,
$$ (\mathfrak{n}^+\otimes \mathbb C[t])v_1\ast\cdots\ast v_k =0, \hspace{.35cm}
(h_{\alpha_i}\otimes t^s).\lambda(h_{\alpha_i}) \delta_{s,0} v_1\ast\cdots\ast v_k =0$$
$$(x_{\alpha_i}^-\otimes 1)^{\lambda(h_{\alpha_i})+1}.v_1\ast\cdots\ast v_k =0. $$  
For $\alpha\in R^+$, if $\sigma$ is an element of the symmetric group $\mathcal{S}_k$ such that $ \blambda_{\alpha}^\downarrow= (\lambda_{\sigma(1)}(h_\alpha)\geq \cdots \geq \lambda_{\sigma(k)}(h_\alpha)),$
then using the fact that $$x_\alpha^-\otimes (t-z_{\sigma(1)})(t-z_{\sigma(2)})\cdots(t-z_{\sigma(\ell)}).ev_{z_{\sigma(1)}}V(\lambda_{\sigma(1)})\otimes \cdots\otimes ev_{z_{\sigma(\ell)}}V(\lambda_{\sigma(\ell)}) =0,$$
and $$(x_\alpha^+\otimes t)^s(x_\alpha^-\otimes 1)^{r+s}v_{\sigma(\ell+1)}\otimes\cdots \otimes v_{\sigma(k)} =0 \hspace{.25cm} \forall \ r+s\geq 1+r\ell+\sum_{j={\ell+1}}^k \lambda_{\sigma(j)}(h_\alpha),$$ the same proof as \cite[Proposition 6.8]{cv15} shows that for each $\alpha\in R^+$,
$${}_{_\ell} x_\alpha^-(r,s). v_1\ast\cdots\ast v_k =0$$ whenever $r,s,\ell \in \mathbb N$ are such that $r+s\geq 1+r\ell+\sum_{j\geq \ell+1} \lambda_j(h_\alpha)^\downarrow.$
Since $V(\lambda_1)^{z_1}\ast \cdots \ast V(\lambda_k)^{z_k}$ is a quotient of $W_{loc}(\sum\limits_{i=1}^{k}\lambda_i)$, it follows that $V(\lambda_1)^{z_1}\ast \cdots \ast V(\lambda_k)^{z_k}$ is a quotient of $\mathcal{F}_{\blambda}$.\endproof

\begin{rem}\label{f-lambda dim} It follows from Remark \ref{sec2.1} and part(iii) of Lemma \ref{prop of F_lambda} that, $\dim \mathcal F_{\blambda} \geq \prod\limits_{i=1}^{k} \dim V(\lambda_i).$
\end{rem}
\section{The CV module ${\mathcal F}_{\lambda,\mu}$ for {$\mathfrak{sl}_3[t]$}}\label{fusion-sec}

Given a dominant integral weight $\lambda=n\omega_1$ of $\mathfrak{sl}_2$, it is easy to see that the set $P^+(\lambda,k)$ is in one-to-one correspondence with the set of partitions of $n.$  In \cite{cv15}, using a series of canonical short exact sequences,  a monomial basis for $V(\xi)$ had been constructed for all $\xi\in P^+(n\omega_1,k)$. Then, using dimension arguments, it was proved that fusion product modules for $\lie{sl}_2[t]$ are isomorphic to Chari-Venkatesh modules whose definition is independent of any evaluation parameter. We 
extend this method in the case when $\mathfrak g =\mathfrak{sl}_3(\mathbb C)$ 
and $\blambda\in P^+(\lambda,2)$ for $\lambda\in P^+$. In this section, we state the main result of the paper and prove preliminary results that help establish it. 
\subsection{} 
For $\lambda= \lambda_1 \omega_1+\lambda_2 \omega_2\in P^+$, let $|\lambda|:=\lambda_1+\lambda_2$.
\begin{defn}Let $\mathfrak g:=\lie {sl}_3(\mathbb C)$. Given $\nu\in P^+$ and $(\lambda,\mu)\in P^+(\nu,2)$ with $\lambda=\lambda_1\omega_1+\lambda_2\omega_2$ and $\mu=\mu_1\omega_1+\mu_2\omega_2$, we  say: 
\begin{enumerate}
\item[i.] $(\lambda, \mu)$ is a partition of $\nu$ of first kind
if $|\lambda|\geq |\mu|$,  $\lambda_i\geq \mu_i$ for $i=1,2$;
\item[ii.]  $(\lambda, \mu)$ is a partition of $\nu$ of second kind if $|\lambda|\geq |\mu|$, $\lambda_1\geq \mu_1$ and $\mu_2>\lambda_2$.
\item[iii.]  $(\lambda, \mu)$ is a partition of $\nu$ of third kind if $|\lambda|\geq |\mu|$, $\lambda_2\geq \mu_2$ and $\mu_1>\lambda_1.$ 
\end{enumerate}
\end{defn}

Notice that if $(\lambda,\mu)$ is a partition of third kind, then the   $(-w_0\lambda,-w_0\mu):=(\lambda_2\omega_1+\lambda_1\omega_2,\mu_2\omega_1+\mu_1\omega_2)$ is a partition of second kind and the  Dynkin diagram automorphism of $\lie{sl}_3(\mathbb C)$ that maps $\alpha_1$ to $\alpha_2$ establishes a $\mathfrak g[t]$-module isomorphism between $\cal F_{\lambda,\mu}$
and $\cal F_{-w_0\lambda,-w_0\mu}$. In the rest of the paper, we shall therefore study the modules $\cal F_{\lambda,\mu}$ for partitions $(\lambda, \mu)$ of $\lambda+\mu$ of first and second kind only.

\subsection{} Given dominant integral weights, $\lambda=\lambda_1\omega_1+\lambda_2\omega_2$, $\mu=\mu_1\omega_1+\mu_2\omega_2$ such that $|\lambda|\geq |\mu|$ and $\lambda_1 \geq \mu_1$, set
$${\mathcal F}_{\lambda,\mu}^+ =\left\lbrace \begin{array}{ll} {\mathcal F}_{\lambda+\omega_2,\mu-\omega_2}& \text{if}\  \lambda_2\geq \mu_2>0\\
{\mathcal F}_{\lambda+\omega_1,\mu-\omega_1}& \text{if}\  \mu_1>0, \mu_2=0\\
{\mathcal F}_{\lambda+(\mu_2-\lambda_2)\omega_2,\mu-(\mu_2-\lambda_2)\omega_2}& \text{if}\ \lambda_2 < \mu_2 .
\end{array}\right.$$

\label{cal.F.+}

\begin{thm}\label{filtration} Let  $\lambda=\lambda_1\omega_1+\lambda_2\omega_2$, $\mu=\mu_1\omega_1+\mu_2\omega_2\in P^+$ with $|\lambda|\geq |\mu|$ and $\lambda_1\geq \mu_1$.
There exists a short exact sequence
$$0\rightarrow \ker \phi(\lambda,\mu)\rightarrow \mathcal F_{\lambda,\mu}\overset{\phi(\lambda,\mu)}\longrightarrow \mathcal F_{\lambda,\mu}^+ \rightarrow 0
$$ where the submodule, $\ker \phi(\lambda,\mu)$ admits a filtration whose successive quotients are:
\begin{itemize}
\item[i.]  $$\begin{array}{c}\bigoplus\limits_{a\in S^{\lambda,\mu}_j} 
	\tau_{|\mu|}^\ast(V(\lambda+w_0\mu+a\alpha_2+(j-a)\alpha_1), \qquad  0\leq j< \mu_1+\mu_2, \\ \tau_{|\mu|}^\ast \mathcal F_{\lambda+\mu_2(\omega_1-\omega_2),\mu_1\omega_1},\end{array}$$
with $S^{\lambda,\mu}_j=\{a\in \mathbb Z:0\leq a\leq \mu_1,\, 0\leq j-a<\mu_2,\mu_2-\lambda_1\leq j-2a\leq\lambda_2-\mu_1\}$,
when $(\lambda,\mu)$ is a partition of first kind and $\mu_2\neq 0$.

\item[ii.]  $$\tau_{\mu_1}^\ast V(\lambda+w_0\mu +j\alpha_2), \qquad \max\{0, \mu_1-\lambda_2\} \leq j \leq \mu_1,\quad$$ when $(\lambda,\mu)$ is a partition of $\lambda+\mu$ of first kind with $\mu_2=0 $.

\item[iii.] $$\begin{array}{c}\bigoplus\limits_{a\in S_{\lambda,\mu}^{(\ell,k)}}
	\tau^\ast_{\mu_1+\lambda_2+\ell} V(\lambda+w_0\mu+(\mu_2-\lambda_2-\ell)\theta+a\alpha_2+(k-a)\alpha_1), \qquad \begin{array}{l} 1 \leq \ell \leq {\mu_2-\lambda_2},\\
		0 \leq k \leq \mu_1+\lambda_2\end{array},
\end{array}$$
with $S_{\lambda,\mu}^{(\ell,k)}=\{a\in \mathbb Z: 0\leq a\leq \mu_1,\, 0\leq k-a\leq \lambda_2,\, \mu_1-\mu_2+\ell\leq 2a-k\leq \lambda_1-\lambda_2-\ell\}$,
when $(\lambda,\mu)$ is a partition of $\lambda+\mu$ of second kind.
\end{itemize}
\end{thm}

Given  $(\lambda,\mu)\in P^+(\lambda+\mu,2)$, such that $\lambda(h_{\theta})=|\lambda| \geq |\mu|=\mu(h_{\theta})$ and $\lambda(h_{\alpha_1})\geq \mu(h_{\alpha_1})$,  
$$(\lambda,\mu) \succeq \left\{\begin{array}{ll} (\lambda+\omega_2,\mu-\omega_2) & \text{when}\, \lambda_2\geq\mu_2>0, \\
(\lambda+\omega_1,\mu-\omega_1) & \text{when}\, \lambda_1\geq\mu_1>0, \mu_2=0, \\ (\lambda+(\mu_2-\lambda_2)\omega_2,\mu-(\mu_2-\lambda_2)\omega_2) & \text{when}\, \lambda_2<\mu_2.
\end{array}\right.$$ Hence if $v_{\lambda,\mu}$ and $v_{\lambda,\mu}^+$ generate the CV-modules $\cal F_{\lambda,\mu}$ and $\cal F_{\lambda,\mu}^+$ respectively, then 
by \lemref{quotient.lem}(2) we know that there exists 
is a surjective morphism $\phi(\lambda,\mu): \cal F_{\lambda,\mu}\rightarrow \cal F_{\lambda,\mu}^+$ such that $\phi(\lambda,\mu)(v_{\lambda,\mu})=v_{\lambda,\mu}^+.$ 
In order to prove the theorem, we thus need to analyse the structure of $\ker\phi(\lambda,\mu)$. We do so in the next two sections. 
\subsection{} The following result determines the generators of the submodule $\ker \phi(\lambda,\mu)$.
\begin{prop} \label{surjective.l.m} For $\lambda,\mu\in P^+$, let $\cal F_{\lambda,\mu}, \cal F_{\lambda,\mu}^+$ be graded $\lie g[t]$-modules as described in \thmref{filtration}. Let $v_{\lambda,\mu}$ and $v_{\lambda,\mu}^+$ be the generators of $\cal F_{\lambda,\mu}$ and $\cal F_{\lambda,\mu}^+$ respectively and $\phi(\lambda,\mu): \cal F_{\lambda,\mu}\rightarrow \cal F_{\lambda,\mu}^+$  a surjective $\mathfrak g[t]$-module morphism such that $\phi(\lambda,\mu)(v_{\lambda,\mu})= v^+_{\lambda,\mu}$. Then  kernel $\phi(\lambda,\mu)$ is the $\mathfrak g[t]$-module generated by the set $\cal K_{\lambda,\mu}$ which is described as follows: 
$$\cal K_{\lambda,\mu} =\left\{ \begin{array}{ll}\{(x_{\alpha_2}^-\otimes t)^{\mu_2}v_{\lambda,\mu}, (x_{\theta}^-\otimes t)^{\mu_1+\mu_2}v_{\lambda,\mu}\}, & \text{when}\ \lambda_2\geq \mu_2>0,  \\\{(x_{\alpha_1}^-\otimes t)^{\mu_1}v_{\lambda,\mu}, 
(x_{\theta}^-\otimes t)^{\mu_1}v_{\lambda,\mu}\}, & \text{when}\ \mu_2=0\, 
\text{and}\, \, \mu_1 >0,\\ 
\{(x_{\theta}^-\otimes t)^{\mu_1+\lambda_2+s}v_{\lambda,\mu}:0<s\leq \mu_2-\lambda_2 \}, & \text{when}\ \lambda_2< \mu_2
\end{array}\right.$$
\end{prop} 
\proof  
\noindent \textbf{Case 1.} Suppose $0<\mu_2<\lambda_2$. By definition(\ref{cal.F.+}), $\cal F_{\lambda,\mu}^+ = \cal F_{\lambda+\omega_2,\mu-\omega_2}$.  Hence, $$\begin{array}{ll}
x_{\alpha_1}^-(r,s).v_{\lambda,\mu}^+=0, & \forall\ r,s\in \mathbb N,\, \text{with}\,\, r+s\geq 1+r+\mu_1\\
x_{\alpha_2}^-(r,s).v_{\lambda,\mu}^+=0, & \forall\ r,s\in \mathbb N,\, \text{with}\,\, r+s\geq 1+r+\mu_2-1\\
x_{\theta}^-(r,s).v_{\lambda,\mu}^+=0, & \forall\ r,s\in \mathbb N,\, \text{with}\,\, r+s\geq 1+r+\mu_1+\mu_2-1
\end{array}$$
whereas,
$$\begin{array}{ll}
x_{\alpha_i}^-(r,s).v_{\lambda,\mu}=0, & \forall\ r,s\in \mathbb N,\, \text{with}\,\, r+s\geq 1+r+\mu_i, \, \text{for}\, i=1,2\\
x_{\theta}^-(r,s).v_{\lambda,\mu}=0, & \forall\ r,s\in \mathbb N,\, \text{with}\,\, r+s\geq 1+r+\mu_1+\mu_2
\end{array}$$
Since $x_\alpha^-\otimes t^2.v_{\lambda, \mu} =0$ for all $\alpha\in R^+$, it follows that in this case, $$x_{\alpha_2}^-(r,s).v_{\lambda,\mu},\, x_{\theta}^-(r',s').v_{\lambda,\mu} \in \ker \phi(\lambda,\mu)$$ for $s=\mu_2$ and $s'=\mu_1+\mu_2$, i.e., 
$(x_{\alpha_2}^-\otimes t)^{\mu_2}v_{\lambda,\mu},\, (x_{\theta}^-\otimes t)^{\mu_1+\mu_2}.v_{\lambda,\mu}\in \ker\phi(\lambda, \mu).$

\noindent \textbf{Case 2.} Suppose $\mu_2=0$ and $\lambda_1\geq \mu_1>0$. By definition, we have $\mathcal{F}_{\lambda,\mu}^+= \mathcal{F}_{\lambda+\omega_1,\mu-\omega_1}$. Following similar arguments, it is easy to see that $$(x_{\alpha_1}^-\otimes t)^{\mu_1}v_{\lambda,\mu}, \, (x_{\theta}^-\otimes t)^{\mu_1}.v_{\lambda,\mu}\in \ker \phi(\lambda, \mu).$$

\noindent \textbf{Case 3.} Suppose $\mu_2>\lambda_2\geq 0$. By definition, we have $\cal F_{\lambda,\mu}^+ = \cal F_{\lambda+(\mu_2-\lambda_2)\omega_2,\mu-(\mu_2-\lambda_2)\omega_2}$, considering the following defining relations of $\cal F_{\lambda,\mu}$ and $\cal F_{\lambda,\mu}^+$, we get,
$$\begin{array}{ll}
x_{\alpha_i}^-(r,s).v_{\lambda,\mu}^+=0, & \forall\ r,s\in \mathbb N,\, \text{with}\,\, r+s\geq 1+r+\mu_i,\,  \, \text{for}\, \, i=1,2,\\
x_{\theta}^-(r,s).v_{\lambda,\mu}^+=0, & \forall\ r,s\in \mathbb N,\, \text{with}\,\, r+s\geq 1+r+\mu_1+\lambda_2
\end{array}$$
whereas,
$$\begin{array}{ll}
x_{\alpha_i}^-(r,s).v_{\lambda,\mu}=0, & \forall\ r,s\in \mathbb N,\, \text{with}\,\, r+s\geq 1+r+\mu_i,  \, \text{for}\, i=1,2\\
x_{\theta}^-(r,s).v_{\lambda,\mu}=0, & \forall\ r,s\in \mathbb N,\, \text{with}\,\, r+s\geq 1+r+\mu_1+\mu_2
\end{array}$$
As, $x_\alpha^-\otimes t^2.v_{\lambda,\mu}=0$, it follows that $ x_{\theta}^-(r,s).v_{\lambda,\mu} \in \ker \phi(\lambda,\mu)$ for $\mu_1+\lambda_2< s\leq \mu_2+\mu_1$ that is $$\{(x_{\theta}^-\otimes t)^{\mu_1+\lambda_2+\ell} .v_{\lambda,\mu}: 0< \ell\leq \mu_2-\lambda_2 \}\subset \ker \phi(\lambda, \mu).$$

To complete the proof of the proposition, we consider $X\in \bold{U}(\mathfrak g[t])$ is such that $X.v_{\lambda,\mu}^+=0$. Then $X$ can be written as $X=Y+Z$ where $Y$ is in the left ideal of $\bold{U}(\mathfrak g[t])$ generated by the set 
$$ I_{Y}(\lambda,\mu)= \left\lbrace x_{\alpha}^+ \otimes t^q, \, (x_{\alpha}^- \otimes 1)^{(\lambda+\mu)(h_{\alpha})+1}, (h\otimes t^q)-\delta_{q,0}(\lambda+\mu)(h).1 : q \in \Z^+, \alpha\in R^+,\, h\in \mathfrak h \right\rbrace, $$ and 
$ Z$ is in the left ideal of  $\bold{U}(\mathfrak{g}[t])$ generated by the set ${I}_{Z}^+(\lambda,\mu)$ where,\\ 
i. For $0<\mu_2 \leq \lambda_2$,  $${I}_{Z}^+(\lambda,\mu):= \begin{array}{l} \{x_{\alpha_1}^-(r,s): r,s\in \mathbb N, \, r+s\geq 1+r+\mu_1\} \\  \cup\ 
\left\lbrace x_{\alpha_2}^-(r,s),x_{\theta}^-(r',s'):  r,s,r',s'\in \mathbb N, \, r+s\geq r+\mu_2,\, \, r'+s'\geq r'+\mu_1+\mu_2
\right\rbrace. \end{array}$$
ii. For $0<\mu_1\leq \lambda_1,\,  \mu_2=0$,  $${I}_{Z}^+(\lambda, \mu):=\{x_{\alpha_1}^-(r, s), x_{\theta}^-(r', s'): r, s, r', s'\in \mathbb N, \, r+s\geq r+\mu_1, \, r'+s'\geq r'+\mu_1\}.$$
iii. For $\mu_2>\lambda_2$,  $${I}_{Z}^+(\lambda, \mu):= \begin{array}{l} \left\lbrace x_{\alpha_1}^-(r, s),x_{\alpha_2}^-(r', s'): r,s,r',s'\in \mathbb N, \, r+s\geq 1+r+\mu_1, r'+s'\geq 1+r'+\lambda_2 \right\rbrace \\ \cup \{x_{\theta}^-(r,s): r,s\in \mathbb N, \, r+s\geq 1+r+\mu_1+\lambda_2\}. \end{array}$$
Since $Y.v_{\lambda,\mu} =0,$ for all $Y \in I_Y(\lambda,\mu),$ and $Z.v_{\lambda,\mu} =0,$ for
$$ Z\in \left\{\begin{array}{l} \{x_{\alpha_1}^-(r,s): r,s\in \mathbb N, \, r+s\geq 1+r+\mu_1\}, \, \text{when}\, \,  0\leq \mu_2\leq \lambda_2\\
\{x_{\alpha_1}^-(r,s),x_{\alpha_2}^-(r', s'): r,s,r',s'\in \mathbb N, \, r+s\geq 1+r+\mu_1, \, r'+s'>1+r'+\lambda_2\} \,  \text{when}\, \, \mu_2>\lambda_2,\end{array}\right.$$  it is clear that $\ker \phi(\lambda,\mu)$ is generated by elements of $\mathcal{K}_{\lambda,\mu}$. Hence the proposition.
\subsection{} The following lemma and its corollary can be easily checked.   
\begin{lem} Given $a, b, p\in \mathbb Z_+$, we have,
\begin{enumerate}
\item[i.] $[x_{\alpha_1}^+\otimes t^a, (x_{\theta}^-\otimes t^b)^{(p)}] = -(x_{\theta}^-\otimes t^b)^{(p-1)} (x_{\alpha_2}^-\otimes t^{a+b})$,
\item[ii.]  $[x_{\alpha_2}^+\otimes t^a, (x_{\theta}^-\otimes t^b)^{(p)}] = (x_{\theta}^-\otimes t^b)^{(p-1)} (x_{\alpha_1}^-\otimes t^{a+b})$,
\item[iii.] $[x_{\alpha_1}^-\otimes t^a, (x_{\theta}^+\otimes t^b)^{(p)}] = (x_{\theta}^+\otimes t^b)^{(p-1)} (x_{\alpha_2}^+\otimes t^{a+b}$),
\item[iv.]  $[x_{\alpha_2}^-\otimes t^a, (x_{\theta}^+\otimes t^b)^{(p)}] = -(x_{\theta}^+\otimes t^b)^{(p-1)} (x_{\alpha_1}^+\otimes t^{a+b})$.
\end{enumerate}\label{e.f}
\end{lem}


\subsection{} In the rest of the paper, given $s,r\in \mathbb N$ and  $\alpha\in R^+,$ we set $\bold X_{\alpha}(r,s) = (x_{\alpha}^+\otimes t)^{(s)}(x_{\alpha}^-\otimes 1)^{(r+s)}.$
\vspace{-.3cm}
\begin{lem}\label{X.r.s} Suppose $V$ is a $\mathfrak{sl}_3[t]$-module and 
$v\in V$ is a non-zero weight vector satisfying,
$$\mathfrak n^+\otimes \mathbb C[t].v=0, \quad \, \mathfrak h\otimes t
\mathbb C[t].v=0, \quad \, x_{\alpha}^-\otimes t^2.v=0. $$
If for each $\alpha\in R^+$, there exists $l_\alpha, l_\alpha^1 \in \mathbb N$ 
such that 
$$(x_\alpha^-\otimes 1)^r.v =0, \, \forall\, r\geq l_\alpha, \qquad \bold X_\alpha(r,s)v=0 \quad \text{for } r,s\in \mathbb N,\, \ s\geq  l_\alpha^1, 
\hspace{.45cm} $$ then \\
i. for all $\alpha,\beta \in R^+$, such that  $\langle\alpha,\beta\rangle\geq 0$, 
\begin{eqnarray} \bold X_{\beta}(r,s) (x_\alpha^-\otimes t)^d v=0, & \quad \forall\ s\geq l_\alpha^1 -\langle \alpha,\beta\rangle d, \,   \label{N.alpha.1}
\end{eqnarray}  

\noindent ii.  for all $\alpha,\beta \in R^+$, such that  $\langle\alpha,\beta\rangle < 0$, 
\begin{equation}\label{N.alpha.3.1}
\bold X_\beta(r,s)(x_{\alpha}^-\otimes t)^d.v =0 \quad \text{for } s\geq l_{\beta}^{1}, 
\qquad (x_{\beta}^-\otimes 1)^r(x_{\alpha}^-\otimes t)^d.v =0 \quad \text{for }\ l\geq l_{\beta}+d
\end{equation}  \end{lem}
\proof \begin{itemize} 
\item[i.] Suppose $\alpha,\beta\in R^+$ such that $\langle \alpha,\beta\rangle \geq 0$.

\noindent Since $\langle \alpha,\beta\rangle=2$ if and only if $\alpha=\beta$, \eqref{N.alpha.1} follows from \cite[Corollary 6.6]{cv15} when $\langle \alpha,\beta\rangle =2$.\\
Now suppose that $\langle \alpha,\beta\rangle =1$. Without loss of generality, assume that $\beta-\alpha\in R^+.$ If $N_{-\alpha,\beta}, N_{\beta-\alpha,-\beta}\in \mathbb Z$ are such that
$[x^-_\alpha,x^+_\beta] =N_{-\alpha,\beta}x^+_{\beta-\alpha}$ and $[x^+_{\beta-\alpha},x^-_\beta] =N_{\beta-\alpha,-\beta}x^-_\alpha,$ then  
$$x_{\alpha}^-.\bold X_\beta(r,s).v = \bold X_\beta(r,s)x_{\alpha}^-.v + N_{-\alpha,\beta}N_{\beta-\alpha,-\beta}\bold X_\beta(r,s-1) (x_{\alpha}^-\otimes t).v.$$
Using \eqref{G.eqn}, and the relations satisfied by $v$, we see that for $\beta\in R^+$ and $r,s\in \N$,
$$\bold X_\beta(r,s)v = (-1)^{s}(x_\beta^-(r,s)h_\beta.v +\sum_{k\geq 1} x_\beta^-(r,s-k)P_\beta(u)_k).v = (-1)^s x_\beta^-(r,s)h_\beta.v.$$ Further using induction on 
$k$ and the condition that $x^-\otimes t^2.v=0$ it is easy to see that $P_\beta(u)_{k+1}.(x_{\alpha}^-)^{k}.v=0$, for all $k\in \mathbb Z_+$. Hence, 
$$\begin{array}{rl} x_{\alpha}^-.\bold X_\beta(r,s).v =& (-1)^s\left(x_{\beta}^-(r,s)h_{\beta}x_{\alpha}^-.v+x_{\beta}^-(r,s-1)(h_{\beta}\otimes t)x_{\alpha}^-.v \right)+ N_{-\alpha,\beta}N_{\beta-\alpha,-\beta}\bold X_\beta(r,s-1).x_{\alpha}^-\otimes t.v .\\
=&(-1)^s\left(x_{\beta}^-(r,s)x_{\alpha}^-.v+x_{\beta}^-(r,s-1)x_{\alpha}^-\otimes t.v \right)+ \bold X_\beta(r,s-1).x_{\alpha}^-\otimes t.v .\\ =& (-1)^s\left(x_{\alpha}^-.x_{\beta}^-(r,s).v+x_{\beta}^-(r,s-1)x_{\alpha}^-\otimes t.v \right)+ \bold X_\beta(r,s-1).x_{\alpha}^-\otimes t.v .\end{array}$$
Since for $\langle \alpha,\beta\rangle >0$, $\alpha+\beta\notin R^+$, $[x_\alpha^-,x_\beta^-\otimes t^s]=0$. Thus,  using the relations satisfied by $v$,  and \lemref{cv.2.5}, for $s\geq l_{\beta}^1$ we have,
\begin{equation} \label{G.eqn.12} \bold X_\beta(r,s-1).x_{\alpha}^-\otimes t.v = (-1)^{s-1}x_{\beta}^-(r,s-1).x_{\alpha}^-\otimes t.v.\end{equation}
On the other hand, given the condition, $(x_{\beta}^-\otimes t)^{(s)}.v=0,$ for $s\geq l_\beta^1$, we see that
$$x_{\beta-\alpha}^+(x_{\beta}^-)^{(r-s)}(x_\beta\otimes t)^{(s)}=0 \qquad \text{for } \ s\geq l_{\beta}^1,$$ 
$$ \Rightarrow \quad N_{\beta-\alpha,-\beta}[(x_{\beta}^-)^{(r-s)}(x_{\beta}^-\otimes t)^{(s-1)}(x_{\alpha}^-\otimes t).v + (x_{\alpha}^-).(x_{\beta}^-)^{(r-s-1)}(x_{\beta}^-\otimes t)^{(s)}].v = 0, \quad \text{for } \ s\geq l_{\beta}^1, \label{e_1.f_12}
$$
$$\Rightarrow \quad N_{\beta-\alpha,-\beta}(x_{\beta}^-)^{(r-s)}(x_{\beta}^-\otimes t)^{(s-1)}(x_{\alpha}^-\otimes t).v .v = 0, \quad \text{for } \ s\geq l_{\beta}^1. $$
Consequently, it follows from \eqref{G.eqn.12} that 
$$ \bold X_\beta(r,s-1).x_{\alpha}^-\otimes t.v = (-1)^{s-1} x_{\beta}^-(r,s-1).x_{\alpha}^-\otimes t.v =0, \quad \text{for } \ s\geq l_{\beta}^1. $$
That is, $ \qquad \qquad \bold X_\beta(r,s).x_{\alpha}^-\otimes t.v = 0, \hspace{.25cm} \text{whenever} \ s\geq l_{\beta}^1-1.\qquad $\\ 
Since for $\alpha,\gamma\in R^+$, $\alpha\neq \gamma$, $x_\alpha^+.(x_{\gamma}^-\otimes t)^q.v =0$, for all $q>0$ and under the given conditions, $h\otimes t^s.(x_{\alpha_2}^-\otimes t)^q.v =0, \ s>0,$ by repeating the above arguments we see that \eqref{N.alpha.1} holds.

\item[ii.] For $\mathfrak g$ of type $A_2$, given $\alpha,\beta\in R^+$ such that $\langle \alpha,\beta \rangle<0$, $\alpha+\beta\in R^+$ and $[x_\alpha^+,x_\beta^+]= -x_{\alpha+\beta}^+$,  $[x_\alpha^-,x_\beta^-]= x_{\alpha+\beta}^-$. Thus, we have
$$\begin{array}{ll} x_{\alpha}^-\otimes t. \bold X_\beta(r,s).v 
&= \bold X_\beta(r,s) x_{\alpha}^-\otimes t + \bold X_\beta(r-1,s).x_{\alpha+\beta}^-\otimes t.v 
\end{array}$$
Now 
using that $[x_{\alpha+\beta}^-\otimes t, x_{\beta}^-\otimes 1]=0$, and  the relation $x_\gamma^-\otimes t^2.v=0$ for all $\gamma\in R^+$, the above relation reduces to
\begin{equation}\begin{array}{ll} x_{\alpha}^-\otimes t. \bold X_\beta(r,s).v 
&= \bold X_\beta(r,s). x_{\alpha}^-\otimes t + x_{\alpha+\beta}^-\otimes t.\bold X_\beta(r-1,s).v.
\end{array}\label{2b+c}\end{equation}
Therefore under the given conditions we have $\quad \bold X_\beta(r,s) x_{\alpha}^-\otimes t.v=0,\hspace{.25cm} \text{for }\ 
s\geq l_{\beta}^1.$\\
Further, applying $x_{\alpha}^-\otimes t$ to the relation $(x_{\beta}^-\otimes 1)^s.v=0$ for $s\geq l_{\beta}^1$, we get
$$ 0=(x_\alpha^-\otimes t)(x_\beta\otimes 1)^{(s)}.v= 
[(x_\beta\otimes 1)^{(s)}.(x_\alpha^-\otimes t)+ (x_{\alpha+\beta}^-\otimes t)(x_\beta\otimes 1)^{(s-1)}].v \quad \text{for } s>l_\beta, $$
$$\Rightarrow (x_{\beta}^-\otimes 1)^s.x_{\alpha}^-\otimes t.v =0, \quad \text{for }\ 
s\geq l_\beta+1. $$ As above, by repeating the above arguments, we see that \eqref{N.alpha.3.1} holds. \endproof\end{itemize}

\begin{cor} Let $V$ be a $\mathfrak{sl}_3[t]$-module and $v\in V$ a non-zero weight vector satisfying the conditions given in \lemref{X.r.s}. If $\alpha,\beta\in R^+$ such that $\alpha+\beta\in R^+$, then,
\begin{enumerate}
\item[i.] \begin{equation}\begin{array}{l} (x_{\beta}^+)(x_{\beta}^-\otimes t)^{(a)}(x_{\alpha+\beta}^-\otimes t)^{(b)} (x_{\alpha}^-\otimes t)^{(a)}.v\\=-N_{\beta+\alpha,\beta}(x_{\beta}^-\otimes t)^{(a)}(x_{\alpha+\beta}^-\otimes t)^{(b-1)}(x_{\alpha}^-\otimes t)^{(c+1)}.v  \end{array}\label{e_i.V_j.1}. \end{equation} 

\item[ii.]  \begin{equation}\sum\limits_{d=0}^l {l \choose d} (x_\beta\otimes 1)^{a-d}
(x_{\alpha+\beta}^-\otimes t)^d(x_{\alpha}^-\otimes t)^{l-d}v =0, \hspace{.35cm}
\forall\ a\geq l_{\beta}. \label{} \end{equation}

\item[iii.] \begin{equation}\sum_{d=0}^l {l \choose d} \bold X_\beta(r-d,s)
(x_{\alpha+\beta}^-\otimes t)^d(x_{\alpha}^-\otimes t)^{l-d}v =0, \hspace{.35cm}
\forall\ s\geq l_{\beta}^1. \label{N.alpha.3.3} \end{equation}
\end{enumerate}
\end{cor}
\proof 
\begin{enumerate}\item[i.] Since $\mathfrak g$ is of type $A$, \eqref{e_i.V_j.1} is a straightforward consequence of the conditions on the vector $v$ and the fact that $[x_\alpha^+,x_{\alpha+\beta}^-]=N_{\alpha,\alpha+\beta}x_\beta^-$.

\item[ii.] Since $(x_\beta^-\otimes 1)^a.v=0$ for $a>l_\beta$, we get
$$\begin{array}{ll} 0&=(x_{\alpha}^-\otimes t)^l(x_\beta^-\otimes 1)^a.v\\
&=\sum\limits_{d=0}^l {l \choose d} (x_\beta\otimes 1)^{a-d}
(x_{\alpha+\beta}^-\otimes t)^d(x_{\alpha}^-\otimes t)^{l-d}v, \quad \text{for } a>l_\beta. \end{array}$$

\item[iii.] Applying $(x_{\alpha}^-\otimes t)^l$ on both sides of \eqref{2b+c} 
and using the relations satisfied by $v$, the fact that $[x_\alpha^-,x_\beta^-]=x_{\alpha+\beta}^-$ and \lemref{G.eqn} , we get 
$$\sum_{d=0}^l {l \choose d} \bold X_\beta(r-d,s)
(x_{\alpha+\beta}^-\otimes t)^d(x_{\alpha}^-\otimes t)^{l-d}v =0, \qquad \forall \  s\geq l_{\beta}^1. $$
This shows that \eqref{N.alpha.3.3} holds. \endproof
\end{enumerate} 
\section{Filtration of $\ker \phi(\lambda, \mu)$}\label{sec S_lambda,mu}
In this section, we study the structure of the kernel of the surjective homomorphism $\phi(\lambda, \mu): \cal F_{\lambda, \mu}\rightarrow \cal F_{\lambda,\mu}^+$. 
In \propref{Case.1} (respectively, \propref{Case.3} )  we consider the case when  $(\lambda, \mu)$ is a partition of $P^+(\lambda+\mu,2)$ of first (respectively second) kind.
\subsection{}\label{lem.case}
\begin{prop} Let $(\lambda,\mu)$ be a partition of first kind of $\lambda+\mu$ with $\mu=\mu_1\omega_1+\mu_2\omega_2$.  Then $\ker \phi(\lambda,\mu) $
admits a filtration $V_0\supset V_1\supset \cdots \supset V_{|\mu|}$ such that the following hold. 
\begin{enumerate}
\item[i.] When $\mu_2>0$, 
$$\phi^{(\lambda,\mu)}_{|\mu|}: \tau_{\mu_2}^\ast(\cal F_{\lambda+\mu_2(\omega_1-\omega_2),\mu_1\omega_1}) \rightarrow V_{|\mu|},$$ and 
$$\phi^{(\lambda,\mu)}_j: \bigoplus\limits_{a\in  S_j^{\lambda,\mu}} 
\tau_{|\mu|}^\ast(V(\lambda+w_0\mu+(j-a)\alpha_1+a\alpha_2) \rightarrow V_j/V_{j+1}, \quad \forall\ \,  0\leq j<|\mu|, $$  are surjective $\mathfrak{sl}_3[t]$-homomorphisms, with $S^{\lambda,\mu}_j=\{a\in \mathbb Z:0\leq a\leq \mu_1,\, 0\leq j-a<\mu_2,\mu_2-\lambda_1\leq j-2a\leq\lambda_2-\mu_1\}$. 
\item[ii.] When $\mu_2=0$, for each $ \max\{0,\mu_1-\lambda_2\}\leq j\leq \mu_1,$
$$\phi^{(\lambda,\mu)}_j: \tau_{\mu_1}^\ast V(\lambda+w_0\mu+j\alpha_2) \rightarrow V_j/V_{j+1}, $$ is a surjective $\mathfrak{sl}_3[t]$-homomorphisms.\end{enumerate}
\label{onto.filter.1}\end{prop} 
\vspace{-5mm}
\proof \label{lem.case1}  \begin{enumerate} \item[i]. 
Suppose that $\mu_2>0$. By \propref{surjective.l.m}, 
$$\ker\phi(\lambda,\mu) = \bu(\mathfrak g[t])(x_\theta^-\otimes t)^{|\mu|}v_{\lambda,\mu}+\bu(\mathfrak g[t])(x_{\alpha_2}^-\otimes t)^{\mu_2}v_{\lambda,\mu}.$$ Set, $$
V_{|\mu|} = \bu(\mathfrak g[t])(x_{\alpha_2}^-\otimes t)^{\mu_2}v_{\lambda,\mu},$$ 
$$\text{and } \quad V_j=\sum_{a_1+a_2=j}\bu(\mathfrak g[t])(x_{\alpha_1}^+\otimes 1)^{a_1}(x_{\alpha_2}^+\otimes 1)^{a_2}(x_{\theta}^-\otimes t)^{|\mu|}v_{\lambda,\mu} +V_{|\mu|}, \quad 0\leq j<|\mu|.$$
Clearly,  $\ker \phi(\lambda,\mu) = V_0$ and $V_{j+1}\subset V_j$ for $0\leq j\leq |\mu|.$
As 
$h\otimes t^sv_{\lambda,\mu}=0$ for all $h\in \lie h, s\in \mathbb N$, $\quad g\otimes t^2v_{\lambda,\mu}=0$ for all $g\in \mathfrak g$,  and 
$(x_{\alpha_i}^-\otimes t)^sv_{\lambda,\mu}=0$ for $s>\mu_i$, i=1,2,
$$(x_{\alpha_1}^+\otimes 1)^{a_1}(x_{\alpha_2}^+\otimes 1)^{a_2}(x_{\theta}^-\otimes t)^{|\mu|}v_{\lambda,\mu}=(-1)^{a_1} (x_{\alpha_1}^-\otimes t)^{a_2}(x_\theta^-\otimes t)^{|\mu|-a_1-a_2}(x_{\alpha_2}^-\otimes t)^{a_1}v_{\lambda,\mu}.$$ Consequently, $$V_j= \sum\limits_{j-a=0}^{\mu_2-1}\sum\limits_{a=0}^{\mu_1}\bu(\mathfrak g[t])(x_{\alpha_1}^-\otimes t)^{a}(x_\theta^-\otimes t)^{|\mu|-j}(x_{\alpha_2}^-\otimes t)^{j-a}v_{\lambda,\mu}+ V_{|\mu|}, \quad \text{for } 1\leq j<|\mu|. $$
Further, since $(x_\theta^-\otimes t)^{s}v_{\lambda,\mu} =0$ for $s\geq |\mu|+1$ implies  $$(x_{\alpha_1}^-\otimes t)^{a_1}(x_\theta^-\otimes t)^{s-a_1-a_2}(x_{\alpha_2}^-\otimes t)^{a_2}v_{\lambda,\mu} =0,\qquad \text{for } s>|\mu|,$$ 
we have  $$g\otimes t.\dfrac{V_j}{V_{j+1}} =0, \quad \forall \ g\in \mathfrak g, \ 0\leq j\leq |\mu|-1,$$ implying that for $0\leq j<|\mu|$, 
$\dfrac{V_j}{V_{j+1}}$ is a graded $\mathfrak g[t]$-module of the form $\tau^\ast_{|\mu|} W_j$ where $W_j$ is a finite-dimensional $\mathfrak g$-module. Hence there exists $\lambda_{kj}\in P^+$ such that 
$W_j$ is isomorphic to $\bigoplus\limits_{k} V(\lambda_{kj})$. 

\noindent Now observe that as $(x_{\alpha_1}^-\otimes 1)^{s_1}v_{\lambda,\mu}=0$ 
for $r=|\mu|-a_2$, $s=a_1$ and  $l= |\mu|-a_1$  in \eqref{N.alpha.3.3} such that $r+s = |\mu|-a_2+a_1 > \lambda_1+\mu_1 $, we have 
$$\sum_{d=0}^{|\mu|-a_1}
{|\mu|-a_1 \choose d} \bold X_1(|\mu|-a_2-d,a_1)
(x_{\theta}^-\otimes t)^d(x_{\alpha_2}^-\otimes t)^{|\mu|-a_1-d}v_{\lambda,\mu} =0.$$
Using the fact that $(x_{\alpha_2}^-\otimes t)^{|\mu|-a_1-d}v_{\lambda,\mu} =0$ for $|\mu|-a_1-d>\mu_2$, it thus follows that
$$ \sum_{i=0}^{\mu_2}
{|\mu|-a_1 \choose \mu_1-a_1+i} \bold X_1(\mu_2-a_2+a_1-i,a_1)
(x_{\theta}^-\otimes t)^{\mu_1-a_1+i}(x_{\alpha_2}^-\otimes t)^{\mu_2-i}v_{\lambda,\mu} 
=0, \quad  \text{whenever } \mu_2-a_2+a_1 > \lambda_1.$$ Consequently, 
for $\mu_2-a_2+a_1 > \lambda_1$,
$$\begin{array}{l} 
\sum_{i=0}^{\mu_2-a_2-1} 
{|\mu|-a_1 \choose \mu_1-a_1+i} \bold X_1(\mu_2-a_2+a_1-i,a_1)
(x_{\theta}^-\otimes t)^{\mu_1-a_1+i}(x_{\alpha_2}^-\otimes t)^{\mu_2-i}
v_{\lambda,\mu}\\ +{|\mu|-a_1 \choose \mu_1-a_1+\mu_2-a_2} \bold X_1(a_1,a_1)
(x_{\theta}^-\otimes t)^{\mu_1-a_1+\mu_2-a_2}(x_{\alpha_2}^-\otimes t)^{a_2}
v_{\lambda,\mu}\\+\sum_{i= \mu_2-a_2+1}^{\mu_2}
{|\mu|-a_1 \choose \mu_1-a_1+i} \bold X_1(\mu_2-a_2+a_1-i,a_1)
(x_{\theta}^-\otimes t)^{\mu_1-a_1+i}(x_{\alpha_2}^-\otimes t)^{\mu_2-i}v_{\lambda,\mu} 
=0.\end{array}$$ Using the relation $x_\alpha^-\otimes t^s.v_{\lambda,\mu}=0$ for $s\geq 2$, we deduce from above that 
$$\begin{array}{l} 
\sum_{i=0}^{\mu_2-a_2-1} 
{|\mu|-a_1 \choose \mu_1-a_1+i} \bold X_1(\mu_2-a_2+a_1-i,a_1)
(x_{\theta}^-\otimes t)^{\mu_1-a_1+i}(x_{\alpha_2}^-\otimes t)^{\mu_2-i}
v_{\lambda,\mu}\\ \quad  =- {|\mu|-a_1 \choose \mu_1-a_1+\mu_2-a_2} (x_{\alpha_1}^-\otimes t)^{(a_1)}(x_{\theta}^-\otimes t)^{\mu_1-a_1+\mu_2-a_2}(x_{\alpha_2}^-\otimes t)^{a_2}
v_{\lambda,\mu}.\end{array}$$
If $(a_1,a_2)$ is a pair of integers such that $a_1+a_2=j$, then using the facts that $$\bold X_1(\mu_2-a_2+a_1-i,a_1)
(x_{\theta}^-\otimes t)^{\mu_1-a_1+i}(x_{\alpha_2}^-\otimes t)^{\mu_2-i}
v_{\lambda,\mu}\in V_{\mu_2-i+a_1},$$ 
we see that, as  $a_2<\mu_2$, for $0\leq i\leq \mu_2-a_2-1$, $V_{\mu_2-i+a_1}\subset V_{j+1}$. Thus, 
$$(x_{\alpha_1}^-\otimes t)^{(a_1)}(x_{\theta}^-\otimes t)^{(|\mu|-j)}(x_{\alpha_2}^-\otimes t)^{(a_2)}v_{\lambda,\mu}\in V_{j+1}, \qquad \text{when } a_1+a_2=j,   \text{ and } |\mu|-a_2+a_1>\lambda_1+\mu_1.$$
\noindent Similarly, using the relation that 
$(x_{\alpha_2}^-\otimes t)^{(s_2)}v_{\lambda,\mu}=0$ for $s_2>\lambda_2+\mu_2$, it can be shown that 
$$(x_{\alpha_1}^-\otimes t)^{(a_1)}(x_{\theta}^-\otimes t)^{(|\mu|-j)}(x_{\alpha_2}^-\otimes t)^{(a_2)}v_{\lambda,\mu}\in V_{j+1}, \qquad \text{when } a_1+a_2=j,   \text{ and } |\mu|-a_1+a_2>\lambda_2+\mu_2.$$
Hence it follows that $\dfrac{V_j}{V_{j+1}}$ is a $\mathfrak g$-module generated by 
$$\{(x_{\alpha_1}^-\otimes t)^{(a)}(x_{\theta}^-\otimes t)^{(|\mu|-j)}(x_{\alpha_2}^-\otimes t)^{(j-a)}v_{\lambda,\mu}:0\leq a\leq \mu_1, 0\leq j-a <\mu_2, 
\mu_2-\lambda_1 \leq j-2a\leq \lambda_2-\mu_1\}.$$
$$\text{Set, } \qquad S^{\lambda,\mu}_j=\{a\in \mathbb Z: 0\leq a\leq \mu_1, 0\leq j-a < \mu_2,  
\mu_2-\lambda_1 \leq j-2a\leq \lambda_2-\mu_1\}.\qquad $$
Note that weight of  $(x_{\alpha_1}^-\otimes t)^{(a)}(x_{\theta}^-\otimes t)^{(|\mu|-j)}(x_{\alpha_2}^-\otimes t)^{(j-a)}v_{\lambda,\mu} = \lambda+w_0\mu+(j-a)\alpha_1+a\alpha_2,$ and 
$\lambda+w_0\mu+(j-a)\alpha_1+a\alpha_2\in P^+$, for all $a\in S^{\lambda,\mu}_j.$
Since by construction $(x_{\alpha_i}^+\otimes 1).(x_{\alpha_1}^-\otimes t)^{(a)}(x_{\theta}^-\otimes t)^{(|\mu|-j)}(x_{\alpha_2}^-\otimes t)^{(j-a)}v_{\lambda,\mu} \in V_{j+1}$ for $i=1,2$, and for a fixed $0\leq j<|\mu|$, the pairs $(a,j)$ are all distinct with $0\leq a\leq \min \{\mu_1,j\}$, it follows that the graded $\mathfrak g$-module $\dfrac{V_j}{V_{j+1}}$ is a quotient of $\bigoplus\limits_{a\in S^{\lambda,\mu}_j} \tau_{|\mu|}^\ast V(\lambda+w_0\mu+(j-a)\alpha_1+a\alpha_2)$, that is, for $0\leq j\leq |\mu|-1$,  there 
exists  a surjective $\mathfrak{sl}_3[t]$-homomorphism $\phi^{(\lambda,\mu)}_j:\bigoplus\limits_{a\in S^{\lambda,\mu}_j} \tau_{|\mu|}^\ast V(\lambda+w_0\mu+(j-a)\alpha_1+a\alpha_2)\rightarrow \dfrac{V_j}{V_{j+1}}.$ 

To see $V_{|\mu|}$ is a quotient of  $\tau_{\mu_2}^\ast \cal F_{\lambda+\mu_2(\omega_1-\omega_2),\mu_1\omega_1}$, observe that the generator $(x_{\alpha_2}^-\otimes t)^{\mu_2}v_{\lambda,\mu}$ of $V_{|\mu|}$ is $\mu_2$-graded element. Further using the relations satisfied by the vector $v_{\lambda,\mu}$ we see,
$$x_{\alpha_i}^+. (x_{\alpha_2}^-\otimes t)^{(\mu_2)}v_{\lambda,\mu} =0, \ \, \text{for}\ i=1,2\quad \quad
h\otimes t^{k+1}. (x_{\alpha_2}^-\otimes t)^{(\mu_2)}v_{\lambda,\mu} =0,  \ \, \forall\ k\geq 0, $$ implying that $V_{|\mu|}$ is a quotient of 
$W_{loc}(\lambda+\mu-\mu_2\alpha_2)$. Further, observing that $\langle \theta,\alpha_i\rangle =1$ for $i=1,2$, $\langle \alpha_1,\alpha_2\rangle =-1$ and using \lemref{X.r.s} with $l_{\gamma}=(\lambda+\mu)(h_\gamma)$, $l_{\gamma}^1=\mu(h_\gamma)$, for each $\gamma\in R^+$ it follows that 
$V_{|\mu|}$ is a quotient of the $CV$-module $\tau_{\mu_2}^\ast V(\bxi)$ where $\bxi=(\xi^\alpha)_{\alpha\in R^+}$ with
$$\xi^{\alpha_1}=(\lambda_1+\mu_2\geq \mu_1\geq 0), \hspace{.25cm} \xi^{\alpha_2}=(\lambda_2-\mu_2\geq  0), \hspace{.25cm}\xi^{\theta}=(\lambda_1+\lambda_2\geq \mu_1\geq 0). $$ Clearly, $V(\bxi)$ is isomorphic to $\cal F_{\lambda+\mu_2(\omega_1-\omega_2),\mu_1\omega_1}$.  \label{Case.1} 

\item[ii.] Suppose $\mu_2=0$ and $\mu_1>0$. Let $$K_1(\lambda,\mu_1) =\bu(\mathfrak g[t])(x_\theta^-\otimes t)^{\mu_1}v_{\lambda,\mu_1} .$$ In this case,  $e_2^ke_1^j (x_{\theta}^-\otimes t)v_{\lambda,\mu_1}=0, \quad \forall j>0,$ as $(x_{\alpha_2}^-\otimes t)v_{\lambda,\mu_1}=0$. Therefore 
setting $$V_j=\bu(\mathfrak g[t])e_2^j(x_{\theta}^-\otimes t)^{\mu_1}v_{\lambda,\mu_1} = \bu(\mathfrak g[t])(x_{\alpha_1}^-\otimes t)^{j}(x_{\theta}^-\otimes t)^{\mu_1-j}v_{\lambda,\mu_1} , \quad \text{for } 0\leq j\leq \mu_1,$$ we see that $V_{j+1}\subset V_{j}$ and $V_{\mu_1} = \bu(\mathfrak g[t])(x_{\alpha_1}^-\otimes t)^{\mu_1}v_{\lambda,\mu_1} $. Hence, by \propref{surjective.l.m}, $$K_1(\lambda,\mu)=\ker(\lambda,\mu).$$ Now using the same arguments as in part(i) of the proposition, we deduce that for each $0\leq j\leq \mu_1$, $\dfrac{V_{j}}{V_{j+1}}$ is a graded $\mathfrak g[t]$-module of the form $\tau^\ast W_j$ where $W_j$ is a finite-dimensional $\mathfrak g$-module generated by $(x_{\alpha_1}^-\otimes t)^{j}(x_{\theta}^-\otimes t)^{\mu_1-j}v_{\lambda,\mu_1}$. 

\noindent As $(x_{\alpha_2}^-\otimes 1)^{k}v_{\lambda,\mu}=0,$ for $k\geq \lambda_2+1,$ we see that 
\begin{equation}\sum_{d=0}^{\mu_1} {\mu_1 \choose d} (x_{\alpha_1}^-\otimes t)^{\mu_2-d}(x_{\theta}^-\otimes t)^d(x_{\alpha_2}^-\otimes 1)^{k-d}v =0, \qquad 
\forall\ k\geq \lambda_1+1. \label{case2.a} \end{equation} Thus 
when $\mu_1>\lambda_2$ and $k=\lambda_2+1$, \eqref{case2.a} reduces to 
$$
\sum_{d=0}^{\lambda_2+1} {\mu_1 \choose d} (x_{\alpha_1}^-\otimes t)^{\mu_1-d}
(x_{\theta}^-\otimes t)^d(x_{\alpha_2}^-\otimes 1)^{\lambda_2+1-d}v =0, $$ implying,
\begin{equation}  \label{mu_2.x_1.rel}
(x_{\theta}^-\otimes t)^{\lambda_1+1}(x_{\alpha_1}^-\otimes t)^{\mu_1-\lambda_2-1}v_{\lambda,\mu} = \sum_{d=0}^{\lambda_2} {\mu_1 \choose d} (x_{\alpha_1}^-\otimes t)^{\mu_1-d}
(x_{\theta}^-\otimes t)^d(x_{\alpha_2}^-\otimes 1)^{\lambda_2+1-d}v_{\lambda,\mu}. \end{equation}  
As the elements, $ (x_{\alpha_1}^-\otimes t)^{\mu_1-d}
(x_{\theta}^-\otimes t)^d(x_{\alpha_1}^-\otimes 1)^{\lambda_2+1-d}v_{\lambda,\mu}\in V_{\mu_1-d}$ and   $V_{\mu_1-d}\subseteq V_{\mu_1-\lambda_2}$ for $0\leq d\leq \lambda_2$, it follows from \eqref{mu_2.x_1.rel}, that $\dfrac{V_j}{V_{j+1}}$ is isomorphic to the 0-module for $0\leq d\leq \lambda_2$. Hence  $\ker\phi(\lambda,\mu_1)=V_{M}$ where  $M=\max \{0,\mu_1-\lambda_2\}$. As weight of $(x_{\alpha_1}^-\otimes t)^{j}(x_{\theta}^-\otimes t)^{\mu_1-j}v_{\lambda,\mu}$ is $\lambda-\mu_1\omega_2+j\alpha_2$ and by construction $x_{\alpha}^+.(x_{\alpha_1}^-\otimes t)^{j}(x_{\theta}^-\otimes t)^{\mu_1-j}v_{\lambda,\mu}\in V_{j+1}$, therefore $\dfrac{V_j}{V_{j+1}}$ is a quotient of $\tau_{\mu_1}^\ast V(\lambda-\mu_1\omega_2+j\alpha_2)$ for $M\leq j\leq \mu_1$. This completes the proof of the proposition. \endproof
\end{enumerate}
\subsection{} 
\begin{prop}\label{case.inv}
Let $(\lambda,\mu)$ be a partition of $\lambda+\mu$ of second kind and $|\rho^\mu_\lambda|=\sum\limits_{i=1}^2\min\{\lambda_i,\mu_i\}$. Then $\ker \phi(\lambda,\mu)$, has a filtration  
whose successive quotients are quotients of 
$$
\bigoplus\limits_{a\in S_{\lambda,\mu}^{(\ell,k)}}
\tau^\ast_{\mu_1+\lambda_2+\ell} V(\lambda+w_0\mu+(\mu_2-\lambda_2-\ell)\theta+a\alpha_2+(k-a)\alpha_1),\qquad \begin{array}{l}  1 \leq \ell \leq \lambda_2-\mu_2,\\ 0 \leq k \leq \mu_1+\lambda_2,\end{array}$$
where $S_{\lambda,\mu}^{(\ell,k)}=\{a\in \mathbb Z: 0\leq a\leq \mu_1,\, 0\leq k-a\leq \lambda_2,\, \mu_1-\mu_2+\ell\leq 2a-k\leq \lambda_1-\lambda_2-\ell\}$. 
\label{onto.filter}\end{prop}
\proof Let $$K_2(\lambda,\mu) = 
\bu(\mathfrak g[t])(x_{\theta}^-\otimes t)^{\mu_1+\lambda_2+1}v_{\lambda,\mu}.$$
Since $\bu(\mathfrak g[t])(x_{\theta}^-\otimes t)^{\mu_1+\lambda_2+\ell}v_{\lambda,\mu}\subseteq \bu(\mathfrak g[t])(x_{\theta}^-\otimes t)^{\mu_1+\lambda_2+1}v_{\lambda,\mu}$ for $\ell\geq 1$, $K_2(\lambda,\mu)$ is equal to the submodule of $\cal F_{\lambda,\mu}$ generated by the set of vectors $\{\bu(\mathfrak g[t])(x_{\theta}^-\otimes t)^{\mu_1+\lambda_2+s}v_{\lambda,\mu}:1\leq s\leq \mu_2-\lambda_2\}.$ Hence by \propref{surjective.l.m}, $K_2(\lambda,\mu)=\ker\phi(\lambda,\mu)$. 
For $1\leq \ell\leq \mu_2-\lambda_2$ and $0\leq k\leq |\rho^\mu_\lambda|$, let
$$V_{\ell,k} = \sum\limits_{\{(a_1,a_2): a_1+a_2=k, 0\leq a_i\leq \min\{\mu_i,\lambda_i\},i=1,2\}}\bu(\mathfrak g[t])(x_{\alpha_1}^+\otimes 1)^{(a_2)}(x_{\alpha_2}^+\otimes 1)^{(a_1)}(x_{\theta}^-\otimes t)^{\mu_1+\lambda_2+\ell}v_{\lambda,\mu}.$$ Notice that for all $1\leq \ell\leq \mu_2-\lambda_2$ and  $0\leq k\leq |\rho^\mu_\lambda|-1$, 
$V_{\ell,k+1}\subset V_{\ell,k}$ and $V_{\ell+1,0}\subset V_{\ell,|\rho^\mu_\lambda|}$. Further,  
as $x_\gamma^-\otimes t^2.v_{\lambda,\mu}=0$ for all $\gamma\in R^+$,  $(x_{\theta}^-\otimes t)^s.v_{\lambda, \mu}=0$, for $s>|\mu|$ and $$(x_{\alpha_1}^+\otimes 1)^{(a_2)}(x_{\alpha_2}^+\otimes 1)^{(a_1)}(x_{\theta}^-\otimes t)^{\mu_1+\lambda_2+\ell}v_{\lambda,\mu} = (x_{\alpha_1}^-\otimes t)^{(a_1)}(x_{\theta}^-\otimes t)^{\mu_1+\lambda_2+\ell-a_1-a_2}(x_{\alpha_2}^-\otimes t)^{(a_2)}v_{\lambda,\mu}\in V_{\ell,a_1+a_2},$$   
it follows that for all $g\in \mathfrak g$,
$$g\otimes t.\dfrac{V_{\ell,k}}{V_{\ell,k+1}}=0,  \qquad  \forall \quad 0\leq k\leq |\rho^\mu_\lambda|-1, \, \, 1\leq \ell\leq \mu_2-\lambda_2,$$ 
$$ g\otimes t.\dfrac{V_{\ell,|\rho^\mu_\lambda|}}{V_{\ell+1,0}}=0, \qquad \forall  \, 1\leq \ell\leq \mu_2-\lambda_2-1, \qquad \text{and } \qquad g\otimes t.V_{\mu_2-\lambda_2,|\rho^\mu_\lambda|}=0.$$
This implies that for all $0\leq k\leq |\rho^\mu_\lambda|-1,$ $1\leq \ell\leq \mu_2-\lambda_2$, $\dfrac{V_{\ell,k}}{V_{\ell,k+1}}$ is isomorphic to a graded $\mathfrak g[t]$-module of the form $\tau_{\mu_1+\lambda_2+\ell}^\ast W_{\ell,k}$ where $W_{\ell,k}$ is a finite-dimensional $\mathfrak g$-module, for $1\leq \ell\leq \mu_2-\lambda_2-1$. Similarly,  $\dfrac{V_{\ell,|\rho^\mu_\lambda|}}{V_{\ell+1,0}}$ and $V_{\mu_2-\lambda_2,|\rho^\mu_\lambda|}$ are  isomorphic to graded $\mathfrak g[t]$-modules of the form $\tau_{\mu_1+\lambda_2+\ell}^\ast W_{\ell,|\rho^\mu_\lambda|}$ and  $\tau_{|\mu|}^\ast W_{\mu_2-\lambda_2,|\rho^\mu_\lambda|}$
respectively, where $W_{\ell,|\rho^\mu_\lambda|}$ is a finite-dimensional $\mathfrak g$-module for $1\leq \ell\leq \mu_2-\lambda_2$.

For $1\leq \ell\leq \mu_2-\lambda_2$, if $k=|\rho^\mu_\lambda|=\mu_1+\lambda_2$ then the only values of $a_1,a_2\in \mathbb{Z}_+$ satisfying the conditions   $a_1+a_2=|\rho^\mu_\lambda|$, $a_1\leq \mu_1$ and $a_2\leq \lambda_2$ is $a_1=\mu_1$ and $a_2=\lambda_2$. Thus, in this case weight of  $(x_{\alpha_1}^-\otimes t)^{(a_1)}(x_{\theta}^-\otimes t)^{\mu_1+\lambda_2+\ell-a_1-a_2}(x_{\alpha_2}^-\otimes t)^{(a_2)}v_{\lambda,\mu}$ is $\lambda+\mu-\ell\theta-\mu_1\alpha_1-\lambda_2\alpha_2\in P^+$ for $1\leq \ell\leq \mu_2-\lambda_2$ as we are considering the case when $\mu_2>\lambda_2$ and $|\lambda|\geq |\mu|$. Consequently, it follows that $V_{\mu_2-\lambda_2,|\rho^\mu_\lambda|}$ and $\dfrac{V_{\ell,|\rho^\mu_\lambda|}}{V_{\ell+1,0}}$,    are respectively quotients of $\tau_{|\mu|}^\ast V(\lambda+w_0\mu-\mu_1\alpha_1-\lambda_2\alpha_2)$ and $\tau_{\mu_1+\lambda_2+\ell}^\ast V(\lambda+\mu-\ell\theta-\mu_1\alpha_1-\lambda_2\alpha_2)$, for $1\leq \ell\leq \mu_2-\lambda_2-1$.

\noindent On the other hand, since $v_{\lambda,\mu}$ satisfies the conditions 
of \lemref{X.r.s},  we see that if $r=\mu_1+\lambda_2-a_2+\ell$, $s=a_1$ and  $l= \mu_1+\lambda_2+\ell-a_1$  in \eqref{N.alpha.3.3}, then for  $r+s = \mu_1+\lambda_2+\ell-a_2+a_1 > \lambda_1+\mu_1 $ we have,
$$\underset{d=0}{\overset{\mu_1+\lambda_2+\ell-a_1}{\sum}}{\mu_1+\lambda_2+\ell-a_1\choose d} \bold X_1(\mu_1+\lambda_2+\ell-a_2-d,a_1)(x_\theta^-\otimes t)^d(x_{\alpha_2}^-\otimes t)^{\mu_1+\lambda_2+\ell-a_1-d} = 0.$$ Similarly, 
if $r=\mu_1+\lambda_2-a_1+\ell$, $s=a_2$ and  $l= \mu_1+\lambda_2+\ell-a_2$,  then for  $r+s = \mu_1+\lambda_2+\ell-a_1+a_2 > \lambda_2+\mu_2 $ we have,
$$\underset{d=0}{\overset{\mu_1+\lambda_2+\ell-a_2}{\sum}}{\mu_1+\lambda_2+\ell-a_2\choose d} \bold X_2(\mu_1+\lambda_2+\ell-a_1-d,a_2)(x_\theta^-\otimes t)^d(x_{\alpha_1}^-\otimes t)^{\mu_1+\lambda_2+\ell-a_2-d} = 0.$$
Thus if $a_1,a_2\in \mathbb{Z}_+$ are such that $a_1+a_2=k$ with $0\leq k<|\rho^\mu_\lambda|$, then similar calculations as in \propref{Case.1}(i), show that when 
$\mu_1+\lambda_2+\ell-a_2+a_1 > \lambda_1+\mu_1$ or $\mu_1+\lambda_2+\ell-a_1+a_2 > \lambda_2+\mu_2$, then
$(x_{\alpha_1}^-\otimes t)^{(a_1)}(x_\theta^-\otimes t)^{\mu_1+\ell-a_1-a_2}(x_{\alpha_2}^-\otimes t)^{(a_2)}v_{\lambda,\mu}\in V_{\ell,k+1}$.
Hence, for $1\leq \ell\leq \mu_2-\lambda_2$ and $0\leq k\leq |\rho^\mu_\lambda|-1$,
$\dfrac{V_{\ell,k}}{V_{\ell,k+1}}$ is a quotient of the finite-dimensional graded $\mathfrak g[t]$-module 
$$\bigoplus\limits_{a\in S_{\lambda,\mu}^{(\ell,k)}}
\tau^\ast_{\mu_1+\lambda_2+\ell} V(\lambda+w_0\mu+(\mu_2-\lambda_2-\ell)\theta+a\alpha_2+(k-a)\alpha_1),$$
with $S_{\lambda,\mu}^{(\ell,k)}=\{a\in \mathbb Z: 0\leq a\leq \mu_1,\, 0\leq k-a\leq \lambda_2,\, \mu_1-\mu_2+\ell\leq 2a-k\leq \lambda_1-\lambda_2-\ell\}$.
This completes the proof of the proposition.\label{Case.3}
\endproof
\vspace*{-1cm}\section{Proof of \thmref{filtration}}\label{proof}
\subsection{} In this section, we complete the proof of \thmref{filtration}.
We have seen (\remref{f-lambda dim}),   
\begin{equation}\label{1}\dim \mathcal F_{\lambda,\mu} \geq \dim V(\lambda).\dim V(\mu),\quad \quad \forall \, \, (\lambda, \mu) \in P^+(\lambda+\mu,2).\end{equation} We now prove the reverse inequality using the results of \secref{sec S_lambda,mu}. Throughout this section, we assume 
$\dim \mathcal F_{\nu,\gamma}=0$, whenever either $\nu$ or $\gamma$ is non-zero, 
non-dominant weight and $\dim \mathcal F_{\nu,\gamma}=\dim V(\nu)$ whenever 
$\gamma=0$, and $\nu \in P^+$. 
\label{l*m.dim}
\begin{prop} Let $\lambda=\lambda_1\omega_1+\lambda_2\omega_2$, $\mu= \mu_1\omega_1+\mu_2\omega_2\in P^+$ and  $\rho^\mu_\lambda=\mu_1\omega_1+\lambda_2\omega_2$. 
\begin{enumerate}
	\item[i.] If $(\lambda,\mu)$ is a partition  of $\lambda+\mu$ of first kind with $\mu_1>0$ and $\mu_2=0$, then the surjective homomorphisms 
	$\{\phi_j^{(\lambda, \mu)}: \max\{0,\mu_1-\lambda_2\}\leq j \leq \mu_1\}$ given in \propref{Case.1}(ii) are isomorphisms. 
	\item[ii.] If $(\lambda,\mu)$ is a partition  of $\lambda+\mu$ of first kind with $\mu_2>0$ and $\mu_1>0$, then the surjective homomorphisms  
	$\{\phi_j^{(\lambda,\mu)}: 0\leq j\leq |\mu|\}$ given in \propref{Case.1}(i) are isomorphisms.
	\item[iii.] If $(\lambda,\mu)$ is a partition  of $\lambda+\mu$ of second kind, then $\ker \phi(\lambda,\mu)$ has a filtration where the successive quotients are isomorphic to $
	\bigoplus\limits_{a\in S_{\lambda,\mu}^{(\ell,k)}}
	\tau^\ast_{\mu_1+\lambda_2+\ell} V(\lambda+w_0\mu+(\mu_2-\lambda_2-\ell)\theta+a\alpha_2+(k-a)\alpha_1),$
	where $S_{\lambda,\mu}^{(\ell,k)}=\{a\in \mathbb Z: 0\leq a\leq \mu_1,\, 0\leq k-a\leq \lambda_2,\, \mu_1-\mu_2+\ell\leq 2a-k\leq \lambda_1-\lambda_2-\ell\}$.
	Consequently, $\dim \cal F_{\lambda,\mu} = \dim V(\lambda)\dim V(\mu).$
\end{enumerate} 
\end{prop}
\subsection{Proof of \propref{l*m.dim}(i).}\label{part1.dim} We prove part(i) by considering two subcases.
\item[Subcase 1.] Suppose $\lambda$ and $\mu$ are both multiples of $\omega_1$. Since $\lambda_1\geq \mu_1$, by \propref{Case.1}(ii), $\ker\phi(\lambda_1\omega_1,\mu_1\omega_1)$ is a quotient of $V((\lambda_1-\mu_1)\omega_1+\mu_1\omega_2)$, implying, \begin{equation}\label{mu.2=0=lambda2} \begin{array}{rl}\dim \cal F_{\lambda,\mu} =& \dim \cal F_{(\lambda_1+1)\omega_1,(\mu_1-1)\omega_1} +\dim \ker \phi(\lambda,\mu)\\ &\leq \dim \cal F_{(\lambda_1+1)\omega_1,(\mu_1-1)\omega_1} +\dim V((\lambda_1-\mu_1)\omega_1+\mu_1\omega_2).\end{array}\end{equation}
By assumption, for $\mu_1=1$, $\cal F_{(\lambda_1+1)\omega_1,(\mu_1-1)\omega_1}$ is isomorphic to $ V((\lambda_1+1)\omega_1)$. Hence we have, $$\begin{array}{rl}\dim \cal F_{\lambda,\mu}\leq &\dim V(\lambda_1+1)\omega_1) +\dim V((\lambda_1-1)\omega_1+\omega_2) = \frac{1}{2}[(\lambda_1+2)(\lambda_1+3)+\lambda_1(\lambda_1+2)]\\ \leq &\frac{3}{2}(\lambda_1+1)(\lambda_1+2)=  \dim V(\lambda_1\omega_1)\dim V(\omega_1).\end{array}$$ Along with \eqref{1}, this implies that, $\dim \cal F_{\lambda_1\omega_1,\omega_1}= \dim V(\lambda_1\omega_1)\dim V(\omega_1)$. 

\noindent Now, by induction hypothesis assume that $\dim \cal F_{\lambda_1\omega_1,\mu_1\omega_1}= \dim V(\lambda_1\omega_1)\dim V(\mu_1\omega_1)$ for all $\mu_1\in \mathbb N$ such that $\mu_1< n$. Then, using  \propref{Case.1}(ii), we have, 
$$\begin{array}{rl}4 \dim \cal F_{\lambda_1\omega_1,n\omega_1} \leq & 4 \dim \cal F_{(\lambda_1+1)\omega_1,(n-1)\omega_1} + 4 \dim V((\lambda_1-n)\omega_1+n\omega_2)\\ \leq & [(\lambda_1+2)(\lambda_1+3)n(n+1)+2(\lambda_1-n+1)(n+1)(\lambda_1+2)]\\ \leq & (\lambda_1+1)(\lambda_1+2)(n+1)(n+2)= 4 \dim V(\lambda_1\omega_1)\dim V(n\omega_1).\end{array}$$ 
Along with \eqref{1}, this implies, $\dim \cal F_{\lambda_1\omega_1,n\omega_1}=\dim V(\lambda_1\omega_1)\dim V(n\omega_1)$. Consequently, 
$$ \ker \phi(\lambda_1\omega_1,\mu_1\omega_1)\cong_{\mathfrak{sl}_3[t]} \tau_{\mu_1}^\ast V((\lambda_1-\mu_1)\omega_1+\mu_1\omega_2) = \tau_{\mu_1}^\ast V(\lambda+w_0\mu+\mu_1\alpha_2).$$

\item[Subase 2.]{Suppose $\lambda=\lambda_1\omega_1+\lambda_2\omega_2, \mu=\mu_1\omega_1 \in P^+$ and $\lambda_1>\mu_1>0$.} Then by \propref{Case.1}(ii), 
\begin{equation}\label{filter1}\dim \ker(\lambda,\mu_1\omega_1)\leq \sum\limits_{j=\max\{0,\mu_1-\lambda_2\}}^{\mu_1} \dim  V(\lambda-\omega_2+j(2\omega_2-\omega_1)).
\end{equation}
In particular, when $\mu_1=1$, 
$$	\begin{array}{rl}
2 \dim \cal F_{\lambda,\omega_1} \leq & 2\dim V(\lambda+\omega_1)  + 2 \dim V(\lambda+(\omega_2-\omega_1)) +  2\dim V(\lambda-\omega_2)\\
\leq &  (\lambda_1+2)(\lambda_2+1)(\lambda_1+\lambda_2+3)+ \lambda_1 (\lambda_2+2)(\lambda_1+\lambda_2+2)+ (\lambda_1+1)(\lambda_2)(\lambda_1+\lambda_2+1) \\
\leq & 3(\lambda_1+1)(\lambda_2+1)(\lambda_1+\lambda_2+2) + (\lambda_2+1)(2\lambda_1+\lambda_2+4) -(\lambda_2-\lambda_1+1)(\lambda_1+\lambda_2+2)\\ &-
(\lambda_1+1)(\lambda_1+2\lambda_2+2)\\
\leq &  3 (\lambda_1+1)(\lambda_2+1)(\lambda_1+\lambda_2+2)= 2 \dim V(\omega_1). \dim V(\lambda)
\end{array}$$
Along with \eqref{1}, this shows $\dim \cal F_{\lambda,\mu_1\omega_1}=\dim V(\lambda)\dim V(\omega_1)$ and 
$$\ker\phi(\lambda, \omega_1)\cong_{\mathfrak{sl}_3[t]} \tau_1^\ast V(\lambda-(\omega_1-\omega_2)) \oplus \tau_1^\ast V(\lambda-\omega_1) = \bigoplus\limits_{j=0}^1 \tau_1^\ast V(\lambda+w_0\mu+j\alpha_2), $$ as claimed in \thmref{filtration}.
Now by induction hypothesis assume that \propref{l*m.dim}(i) holds for all CV modules $ \cal F_{\lambda,\mu_1\omega_1}$ with $\mu_1< n$, i.e, 
$$\dim \cal F_{\lambda,\mu_1\omega_1}=\dim V(\lambda) \dim V(\mu_1\omega_1)\quad \text{whenever } \mu_1<n.$$  
Then  using \propref{Case.1}(ii) and applying induction hypothesis,  we see that
for $\mu_1=n$, 
$$\begin{array}{rl}
\dim \ker \phi(\lambda,n\omega_1) 
\leq & \dim V(\lambda+n(\omega_2-\omega_1)) + \sum\limits_{j=\max\{0,n-\lambda_2\}}^{n-1} \dim  V(\lambda-\omega_2+-(n-1)\omega_2+j(2\omega_2-\omega_1))\\
\leq & \dim V(\lambda+n(\omega_2-\omega_1))+ \dim \ker \phi(\lambda-\omega_2,(n-1)\omega_1).\end{array}$$
\noindent By \propref{Case.1}, 
$\dim  \ker \phi(\lambda,\mu_1\omega_1)= \dim \cal F_{\lambda,\mu_1\omega_1} -\dim \cal F_{\lambda+\omega_1, (\mu_1-1)\omega_1}$. Hence,  
\begin{equation} \dim \cal F_{\lambda,n\omega_1} \leq \dim \cal F_{\lambda+\omega_1, (n-1)\omega_1}  +  \dim \cal F_{\lambda-\omega_2,(n-1)\omega_1} -\dim \cal F_{\lambda+\omega_1-\omega_2, (n-2)\omega_1}+  \dim V(\lambda+n(\omega_2-\omega_1)) .\label{dim.mu2=00}
\end{equation} 
Consequently,
\begin{equation}\begin{array}{rl} 4 \dim \cal F_{\lambda,n\omega_1} \leq & 
(n)(n+1) [(\lambda_1+2)(\lambda_2+1)(\lambda_1+\lambda_2+3)+(\lambda_1+1)(\lambda_2)(\lambda_1+\lambda_2+1)]\\
&-(\lambda_1+2)(\lambda_2)(\lambda_1+\lambda_2+2)(n-1)n+ 2 (\lambda_1-n+1)(\lambda_2+n+1)(\lambda_1+\lambda_2+2)\\
\leq & (n)(n+1) [3(\lambda_1+1)(\lambda_2+1)(\lambda_1+\lambda_2+2)-\lambda_1(\lambda_2+2)(\lambda_1+\lambda_2+2)]\\
& - (\lambda_1+2)(\lambda_2)(\lambda_1+\lambda_2+2)(n-1)n + 
2 (\lambda_1-n+1)(\lambda_2+n+1)(\lambda_1+\lambda_2+2)\\
&\hspace{7cm}	(\text{using induction step for $\mu_1=1$}) 
\\ \leq & (n)(n+1)(\lambda_1+\lambda_2+2) [3(\lambda_1+1)(\lambda_2+1)-\lambda_1(\lambda_2+2)-\lambda_2(\lambda_1+2)]\\ & +2(\lambda_1+\lambda_2+2) [(\lambda_1+2)\lambda_2n+ (\lambda_1-n+1)(\lambda_2+n+1)] \\
\leq & (\lambda_1+\lambda_2+2)n(n+1)[(\lambda_1+1)(\lambda_2+1)+2]
\\ & +2(\lambda_1+\lambda_2+2) [n(\lambda_1\lambda_2 +2\lambda_2-\lambda_2 +\lambda_1+1)-n(n+1)+(\lambda_1+1)(\lambda_2+1)]\\
\leq & (\lambda_1+\lambda_2+2)n(n+1) [(\lambda_1+1)(\lambda_2+1)]
+2(\lambda_1+\lambda_2+2)[(n+1)(\lambda_1+1)(\lambda_2+1)]\\
\leq & (\lambda_1+\lambda_2+2)(\lambda_1+1)(\lambda_2+1)(n+1)(n+2)  =  4\dim V(\lambda). \dim V(\mu_1 \omega_1), \label{lambda.mu1}
\end{array}\end{equation}
Along with \eqref{1}, it now follows that 
$\dim \cal F_{\lambda,\mu_1\omega_1} =\dim V(\lambda)\dim V(\mu_1\omega_1)$, which further establishes \propref{l*m.dim}(i). Thus, we see that \thmref{filtration} holds when $(\lambda, \mu)$ is partition of $\lambda+\mu$ of first kind with $\mu_2=0$. 
\vspace*{-2cm}\subsection{Proof of \propref{l*m.dim}(ii)} Given $(\lambda,\mu)$ a partition of $\lambda+\mu$ of first kind with $\mu_1, \mu_2 > 0$, we prove part (ii) by applying induction on $\mu_2.$  Set, 
$\mathbb{S}^{\lambda,\mu} = \bigcup\limits_{j=0}^{|\mu|} \{(a_1,a_2): a_1\in S_j^{\lambda,\mu}, a_2=j-a_1\} $. By \propref{onto.filter.1}, 
\begin{equation}\label{dim.l.m} \begin{array}{rl}
\dim \cal F_{\lambda,\mu}= & \dim \cal F_{\lambda+\omega_2,\mu-\omega_2} +\, 
\dim \ker(\lambda,\mu)\\ 
\leq & \dim \cal F_{\lambda+\omega_2,\mu-\omega_2} + \dim \cal F_{\lambda+\mu_2(\omega_1-\omega_2), \mu_1\omega_1} + \underset{(a_1,a_2)\in {\mathbb S}^{\lambda,\mu}}{\sum} \dim V(\lambda+w_0\mu+a_2\alpha_1+a_1\alpha_2).\end{array}\end{equation}
In particular, for $\mu_2=1$, if $(a_1,a_2)\in {\mathbb S}^{\lambda,\mu_1\omega_1+\omega_2}$, then by definition of $S^{\lambda,\mu}_j$,  $$a_2=0, \qquad \mu_1-\lambda_2\leq a_1\leq \min\{\mu_1, \lambda_1-1\}.$$
Hence,
\begin{equation} \dim \cal F_{\lambda,\mu_1\omega_1+\omega_2}
\leq  \dim \cal F_{\lambda+\omega_2,\mu_1\omega_1} + \dim \cal F_{\lambda+(\omega_1-\omega_2), \mu_1\omega_1} + \underset{j=\max\{0,\mu_1-\lambda_2\}}{\overset{\min\{\mu_1, \lambda_1-1\}}{\sum}} \dim V(\lambda+w_0\mu+j\alpha_2).\label{dim.l.m.1}\end{equation}
On the other hand, using \propref{Case.1}(ii) and \propref{l*m.dim}(i), we have, 
$$\begin{array}{rl}\dim \ker(\lambda-\theta,\mu-\theta) = & \dim \ker(\lambda-\theta,(\mu_1-1)\omega_1) = \underset{j=\max\{0,\mu_1-\lambda_2\}}{\overset{\mu_1-1}{\sum}} \dim V(\lambda-\theta+w_0(\mu-\theta)+j\alpha_2)  \end{array}$$
Hence \eqref{dim.l.m.1} can be rewritten as,
\begin{equation}\begin{array}{rl}\dim \cal F_{\lambda,\mu_1\omega_1+\omega_2}
\leq & \dim \cal F_{\lambda+\omega_2,\mu_1\omega_1} + \dim \cal F_{\lambda+(\omega_1-\omega_2), \mu_1\omega_1} + \dim \ker(\lambda-\theta,(\mu_1-1)\omega_1) \\ &+ (1-\delta_{\mu_1,\lambda_1})\dim V(\lambda-\omega_1-\mu_1\omega_1+\mu_1\omega_2)\end{array}\label{dim.l.m.1.a}
\end{equation}
Denoting the right-hand side of the inequality \eqref{dim.l.m.1.a} by $R_{\lambda_1>\mu_1}$ in the case when $\lambda_1>\mu_1$, we have,
$$\begin{array}{l}
R_{\lambda_1>\mu_1} = 4\dim \cal F_{\lambda+\omega_2,\mu_1\omega_1}+ 4\dim \cal F_{\lambda+(\omega_1-\omega_2), \mu_1\omega_1}  + 4 \dim \cal F_{\lambda-\theta,(\mu_1-1)\omega_1} \\ \qquad ~ + 4 \dim V((\lambda_1-\mu_1-1)\omega_1+(\lambda_2+\mu_1)\omega_2) - 4\dim \cal F_{(\lambda_1+\mu_2-2)\omega_1+(\lambda_2-\mu_2)\omega_2+\omega_1,(\mu_1-2)\omega_1}\\
=(\mu_1+1)(\mu_1+2)[(\lambda_1+1)(\lambda_2+2)(\lambda_1+\lambda_2+3)(\mu_1+1)(\mu_1+2) + 
(\lambda_1+2)(\lambda_2)(\lambda_1+\lambda_2+2)] \\
\quad +\lambda_1\lambda_2(\lambda_1+\lambda_2)\mu_1(\mu_1+1)+(\lambda_1-\mu_1)(\lambda_2+\mu_1+1)(\lambda_1+\lambda_2+1)2-(\lambda_1+1)\lambda_2(\lambda_1+\lambda_2+1)\mu_1(\mu_1-1) \\ \hspace{12cm}  (\text{Using \propref{l*m.dim}.(i)})\\
=(\mu_1+1)(\mu_1+2)[3(\lambda_1+1)(\lambda_2+1)(\lambda_1+\lambda_2+2) - 
(\lambda_2+1)\lambda_1(\lambda_1+\lambda_2+1)] 
+\lambda_1\lambda_2(\lambda_1+\lambda_2)\mu_1(\mu_1+1)\\ \quad +2(\lambda_1-\mu_1)(\lambda_2+\mu_1+1)(\lambda_1+\lambda_2+1) -(\lambda_1+1)\lambda_2(\lambda_1+\lambda_2+1)\mu_1(\mu_1+1-2)
\end{array}$$
$$\begin{array}{rl}
R_{\lambda_1>\mu_1} =&3(\lambda_1+1)(\lambda_2+1)(\lambda_1+\lambda_2+2)(\mu_1+1)(\mu_1+2)\\ &- (\mu_1+1)\mu_1\{(\lambda_2+1)(\lambda_1)(\lambda_1+\lambda_2+1)-\lambda_1\lambda_2(\lambda_1+\lambda_2)+(\lambda_1+1)\lambda_2(\lambda_1+\lambda_2+1)\} \\
& +2(\lambda_1+\lambda_2+1)\{(\lambda_1-\mu_1)(\lambda_2+\mu_1+1)+(\lambda_1+1)\lambda_2\mu_1-(\lambda_2+1)\lambda_1(\mu_1+1)\}\\
=&3(\lambda_1+1)(\lambda_2+1)(\lambda_1+\lambda_2+2)(\mu_1+1)(\mu_1+2) - (\mu_1+1)\mu_1\{(\lambda_1+\lambda_2+\lambda_1\lambda_2)(\lambda_1+\lambda_2+1)+\lambda_1\lambda_2\} \\
&-2(\lambda_1+\lambda_2+1)\mu_1(\mu_1+1) \\
=&3(\lambda_1+1)(\lambda_2+1)(\lambda_1+\lambda_2+2)(\mu_1+1)(\mu_1+2) -\mu_1(\mu_1+1)(\lambda_1+1)(\lambda_2+1)(\lambda_1+\lambda_2+2)\\
=& 2(\lambda_1+1)(\lambda_2+1)(\lambda_1+\lambda_2+2)(\mu_1+1)(\mu_1+3) = 4 \dim V(\lambda) \dim V(\mu_1\omega_1+\omega_2)
\end{array}$$ 
Using similar computations, it can be verified that for $\mu_1=\lambda_1$, the right-hand side of \eqref{dim.l.m.1.a} is equal to $\dim V(\lambda) \dim V(\mu_1\omega_1+\omega_2)$. Thus, 
it follows that
$\dim \cal F_{\lambda,\mu_1\omega_1+\omega_2} = \dim V(\lambda)\dim V(\mu_1\omega_1+\omega_2),$ and \propref{l*m.dim}(ii) holds in this case.\\ 
By induction hypothesis, assume \propref{l*m.dim}(ii) holds for $\mu_2\in \mathbb N$ with  $\mu_2\leq k$. From the definition of the set $S_j^{\lambda,\mu}$ (ref. 
\propref{onto.filter.1}), it follows that, 
\begin{equation}\label{union.ninv} 
\begin{array}{rl}{\mathbb S}^{\lambda,\mu} =& {\mathbb S}^{\lambda-\theta,\mu-\theta} \cup \{(\mu_1,b):  b\in \Z_+, 
0\leq b \leq \mu_2-1, \, \mu_2+\mu_1-b\leq \lambda_1, b\leq \lambda_2 \} \\
&\cup \,  \{(a,\mu_2-1): a\in \Z_+, 0\leq a\leq \mu_1-1, \mu_2+a-\mu_2+1\leq \lambda_1, \, |\mu|-1-a\leq \lambda_2\}
\end{array}\end{equation} Hence using \propref{onto.filter.1} and applying induction hypothesis on $\mu_2$ we get: 
\begin{equation}\begin{array}{l}
\underset{(a_1,a_2)\in {\mathbb S}^{\lambda-\theta,\mu-\theta}}{\sum} \dim V(\lambda+w_0\mu+a_2\alpha_1+a_1\alpha_2) 
=  \dim \ker(\lambda-\theta,\mu-\theta) - \dim \cal F_{\lambda-\theta+(\mu_2-1)(\omega_1-\omega_2), (\mu_1-1)\omega_1}\\\hspace{1.5cm}  = 
\dim \cal F_{\lambda-\theta,\mu-\theta} - \dim \cal F_{\lambda-\theta+\omega_2,\mu-\theta-\omega_2} - \dim \cal F_{\lambda-\theta+(\mu_2-1)(\omega_1-\omega_2), (\mu_1-1)\omega_1}.
\end{array}\label{dim.l.m.2}\end{equation}
Further setting, $\lambda'=\lambda-\omega_1+\mu_1(\omega_2-\omega_1)$ and $\lambda''= \lambda-\theta+(\mu_2-1)(\omega_1-\omega_2)$, we have
\begin{equation}\begin{array}{l}
\underset{b=\max\{0,|\mu|-1-\lambda_2\}}{\overset{\mu_2-1}{\sum}}\dim V(\lambda+w_0\mu+b\alpha_1+\mu_1\alpha_2) = \sum\limits_{b=\max\{0,|\mu|-1-\lambda_2\}}^{\mu_2-1}
\dim V(\lambda'+w_0(\mu_2-1)\omega_2 +b\alpha_1)\\
\hspace{1.5cm} = \dim \ker(\lambda',(\mu_2-1)\omega_2)=
\dim \cal F_{\lambda',(\mu_2-1)\omega_2}- \dim \cal F_{\lambda'+\omega_2,(\mu_2-2)\omega_2},
\end{array}\label{mu.12'}\end{equation} and 
\begin{equation}\begin{array}{l}
\underset{a=\max\{0,|\mu|-1-\lambda_2\}}{\overset{\mu_1-1}{\sum}}\dim V(\lambda+w_0\mu+a\alpha_2+(\mu_2-1)\alpha_1) =
\underset{a=\max\{0,|\mu|-1-\lambda_2\}}{\overset{\mu_1-1}{\sum}}\dim V(\lambda''+w_0(\mu_1-1)\omega_1+a\alpha_2)
\\ \hspace{1.5cm} = \dim \ker(\lambda'',(\mu_1-1)\omega_1)
= \dim \cal F_{\lambda'',(\mu_1-1)\omega_1}- 
\dim \cal F_{\lambda''+\omega_1,(\mu_1-2)\omega_1}.\end{array}\label{mu.12}
\end{equation} 

\noindent Using \eqref{union.ninv}-
\eqref{mu.12}, for $\mu_2=k+1$, we can thus rewrite the inequality \eqref{dim.l.m} as follows, 
\begin{equation}\begin{array}{rl}4\dim \cal F_{\lambda,\mu} \leq & 4\biggl( \dim \cal F_{\lambda+\omega_2,\mu_1\omega_1+k\omega_2} +  \dim \cal F_{\lambda+(k+1)(\omega_1-\omega_2), \mu_1\omega_1} + \dim \cal F_{\lambda-\theta,(\mu_1-1)\omega_1+k\omega_2}\\ & -\dim \cal F_{\lambda-\theta+\omega_2,(\mu_1-1)\omega_1+(k-1)\omega_2} - \dim \cal F_{\lambda-\theta+k(\omega_1-\omega_2), (\mu_1-1)\omega_1} + \dim \cal F_{\lambda',k\omega_2}\\ &- \dim \cal F_{\lambda'+\omega_2,(k-1)\omega_2}+ \dim \cal F_{\lambda'',(\mu_1-1)\omega_1}- \dim \cal F_{\lambda''+\omega_1,(\mu_1-2)\omega_1}\biggr)\\
\leq & 4\biggl(\dim \cal F_{\lambda+\omega_2,\mu_1\omega_1+k\omega_2}+ \dim \cal F_{\lambda+(k+1)(\omega_1-\omega_2), \mu_1\omega_1}  + \dim \cal F_{\lambda-\theta,(\mu_1-1)\omega_1+k\omega_2}+ \dim \cal F_{\lambda',k\omega_2}\\&- \dim \cal F_{\lambda-\omega_1,(\mu_1-1)\omega_1+(k-1)\omega_2}  - \dim \cal F_{\lambda'+\omega_2,(k-1)\omega_2}- \dim \cal F_{(\lambda_1+k)\omega_1+(\lambda_2-k-1)\omega_2,(\mu_1-2)\omega_1}\biggr)
\end{array}\label{case2.rec}\end{equation}
where  $\lambda'$ and $\lambda''$ are as defined above. 
Now, using induction hypothesis, we have,
\begin{equation}
\begin{array}{rl}4\dim \cal F_{\lambda,\mu}\\
\leq &(\lambda_1+1)(\lambda_2+2)(|\lambda|+3)(\mu_1+1)(k+1)(\mu_1+k+2)\\& +
(\lambda_1+k+2)(\lambda_2-k)(|\lambda|+2)(\mu_1+1)(\mu_1+2)
+\lambda_1\lambda_2(|\lambda|)\mu_1(k+1)(\mu_1+k+1)\\ &-
\lambda_1(\lambda_2+1)(|\lambda|+1)\mu_1 k(\mu_1+k)
+(\lambda_1-\mu_1)(\lambda_2+\mu_1+1)(|\lambda|+1)(k+1)(k+2)\\ &-
(\lambda_1-\mu_1)(\lambda_2+\mu_1+2)(|\lambda|+2)(k)(k+1)
-(\lambda_1+k+1)(\lambda_2-k)(|\lambda|+1)(\mu_1-1)\mu_1\\
\leq &(\lambda_1+1)(\lambda_2+1)(|\lambda|+2)(\mu_1+1)(k+1)(\mu_1+k+2)\\ &+
(\lambda_1+1)(\lambda_2+2)(|\lambda|+3)(\mu_1+1)(\mu_1+2k+2) 
+(\lambda_2-\lambda_1-2k-1)(|\lambda|+2)(\mu_1+1)(\mu_1+2)\\ &
+\lambda_1\lambda_2(|\lambda|)\mu_1(\mu_1+2k+1)
-\lambda_1(\lambda_2+1)(|\lambda|+1)\mu_1 (\mu_1+2k-1)\\ &
+(\lambda_1-\mu_1)(\lambda_2+\mu_1+1)(|\lambda|+1)(2k+2)-
(\lambda_1-\mu_1)(\lambda_2+\mu_1+2)(|\lambda|+2)2k\\ &
-(\lambda_2-\lambda_1-2k)(|\lambda|+1)(\mu_1-1)\mu_1\\
\leq &(\lambda_1+1)(\lambda_2+1)(|\lambda|+2)(\mu_1+1)(k+1)
(\mu_1+k+2)\\ &+(\lambda_1+1)(\lambda_2+1)(|\lambda|+2)(\mu_1+1)
(\mu_1+2k+4) -2(\lambda_1+1)(\lambda_2+1)(|\lambda|+2)(\mu_1+1)\\ &+
(\lambda_1+1)(\lambda_1+2\lambda_2+4)(\mu_1+1)(\mu_1+2k+2)+2\lambda_1\lambda_2
(|\lambda|)\mu_1-\lambda_1\mu_1(\lambda_1+2\lambda_2+1)(\mu_1+2k-1) \\ &
+2(\lambda_1-\mu_1)(\lambda_2+\mu_1+1)(|\lambda|+1)
-(\lambda_1-\mu_1)(\lambda_1+\mu_1+2\lambda_2+3)2k\\ &+(|\lambda|+1)(\lambda_2-\lambda_1-2k)(4\mu_1+2)-2(\lambda_1+k+1)(\mu_1^2+3\mu_1+2)\\
\leq &(\lambda_1+1)(\lambda_2+1)(|\lambda|+2)(\mu_1+1)(k+2)
(\mu_1+k+3)-2(\lambda_1+1)(\lambda_2+1)(|\lambda|+2)(\mu_1+1)\\ &+
(\lambda_1+1)(\lambda_1+2\lambda_2+4)(\mu_1+1)(\mu_1+2k+2)
+2\lambda_1\lambda_2(|\lambda|)\mu_1-\lambda_1\mu_1(\lambda_1+
2\lambda_2+1)(\mu_1+2k-1) \\ &+2(\lambda_1-\mu_1)(\lambda_2+\mu_1+1)
(|\lambda|+1) -(\lambda_1-\mu_1)(\lambda_1+\mu_1+2\lambda_2+3)2k\\
&+(|\lambda|+1)(\lambda_2-\lambda_1-2k)(4\mu_1+2)-
2(\lambda_1+k+1)(\mu_1^2+3\mu_1+2)
\end{array}
\end{equation}
$$ \begin{array}{rl} \text{ Set } A =& -2(\lambda_1+1)(\lambda_2+1)(|\lambda|+2)(\mu_1+1)+ (\lambda_1+1)(\lambda_1+2\lambda_2+4)(\mu_1+1)(\mu_1+2k+2)
\\ &+2\lambda_1\lambda_2(|\lambda|)\mu_1-\lambda_1\mu_1(\lambda_1+
2\lambda_2+1)(\mu_1+2k-1) +2(\lambda_1-\mu_1)(\lambda_2+\mu_1+1)
(|\lambda|+1)\\& -(\lambda_1-\mu_1)(\lambda_1+\mu_1+2\lambda_2+3)2k +(\lambda_1+\lambda_2+1)(\lambda_2-\lambda_1-2k)(4\mu_1+2)\\ & -
2(\lambda_1+k+1)(\mu_1^2+3\mu_1+2) 
\end{array}$$
Then, $ 4 \dim \cal F_{\lambda, \mu_1\omega_1+(k+1)\omega_1}\leq (\lambda_1+1)(\lambda_2+1)(|\lambda|+2)(\mu_1+1)(k+2)(\mu_1+k+3)+A.$
Observe that the coefficient of $2k$ in $A$ is :
$$\begin{array}{rl}
&(\lambda_1+1)(\mu_1+1)(\lambda_1+2\lambda_2+4)-\lambda_1\mu_1(\lambda_1+2\lambda_2+1)-\lambda_1(\lambda_1+\mu_1+2\lambda_2+3) 
+\mu_1(\lambda_1+\mu_1+2\lambda_2+3)\\ &-(|\lambda|+1)(4\mu_1+2)
-(\mu_1^2+3\mu_1+2)\\
=& \lambda_1\mu_1(\lambda_1+2\lambda_2+4-\lambda_1-2\lambda_2-1-4)+(\lambda_1+\mu_1+1)(\lambda_1+2\lambda_2+4)-(\lambda_1^2-\mu_1^2)-(\lambda_1-\mu_1)(2\lambda_2+3)\\ &-2\lambda_1-2(\lambda_2+1)(2\mu_1+1)-(\mu_1^2+3\mu_1+2)\\
=& -\lambda_1\mu_1+\lambda_1(\lambda_1+\mu_1+1+2\lambda_2+4-2\lambda_2-3-2-\lambda_1)+2\lambda_2+4+\mu_1(2\lambda_2+3+2\lambda_2+4-4\lambda_2-4-3)\\ &-2(\lambda_2+1)-2=0\end{array}$$
Using the above relation we see that $A$ reduces to the following:  
$$\begin{array}{rl}A =& -2(\lambda_1+1)(\lambda_2+1)(|\lambda|+2)(\mu_1+1)+ (\lambda_1+1)(\lambda_1+2\lambda_2+4)(\mu_1+1)(\mu_1+2)\\ &
+2\lambda_1\lambda_2(|\lambda|)\mu_1-\lambda_1\mu_1(\lambda_1+2\lambda_2+1)(\mu_1-1) +2(\lambda_1-\mu_1)(\lambda_2+\mu_1+1)(|\lambda|+1) \\
&+(|\lambda|+1)(\lambda_2-\lambda_1)(4\mu_1+2)-2(\lambda_1+1)(\mu_1^2+3\mu_1+2) 
\\=& (\lambda_1+1)(\mu_1+1)[(\mu_1+2)(|\lambda|+2)+
(\mu_1+2)\lambda_2-2(\lambda_2+1)(|\lambda|+2)] + 
2\lambda_1\lambda_2(|\lambda|)\mu_1\\
&-\lambda_1\mu_1(1+2\lambda_2+\lambda_1)(\mu_1-1) +2(\lambda_1-\mu_1)(\lambda_2+\mu_1+1) (|\lambda|+1) +(|\lambda|+1)(\lambda_2-\lambda_1)(4\mu_1+2)\\
=&(\lambda_1+1)(\mu_1+1)[(|\lambda|+2)(\mu_1-2\lambda_2)-\lambda_2(\mu_1+2)] + 
\lambda_1\mu_1[2\lambda_2(|\lambda|)-(1+2\lambda_2+\lambda_1)(\mu_1-1)]\\ 
&+2(\lambda_1-\mu_1)(\lambda_2+\mu_1+1)
(|\lambda|+1) +(|\lambda|+1)(\lambda_2-\lambda_1)(4\mu_1+2)\\
=&\lambda_1\mu_1[\mu_1(|\lambda|)-2\lambda_2(|\lambda|)-4\lambda_2+\lambda_2\mu_1+2\lambda_2+2\lambda_2(|\lambda|)-\mu_1(\lambda_1+2\lambda_2)+2\lambda_2]+ 2\lambda_1\mu_1^2\\ &-\lambda_1\mu_1^2 +\lambda_1^2\mu_1+\lambda_1\mu_1 + (\lambda_1+\mu_1+1)[(\lambda_1+\lambda_2+2)(\mu_1-2\lambda_2)+\lambda_2(\mu_1+2)]\\ 
&+2(\lambda_1-\mu_1)(\lambda_2+\mu_1+1)
(|\lambda|+1) +(|\lambda|+1)(\lambda_2-\lambda_1)(4\mu_1+2)\\
= & (\lambda_1+\mu_1+1)[(|\lambda|+2)(\mu_1-2\lambda_2)+\lambda_2(\mu_1+2)+
\lambda_1\mu_1]\\ &+2(\lambda_1-\mu_1)(\lambda_2+\mu_1+1)
(|\lambda|+1) +(|\lambda|+1)(\lambda_2-\lambda_1)(4\mu_1+2)\\
=&2(|\lambda|+1)[(\lambda_1+\mu_1+1)(\mu_1-\lambda_2) + (\lambda_1-\mu_1)(\lambda_2+\mu_1+1)  + (\lambda_2-\lambda_1)(2\mu_1+1)]=0
\end{array}$$
Hence, \quad  $4 \dim \cal F_{\lambda, \mu_1\omega_1+(k+1)\omega_1}\leq (\lambda_1+1)(\lambda_2+1)(\lambda_1+\lambda_2+2)(\mu_1+1)(k+2)(\mu_1+k+3),$ which
along with \eqref{1} implies 
$\dim \cal F_{\lambda,\mu_1\omega_1+(k+1)\omega_1} = \dim V(\lambda)\dim V(\mu_1\omega_1+(k+1)\omega_1).$ 
This completes the proof of the proposition when $(\lambda, \mu)$ is partition of first kind with $\mu_2 \neq 0$.\endproof
\subsection{Proof of \propref{l*m.dim}(iii)}  Given a partition $(\lambda,\mu)$ of $\lambda+\mu$ of second kind, if $$\rho^\mu_\lambda=\mu_1\omega_1+\lambda_2\omega_2, \qquad \rho^\lambda_\mu=\lambda_1\omega_1+\mu_2\omega_2, $$ then observe that
$(\rho^\lambda_\mu, \rho^\mu_\lambda)\in P^+(\lambda+\mu,2)$ is a partition of first kind. 
\label{Case1.rec}
Using \propref{onto.filter}, and setting $${}_{_\ell}\mathbb{S}_{\lambda,\mu}= \bigcup\limits_{k=0}^{|\rho^\mu_\lambda|}\{(a_1, a_2): a_1 \in {S}_{\lambda,\mu}^{(\ell,k)}, a_2= k-a_1\}$$ 
we have, 
\begin{equation}\label{eq.inv}
\begin{array}{c}
\dim \mathcal F_{\lambda,\mu} = \dim \mathcal F_{\rho^\lambda_\mu,\rho^\mu_\lambda}+\dim \ker(\lambda,\mu),\\
\dim \ker (\lambda,\mu)
\leq 
\sum\limits_{\ell=1}^{\mu_2-\lambda_2}
\sum\limits_{(a_1,a_2)\in {_\ell}{\mathbb S}_{\lambda,\mu}} \dim V(\rho^\lambda_\mu+w_o\rho^\mu_\lambda-\ell\theta+a_2\alpha_1+a_1\alpha_2).
\end{array}\end{equation}
We now prove the result by applying induction on $\mu_2-\lambda_2.$
For $\mu_2-\lambda_2=1$, we have  $$\dim \ker(\lambda,\mu)= \dim \mathcal{F}_{\lambda,\mu}-\dim \mathcal F_{\rho^\lambda_\mu,\rho^\mu_\lambda},$$
$$\text{ and }\dim \ker(\lambda,\mu)\leq  \sum\limits_{(a_1,a_2)\in {_1}{\mathbb S}_{\lambda,\mu}} \dim V(\rho^\lambda_\mu-\theta+w_0\rho^\mu_\lambda+a_2\alpha_1+a_1\alpha_2).$$

\noindent If $\lambda, \mu\in P^+$ are such that $|\lambda|\geq |\mu|$ and $\mu_2=\lambda_2+1$, clearly  $\lambda_1>\mu_1$. Since $\lambda_1-1\geq \mu_1$, we see that $(\lambda-\omega_1,\mu-\omega_2)$ is a partition of $\lambda+\mu-\theta$ of first kind. On the other hand, we have,
$${_1}{\mathbb S}_{\lambda,\mu_1\omega_1+(\lambda_2+1)\omega_2}
=\left\lbrace (a_1,a_2)\in \mathbb Z^2_{\geq 0}: \begin{array}{l}0\leq a_2 \leq \lambda_2,\  0\leq a_1 \leq \mu_1,\\ 
\lambda_2-(\lambda_1-1)\leq a_2-a_1\leq (\mu_2-1)-\mu_1\end{array}\right\rbrace$$
$$\begin{array}{l}
=\mathbb S^{\lambda-\omega_1,\mu-\omega_2}\cup \{(a_1,\lambda_2)\in \mathbb Z_{>0}^2: 0\leq a_1\leq \mu_1, \lambda_2-(\lambda_1-1)\leq \lambda_2-a_1\leq \lambda_2-\mu_1 \}\\
=\mathbb S^{\lambda-\omega_1,\mu-\omega_2}\cup \{(\mu_1,\lambda_2)\} 
\end{array}$$
Putting $\mu_2=\lambda_2+1$ and  using  \propref{l*m.dim}(i)-(ii), we get,
$$\begin{array}{ll} \hspace{1cm}\sum\limits_{(a_1,a_2)\in {_1}{\mathbb S}_{\lambda,\mu}} \dim V(\rho^\lambda_\mu+w_0\rho^\mu_\lambda-\theta+a_2\alpha_1+a_1\alpha_2)\\
= \sum\limits_{(a_1,a_2)\in {\mathbb S}^{\lambda-\omega_1,\mu-\omega_2}} \dim V(\lambda-\omega_1+w_0(\mu-\omega_2)+a_2\alpha_1+a_1\alpha_2)+ \quad \dim  V((\lambda_1-\mu_1+\lambda_2-1)\omega_1+\mu_1\omega_2)\\
= \hspace{.5cm} \dim V((\lambda_1-\mu_1+\lambda_2-1)\omega_1+\mu_1\omega_2)+\dim \ker (\lambda-\omega_1,\mu-\omega_2). \end{array}$$
Note,  
$\dim \ker(\lambda-\omega_1,\mu-\omega_2)=\left\lbrace \small{\begin{array}{ll} 
	\dim \mathcal F_{\lambda-\omega_1,\mu-\omega_2} -\dim \mathcal F_{\lambda-\omega_1+\omega_2,\mu-2\omega_2} 
	-\dim \mathcal F_{\lambda-\omega_1+\lambda_2(\omega_1-\omega_2),\mu_1\omega_1}, & \text{ if }\, \lambda_2\geq 1,\\
	\dim \mathcal F_{\lambda-\omega_1,\mu-\omega_2} -\dim \mathcal F_{\lambda-\omega_1+\omega_1,\mu-\omega_2-\omega_1}, & \text{ if }\, \lambda_2=0.
	\end{array}}\right.\vspace{3mm}$
	Using \propref{l*m.dim}(ii), it thus follows from above that when $\lambda_2\geq 1$ and $\mu_2=\lambda_2+1$, 
	$$\begin{array}{rl}
\dim \mathcal{F}_{\lambda,\mu} \leq & 
\dim V(\lambda_1\omega_1+ (\lambda_2+1)\omega_2). \dim V(\mu_1\omega_1+\lambda_2\omega_2) 
+ \dim V((\lambda_1-1)\omega_1+\lambda_2\omega_2)
\dim V(\mu_1\omega_1+ \lambda_2\omega_2)\\
&- \dim V((\lambda_1-1)\omega_1+(\lambda_2+1)\omega_2). \dim V(\mu_1\omega_1+(\lambda_2-1)\omega_2)\\
&- \dim V((\lambda_1+\lambda_2-1)\omega_1) \dim V(\mu_1\omega_1) +\dim V((\lambda_1+\lambda_2-\mu_1-1)\omega_1+\mu_1\omega_2)\\ 
\leq & \biggl((\lambda_1+1)(\lambda_2+2)(|\lambda|+3)+
\lambda_1 (\lambda_2+1)(|\lambda|+1)\biggr)(\mu_1+1)(\lambda_2+1)(|\rho^\mu_\lambda|+2)  \\
&- \lambda_1(\lambda_2+2)(|\lambda|+2)(\mu_1+1)(\lambda_2)(|\rho^\mu_\lambda|+1))
+\biggl(2(|\lambda|-\mu_1)-(|\lambda|)(\mu_1+2)\biggr)(|\lambda|+1)(\mu_1+1)\\
\leq & \biggl( 3(\lambda_1+1)(\lambda_2+1)-(\lambda_1+2)\lambda_2\biggr)(|\lambda|+2)(\mu_1+1)(\lambda_2+1)(|\rho^\mu_\lambda|+2)  \\
&- (\lambda_1)(\lambda_2+2)(|\lambda|+2)(\mu_1+1)(\lambda_2)(|\rho^\mu_\lambda|+1))- (|\lambda|+1)(|\lambda|+2) (\mu_1)(\mu_1+1)\\
&\hspace{6.75cm} [\text{using induction step for $\mu_1=1$ in \secref{part1.dim}}]\\
\leq & \biggl( 3(\lambda_1+1)(\lambda_2+1)-(\lambda_1+2)\lambda_2-\lambda_1(\lambda_2+2)\biggr)(|\lambda|+2)(\mu_1+1)(\lambda_2+1)(|\rho^\mu_\lambda|+2)  \\
& + \lambda_1(\lambda_2+2)(|\lambda|+2)(\mu_1+1)(2\lambda_2+\mu_1+2)- 
(|\lambda|+1)(|\lambda|+2)\mu_1(\mu_1+1)
\\ \leq & \biggl((\lambda_1+1)(\lambda_2+1)+2\biggr)(|\lambda|+2)(\mu_1+1)(\lambda_2+1)(|\rho^\mu_\lambda|+2)  \\
& + \biggl( (\lambda_1)(\lambda_2+2)(2\lambda_2+\mu_1+2)- 
(|\lambda|+1)\mu_1\biggr)(|\lambda|+2) (\mu_1+1)\\
\leq & (|\lambda|+2)(\mu_1+1)(\lambda_2+1)(|\rho^\mu_\lambda|+2) [(\lambda_1+1)(\lambda_2+1)+2] \\
&+ (|\lambda|+2) (\mu_1+1) \biggl([(\lambda_1+1)(\lambda_2+1)+(\lambda_1-\lambda_2-1)](2\lambda_2+\mu_1+2)- (|\lambda|+1)\mu_1\biggr)\\
\leq & (|\lambda|+2)(\mu_1+1)[(\lambda_1+1)(\lambda_2+1)(\lambda_2+2)(|\rho^\mu_\lambda|+3)-2(\lambda_1+1)(\lambda_2+1)  \\
& +(\lambda_1-\lambda_2-1)(2\lambda_2+\mu_1+2) +2(\lambda_2+1)(|\rho^\mu_\lambda|+2)-\mu_1(|\lambda|+1)]\\
\leq & (\lambda_1+1)(\lambda_2+1)(|\lambda|+2)(\mu_1+1)(\lambda_2+2)(|\rho^\mu_\lambda|+3)
= 4\dim V(\lambda).\dim V(\mu_1\omega_1+(\lambda_2+1)\omega_2)
\end{array}$$

Along with \eqref{1}, the same arguments as used \propref{l*m.dim}(ii) show that  \propref{l*m.dim}(iii) is true when $\mu_2=\lambda_2+1$ and $\lambda_2\geq 1$. Similar computations show that the result holds when $\lambda_2=0$ and $\mu_2=1$. Now, by induction hypothesis assume that  \propref{l*m.dim}(iii) is true for all $(\lambda, \mu)\in P^+(\lambda+\mu,2)$ of second kind  with $\mu_2-\lambda_2\leq k$. 

Consider $(\lambda,\mu)\in P^+(\lambda+\mu,2)$ such that it is of second kind and  $\mu_2-\lambda_2=k+1$. Then  $(\lambda-\omega_1,\mu-\omega_2)\in P^+(\lambda+\mu-\theta,2)$ is a partition of second kind with $\mu_2-1-\lambda_2=k$.
Since for $1<\ell< \mu_2-\lambda_2$,
$$ \small{ \begin{array}{ll}
	{_{\ell+1}}{\mathbb S}_{\lambda,\mu} &= 
	\left\lbrace (a_1,a_2) \in \Z^2_{\geq 0}: 
	0\leq a_1 \leq \mu_1, \, \ 0\leq a_2 \leq \lambda_2, \mu_1-\mu_2+\ell+1\leq a_1-a_2\leq \lambda_1-\lambda_2-\ell-1 
	\right\rbrace,\\
	&={_{\ell}}{\mathbb S}_{\lambda-\omega_1,\mu-\omega_2},\end{array}} $$  applying induction hypothesis, we see that for all $\ell\leq k$,  
	\begin{equation}\begin{array}{r} \sum\limits_{\ell=2}^{\mu_2-\lambda_2}\sum\limits_{(a_1,a_2)\in {_\ell}{\mathbb S}_{\lambda, \mu}} \dim V(\rho^\lambda_\mu+w_0\rho^\mu_\lambda-(\ell+1)\theta+a_2\alpha_1+a_1\alpha_2)\\
	= \hspace{.3cm}\dim \, \ker(\lambda-\omega_1,\mu-\omega_2) \hspace{.4cm}= \hspace{.3cm} \dim \,  \mathcal F_{\lambda-\omega_1,\mu-\omega_2}-\, \dim\, \mathcal F_{\rho^\lambda_\mu-\theta,\rho^\mu_\lambda}\,.\end{array}\label{l>1}\end{equation}
	As $\lambda_1>\mu_1$ and $\mu_2-\lambda_2\geq 1$, we have,
	$${_1}{\mathbb S}_{\lambda,\mu}
	=\lbrace (a_1,a_2) \in \Z^2_{\geq 0}: 
	0\leq a_1\leq \mu_1,\, 0\leq a_2\leq \lambda_2,\,
	\lambda_2-(\lambda_1-1)\leq a_2-a_1 \leq (\mu_2-1)-\mu_1 \rbrace$$
	$${_1}{\mathbb S}_{\lambda,\mu}
	=\left \lbrace\small{\begin{array}{ll}
	\mathbb S^{\rho^\lambda_\mu-\theta,\rho^\mu_\lambda}
	\cup \,\, \{(a_1,\lambda_2):\max\{0,\mu_1+\lambda_2-\mu_2+1\}\leq a_1\leq \min \{\mu_1,\lambda_1-1\}\} , & \text{if}\, \, \lambda_2\neq 0,\\
	
	\mathbb S^{\rho^\lambda_\mu-\theta,\rho^\mu_\lambda}, & \text{if}\, \, \lambda_2= 0.\end{array}}\right.$$
	Further, using \propref{Case.1} and \propref{l*m.dim}(ii), we get
	\begin{equation}\label{l=1.a}
\begin{array}{l}\sum\limits_{(a_1,a_2)\in  {\mathbb S}^{\rho^\lambda_\mu-\theta, \rho^\mu_\lambda}} \dim V(\rho^\lambda_\mu-\theta+w_0\rho^\mu_\lambda+a_2\alpha_1+a_1\alpha_2) = \dim \ker (\rho^\lambda_\mu-\theta,\rho^\mu_\lambda)-(1-\delta_{\lambda_2,0})\dim \mathcal F_{\rho^\lambda_\mu-\theta+\lambda_2(\omega_1-\omega_2),\mu_1\omega_1}\\
	\hspace{1.4cm}= \quad \dim \mathcal F_{\rho^\lambda_\mu-\theta,\rho^\mu_\lambda}-\dim \mathcal F_{\rho^\lambda_\mu-\theta+\omega_2,\rho^\mu_\lambda-\omega_2}-(1-\delta_{\lambda_2,0})\dim \mathcal F_{\rho^\lambda_\mu-\theta+\lambda_2(\omega_1-\omega_2),\mu_1\omega_1}.
\end{array}
\end{equation}
In particular when $a_2=\lambda_2$, using \propref{onto.filter.1}(ii) and \propref{l*m.dim}(i), \eqref{l=1.a} reduces to :
\begin{equation}\begin{aligned} \label{l=1.b}
	\sum\limits_{a_1=\max\{0,|\rho^\mu_\lambda|-\mu_2+1\}}^{\mu_1} &\dim V(\rho^\lambda_\mu-\theta+\lambda_2(\omega_1-\omega_2)+w_0\mu_1\omega_1+a_1\alpha_2)
	= \dim \ker (\rho^\lambda_\mu-\theta+\lambda_2(\omega_1-\omega_2),\mu_1\omega_1)\\ 
	= &\dim \mathcal F_{\rho^\lambda_\mu-\theta+\lambda_2(\omega_1-\omega_2),\mu_1\omega_1}
	-\dim \mathcal F_{\rho^\lambda_\mu-\omega_2+\lambda_2(\omega_1-\omega_2),(\mu_1-1)\omega_1}.
\end{aligned}
\end{equation} 
Applying induction hypothesis and using \eqref{l>1}--\eqref{l=1.b} in 
\eqref{eq.inv}, we finally obtain the following inequalities:
\small{\begin{equation}\label{inv.inq}
	\dim \mathcal F_{\lambda,\mu}
	\leq \left\lbrace\begin{array}{ll}
		\dim \mathcal  F_{\rho^\lambda_\mu,\rho^\mu_\lambda} + 
		\dim \mathcal F_{\lambda-\omega_1,\mu-\omega_2}-
		\dim \mathcal F_{\rho^\lambda_\mu-\omega_2,(\mu_1-1)\omega_1},  & \text{if}\, \lambda_2=0,\\
		\dim \mathcal F_{\rho^\lambda_\mu,\rho^\mu_\lambda} + \dim \mathcal F_{\lambda-\omega_1,\mu-\omega_2}  -\dim \mathcal F_{\rho^\lambda_\mu-\omega_1,(\lambda_2-1)\omega_2},  & \text{if}\, \mu_1=0,\\
		\dim \mathcal F_{\rho^\lambda_\mu,\rho^\mu_\lambda} + \dim \mathcal F_{\lambda-\omega_1,\mu-\omega_2}-\dim \mathcal F_{\rho^\lambda_\mu-\omega_1,  \rho^\mu_\lambda-\omega_2}  
		-\dim \mathcal F_{\rho^\lambda_\mu-\lambda_2(\omega_2-\omega_1)-\omega_2,  (\mu_1-1)\omega_1} & \text{if}\, \mu_1, \lambda_2\neq 0.
	\end{array}\right.
	\end{equation}}
	\normalsize
	\noindent Since the cases when $\lambda_2=0$ and $\mu_1=0$ are can be obtained from the last inequality, it suffices to consider the case when $\mu_1,\lambda_2\neq 0$.
	
	Using 
	\propref{l*m.dim}(ii), it follows from \eqref{inv.inq} that
	when $\mu_1,\lambda_2\neq 0$ and $\mu_2-\lambda_2=k+1$, we have,
	\small{\begin{equation}\begin{aligned} 
		4  \dim  &{\mathcal  F}_{\lambda,\mu}\\
		&\leq  4 \biggl( \dim V(\lambda_1\omega_1+(\lambda_2+k+1)\omega_2)\dim V(\mu_1\omega_1+\lambda_2\omega_2)+ \dim V((\lambda_1-1)\omega_1+\lambda_2\omega_2) \dim V(\mu_1\omega_1+(\lambda_2+k)\omega_2)\\
		&\quad -\dim V((\lambda_1-1)\omega_1+(\lambda_2+k+1)\omega_2) \dim V(\mu_1\omega_1+(\lambda_2-1)\omega_2)-\dim V(|\lambda|\omega_1+k\omega_2)\dim V((\mu_1-1)\omega_1) \biggr)\\
		&	\leq  (\lambda_{1}+1)(\lambda_{2}+k+2)(|\lambda|+k+3)(\mu_1+1)(\lambda_{2}+1)(|\rho^\mu_\lambda|+2)\\
		&\quad +\lambda_1(\mu_1+1)
		\biggl((\lambda_{2}+1)(|\lambda|+1)(\lambda_{2}+k+1)(|\rho^\mu_\lambda|+k+2)
		-(\lambda_{2}+k+2)(\lambda_{1}\lambda_{2}+k+2)\lambda_{2}
		(|\rho^\mu_\lambda|+1)\biggr)\\ &\quad-(|\lambda|+1)(k+1)(|\lambda|+k+2)\mu_1(\mu_1+1)\\
		&	\leq  (\lambda_{1}+1)(\mu_1+1)(\lambda_{2}+1)(|\rho^\mu_\lambda|+2)\biggl((\lambda_{2}+k+1)(|\lambda|+k+2)+\lambda_{1}+2\lambda_{2}+2k+4\biggr)\\
		&\quad + \lambda_{1}(\lambda_{2}+1)(|\lambda|+1)(\mu_1+1)\biggl((\lambda_{2}+k)(|\rho^\mu_\lambda|+k+1)+ \mu_1+2\lambda_{2}+2k+2\biggr)\\
		&\quad-\lambda_1(\mu_1+1)\lambda_{2}(|\rho^\mu_\lambda|+1)\biggl((\lambda_{2}+k+1)(|\lambda|+k+1)+\lambda_{1}+2\lambda_{2}+2k+3\biggr)\\
		& \quad -(|\lambda|+1)\mu_1(\mu_1+1)\biggl(k(|\lambda|+k+1)+|\lambda|+2k+2\biggr)
	\end{aligned}\label{inv.inq.2}\end{equation}}
	\normalsize
	By induction hypothesis, when $\mu_2=\lambda_2+k$, 
	$$\begin{aligned}
\dim \mathcal F_{\lambda,\mu-\omega_2} =& \dim V(\lambda)\dim V(\mu-\omega_2)
= \dim V(\lambda+k\omega_2)\dim V(\rho^\mu_\lambda)
+  \dim V(\lambda-\omega_1) \dim V(\rho^\mu_\lambda+(k-1)\omega_2)\\
&-\dim V(\lambda-\omega_1+k\omega_2).\dim V(\rho^\mu_\lambda-\omega_2)
-\dim V(|\lambda|\omega_1+(k-1)\omega_2)\dim V((\mu_1-1)\omega_1)\\
\text{and } 2\dim V(\mu) = & 2 \dim V(\mu-\omega_2) + (\mu_1+1)(\mu_1+2\lambda_2+2k+4)\end{aligned} .$$ 

\noindent Thus, rearranging the coefficients in  inequality \eqref{inv.inq.2}, we get:
$$\begin{array}{rl}4 \dim \mathcal F_{\lambda,\mu}
\leq & 4\dim V(\lambda)\dim V(\mu-\omega_2)\\
&+(\lambda_1+1)(\lambda_2+1)(\mu_1+1)\biggl((|\lambda|+2)(\mu_1+2\lambda_2+2k+4)-(\lambda_1-\mu_1)(\lambda_2+2k+2)\biggr)\\
&-\lambda_1\lambda_2\biggl((|\lambda|+1)(\mu_1+2\lambda_2+2k+3)-(\lambda_1-\mu_1)(\lambda_2+2k+2)\biggr)(\mu_1+1)\\
&+(|\lambda|+1)\biggl(\lambda_1(\lambda_2+1)(\mu_1+2\lambda_2+2k+2)
-\mu_1(|\lambda|+2k+2)\biggr)(\mu_1+1)\\
\leq & 4 \, \dim V(\lambda)\dim V(\mu)-\biggl((\lambda_1+1)(\lambda_2+1)-\lambda_1\lambda_2\biggr)(\lambda_1-\mu_1)(\lambda_2+2k+2)(\mu_1+1)\\
&-\lambda_1\lambda_2(|\lambda|+1)(\mu_1+2\lambda_2+2k+3-\mu_1-2\lambda_2-2k-2)(\mu_1+1)\\
&+(|\lambda|+1)[\lambda_1(\mu_1+2\lambda_2+2k+2)-\mu_1(|\lambda|+2k+2)](\mu_1+1)\\
\leq &4 \, \dim V(\lambda)\dim V(\mu)-\biggl((\lambda_1+1)(\lambda_2+1)-\lambda_1\lambda_2\biggr)(\lambda_1-\mu_1)(\lambda_2+2k+2)(\mu_1+1)\\
&-\lambda_1\lambda_2(|\lambda|+1)(\mu_1+2\lambda_2+2k+3-\mu_1-2\lambda_2-2k-2)(\mu_1+1)\\
&+(|\lambda|+1)\biggl((\lambda_1-\mu_1)(\lambda_2+2k+2)+
\lambda_1\lambda_2\biggr)(\mu_1+1) = 4\, \dim V(\lambda)\dim V(\mu)
\end{array}$$
Along with \eqref{1}, the above inequality shows that \propref{l*m.dim}(iii)  holds when $\mu_2=\lambda_2+k+1$ and $\lambda_2\geq 1$ and  this completes the proof of the proposition.\endproof
\subsection{Proof of \thmref{filtration}} 
\noindent (i) Observe that part(i) of the theorem follows from \propref{surjective.l.m}, \propref{lem.case1}-\propref{case.inv} and \propref{l*m.dim}. \qed

\begin{cor} For any pair of distinct complex numbers $(z_1, z_2) \in \mathbb C^2$, there exists a $\mathfrak {sl}_3[t]$-module
isomorphism between ${\mathcal F}_{\lambda,\mu}$ and the fusion product $V(\lambda)^{z_1}\ast V(\mu)^{z_2} $.
\end{cor}
\proof By \lemref{prop of F_lambda}, for $(\lambda,\mu)\in P^+(\lambda+\mu,2)$, the fusion product $V(\lambda)^{z_1}\ast V(\mu)^{z_2}$ is a $\lie{sl}_3[t]$-quotient of $\cal F_{\lambda,\mu}$ and by \propref{l*m.dim}, $\dim \cal F_{\lambda,\mu}=\dim V(\lambda)\dim V(\mu)$. So using \remref{prop of F_lambda} we conclude that for any distinct pair of complex numbers 
$(z_1,z_2)$, $\cal F_{\lambda,\mu}$ is isomorphic to $V(\lambda)^{z_1}\ast V(\mu)^{z_2}$ as a $\lie{sl}_3[t]$-module. \qed 

\subsection{} \label{example} In this section, we show through an example that the methods used to study the modules $\mathcal F_{\lambda,\mu}$ can be extended to study the CV modules $\mathcal F_\blambda$ for $\blambda\in P^+(\lambda,k)$ with $k>2$.

Consider the case when $\lambda=(m+2)\omega_1$ and $\blambda = (m\omega_1, \omega_1,\omega_2)$. By definition,  $\mathcal F_\blambda = V(\bxi)$ where $\bxi=(\xi_{\alpha_1}, \xi_{\theta}, \xi_{\alpha_2})$ with $\xi_{\alpha_1}=\xi_\theta =(m,1,1)$ and $\xi_{\alpha_2}=\emptyset.$ Since for any distinct triple of complex 
numbers $(z_1,z_2,z_3)$, the fusion product  $V(m\omega_1)^{z_1}\ast V(\omega_1)^{z_2} \ast V(\omega_1)^{z_3}$ is a quotient of $\mathcal F_\blambda$, $$\dim \mathcal F_\blambda \geq \dfrac{9(m+1(m+2)}{2}.$$

On the other hand, using the defining relations of the CV modules, it is easy to see that there exists a $\mathfrak{sl}_3[t]$-module homomorphism $\phi: \mathcal F_\blambda {\rightarrow} \mathcal F_{m\omega_1,2\omega_1}$ such that $\ker \phi = \bu(\mathfrak g[t])(x_1^-\otimes t^2)v_{\blambda}.$ Since  $x_\alpha^+.(x_1^-\otimes t^2)v_{\blambda}=0,$ for all $\alpha^+\in R^+$ and $x_{\alpha_1}^-\otimes t^3.v_\blambda=0$, $\ker \phi$ is a quotient of $\tau_2^\ast W_{loc}(m\omega_1+\omega_2)$. Further similar computations as in \propref{X.r.s} show that $(x_\theta^-\otimes t)(x_{\alpha_1}^-\otimes t^2)v_\blambda\neq 0.$ Hence $\ker \phi$ is a quotient of $\mathcal F_{m\omega_1,\omega_2}$. Consequently,
$$\dim \cal F_\blambda \leq \dim \cal F_{m\omega_1,2\omega_2}+\dim \cal F_{m\omega_1,\omega_2}= \dfrac{6(m+1)(m+2)}{2}+\dfrac{3(m+1)(m+2)}{2}= \dfrac{9(m+1)(m+2)}{2}.$$
Thus we conclude that $\dim \mathcal F_\blambda=\dim V(m\omega_1)^{z_1}\ast V(\omega_1)^{z_2} \ast V(\omega_1)^{z_3}$ which implies that as $\mathfrak{sl}_3[t]$-modules $\mathcal F_\blambda$ is isomorphic to $V(m\omega_1)^{z_1}\ast V(\omega_1)^{z_2} \ast V(\omega_1)^{z_3}$. Since $\mathcal F_\blambda$ is defined by generators and relations, our method helps establish the Feigin-Loktev conjecture for triples of the form discussed in the above example. 




\section{Graded Character of Fusion product modules $\mathcal F_{\lambda,\mu}$}\label{graded char}
\subsection{} For $\lambda,\mu\in P^+$ the CV-module  $\mathcal F_{\lambda,\mu}$ is a $\mathbb Z_+$-graded vector space  and for each $s>0$, the subspace $\mathcal F_{\lambda, \mu}[s]$ is a finite-dimensional $\lie g$-module.
%
Hence by Weyl's theorem, $\mathcal F_{\lambda,\mu}[s]$ can be written as the direct sum of  $\tau^\ast_sV(\nu)$, with $\nu\in P^+$. The graded character of $\mathcal F_{\lambda,\mu}$ is the polynomial in indeterminate $q$ with coefficient in $\mathbb Z[P]$  given by  
$$\ch_{ gr}\mathcal F_{\lambda,\mu}= \sum\limits_{s \in \mathbb Z}\ch_{\lie g}\mathcal F_{\lambda,\mu}[s]\, q^s = \underset{(\nu,s)\in P^+\times \mathbb Z_+}{\sum} 
\ch_{\lie g}
(\tau^\ast_s V(\nu))\,  q^s. $$
Define $[\mathcal F_{\lambda, \mu} : V(\nu)]_q$ as the  polynomial in  indeterminate q given by,
$$[\mathcal F_{\lambda,\mu} : V(\nu)]_q =
\underset{p\geq 0}{\sum} \, 
[\mathcal F_{\lambda,\mu} : \tau^\ast_p(V(\nu))]\,  q^p,
$$
where $[\mathcal F_{\lambda, \mu} : \tau^\ast_p V(\nu)]$ is
the multiplicity of $V(\nu)$ in the graded component $\cal F_{\lambda,\mu}[p]$ of
$\mathcal F_{\lambda, \mu}$. The polynomial 
$[\mathcal F_{\lambda, \mu} : V(\nu)]_q $
is called the graded multiplicity of $V(\nu)$
in $\mathcal F_{\lambda,\mu}$
and at q = 1 it gives the numerical multiplicity of $V(\nu)$ in the $\mathfrak g$-module $V(\lambda)\otimes V(\mu)$, which is   equal to the Littlewood Richardson coefficients $c^\nu_{\lambda,\mu}$. In this section, we use the graded character of the module $\mathcal F_{\lambda,\mu}$ to obtain an algebraic characterization of the Littlewood Richardson coefficients in type $A_2$. 
\label{section littlewood}

\begin{thm}\label{littlewood coeff} For $\lambda,\mu,\eta\in P^+$,  let $c_{\lambda,\mu}^\eta$ be the multiplicity of $V(\eta)$ in the $\mathfrak{sl}_3(\mathbb C)$-module $V(\lambda)\otimes V(\mu)$. 
\begin{itemize}
	\item[i.] Suppose $(\lambda,\mu)$ is a partition of  $\lambda+\mu$ of first kind with $\mu=\mu_1\omega_1$. Then
	$$c_{\lambda,\mu}^\eta = \left\lbrace \small{\begin{array}{ll} 1,& \text{if} \,\, \nu=\lambda+w_0\mu+j\theta+a_j\alpha_2 \text{ for } 0\leq j\leq \mu_1, \max\{0,\mu_1-j-\lambda_2\}\leq a_j\leq \mu_1-j, \\
			0, &\text{ otherwise}, \end{array}}\right.$$ 
	\item[ii.] Suppose $(\lambda,\mu)$ is a partition of $\lambda+\mu$ of first kind with $\mu_i>0$ for $i=1,2$. 
	Then 
	$$c_{\lambda,\mu}^\eta= \left\lbrace
	\begin{array}{ll} 
		a+1, & \text{if} \, \, \nu=\lambda+w_0\mu+(\ell+j)\theta+(\mu_2-j)\alpha_1+a\alpha_2,\,  \text{with}\, (j,\ell,a)\in A^{\lambda,\mu} \\
		\min\{a,b\}+1, & \text{if} \, \, \nu=\lambda+w_0\mu+j\theta+b\alpha_1+a\alpha_2,\, \text{with}\, (j,a,b)\in B^{\lambda,\mu} \\
		0 & \text{ otherwise}
	\end{array}\right.$$
	where
	$$\begin{array}{rll} A^{\lambda,\mu} & = &\{(j,\ell,a)\in \mathbb Z^3 : j=0,\, 0\leq \ell\leq \mu_1,\, \max\{0,|\mu|-\lambda_2-\ell\} \leq a\leq \mu_1-\ell\}\\ & \sqcup & \{(j,\ell, a)\in \mathbb Z^3 : 1\leq j\leq \mu_2,\, 0\leq \ell\leq \mu_1,\,  a = \mu_1-\ell\}\\
		& \sqcup &\{(j,\ell,a)\in \mathbb Z^3 : 1\leq j\leq \mu_2,\, 0\leq \ell\leq \mu_1, \, \ell+2j \leq |\mu|-\lambda_2,\, a = |\mu|-\lambda_2-(\ell+2j) \}\\
		B^{\lambda,\mu}&=&\{(j,a,b)\in \mathbb Z^3 : j=0, 0\leq a\leq \mu_1,\, 0\leq b\leq \mu_2-1, \,
		\mu_2-\lambda_1\leq b-a\leq \lambda_2-\mu_1 \}\\
		& \sqcup & \{(j,a,b)\in \mathbb Z^3 : 1\leq j\leq \mu_2-1,\, 0\leq a\leq \mu_1, 0\leq b < \mu_2-j, b-a=\lambda_2-\mu_1+j \}\\
		& \sqcup & \{(j,a,b)\in \mathbb Z^3 : 1\leq j\leq \mu_2-1,\, 0\leq a\leq \mu_1, 0\leq b < \mu_2-j, b-a= \mu_2-\lambda_1-j \}.
	\end{array}$$
	\item[iii.]  Suppose $(\lambda,\mu)$ is a partition of  $\lambda+\mu$ of second kind. 
	Then
	$$c_{\lambda,\mu}^\eta= \left\lbrace \begin{array}{ll}
		\min\{a,b\}+1, & \text{if } \nu= \rho^\lambda_\mu+ w_0 \rho^\mu_\lambda- \ell \theta +a \alpha_2+ b \alpha_1, (\ell,a,b) \in C_{\lambda,\mu}\\
		\min\{b_1,b_2\}+1, & \text{if } \nu= \lambda+ w_0 \mu+j \theta +b_1 \alpha_2+ b_2 \alpha_1, (j,b_1,b_2) \in B_{\lambda,\mu} \\
		a+1,  & \text{if } \nu=\rho^\lambda_\mu+w_0\rho^\mu_\lambda+(\ell+j)\theta+(\lambda_2-j)\alpha_1+a \alpha_2, (j,\ell,a) \in A^{\rho^\lambda_\mu, \rho^\mu_\lambda} \\
		0 & \text{otherwise }
	\end{array}\right.$$
	$$\begin{array}{rll}\text{where }A^{\rho^\lambda_\mu, \rho^\mu_\lambda} 	&  =& \{(j,\ell,a)\in \mathbb Z^3 : j=0,\, 0\leq \ell\leq \mu_1,\, \max\{0,|\rho^\mu_\lambda|-\mu_2-\ell\} \leq a\leq \mu_1-\ell\}\\ & \sqcup & \{(j,\ell, a)\in \mathbb Z^3 : 1\leq j\leq \lambda_2,\, 0\leq \ell\leq \mu_1,\,  a = \mu_1-\ell\}\\
		& \sqcup &\{(j,\ell,a)\in \mathbb Z^3 : 1\leq j\leq \lambda_2,\, 0\leq \ell\leq \mu_1, \, \ell+2j \leq |\rho^\mu_\lambda|-\mu_2,\, a = |\rho^\mu_\lambda|-\mu_2-(\ell+2j) \}\\
		B_{\lambda,\mu}&=&\{(r,a_1,a_2)\in \mathbb Z^3 : r=0, 0\leq a\leq \mu_1,\, 0\leq b\leq \lambda_2, \, ~
		\lambda_2-\lambda_1\leq b-a\leq \mu_2-\mu_1 \}\\
		& \sqcup&  \{(r,a_1,a_2)\in \mathbb Z^3 : 1\leq r\leq \mu_2-\lambda_2-1,\, 0\leq a\leq \mu_1, 0\leq b \leq \lambda_2, ~ b-a=\mu_2-\lambda_1-r \}\\
		& \sqcup&  \{(r,a_1,a_2)\in \mathbb Z^3 : 1\leq j\leq \mu_2-\lambda_2-1,\, 0\leq a_1\leq \mu_1, ~ 0\leq a_2 < \lambda_2, ~ a_2-a_1= \lambda_2-\mu_1+r \}\\
		C_{\lambda,\mu} &= &\{(\ell, a,b): \ell=0, 0 \leq a \leq \mu_1, 0 \leq b\leq \lambda_2, \mu_2-\lambda_1 \leq a_2-a_1 \leq \lambda_2-\mu_1\}\\
		& \sqcup& \{(\ell, a,b): 1 \leq \ell \leq \mu_2-\lambda_2-1, 0 \leq a \leq \mu_1, 0 \leq b\leq \lambda_2, b-a= \mu_2-\lambda_1-\ell\}\\
		& \sqcup & \{(\ell, a,b): 1 \leq \ell \leq \mu_2-\lambda_2-1, 0 \leq a \leq \mu_1, 0 \leq b\leq \lambda_2, b-a= \lambda_2-\mu_1+\ell\}.
	\end{array}$$
\end{itemize}
\end{thm}
\subsection{Proof of \thmref{littlewood coeff}(i):}  \label{first.kind.1} Let $(\lambda, \mu)\in P^+(\lambda+\mu,2)$, be a partition 
of first kind and $\mu_2 =0$. Using \propref{Case.1}(ii) and  \propref{l*m.dim}(i), we have
$${\mathcal{F}}_{\lambda, \mu_1 \omega_1}={\mathcal{F}}_{\lambda+\omega_1,(\mu_1-1)\omega_1}
\oplus \ker \phi (\lambda, \mu_1\omega_1),$$
$$\ker \phi(\lambda,\mu_1\omega_1) \cong_{\mathfrak{sl}_3[t]} {\overset{\mu_1}{\underset{a=\max\{0,\mu_1-\lambda_2\}}{\bigoplus}}} 
\tau^\ast_{\mu_1} V(\lambda-a\omega_1-(\mu_1-2a)\omega_2) = \bigoplus\limits_{a \in \max\{0,\mu_1-\lambda_2\}}^{\mu_1} \tau^\ast_{\mu_1}V(\lambda+w_0\mu+a\alpha_2)$$
Further, using the $\mathfrak{sl}_3[t]$-module decomposition for the successive quotients, $\mathcal F_{\lambda+j\omega_1,(\mu_1-j)\omega_1}$, $1\leq j\leq \mu_1$, we gsee that $$ {\mathcal{F}}_{\lambda,\mu}\cong_{\lie{sl}_3[t]} {\overset{\mu_1}{\underset{j=0}{\bigoplus}}}
{\overset{\mu_1-j}{\underset{a_j=\max\{0,\mu_1-j-\lambda_2\}}{\oplus}}} \tau^\ast_{\mu_1-j} V(\lambda+w_0\mu+j \theta+a_j\alpha_2).$$
As a consequence  the graded character of $\cal F_{\lambda, \mu_1\omega_1}$ is given as follows:
$$\ch_{gr} \cal F_{\lambda, \mu_1\omega_1} = {\overset{\mu_1}{\underset{j=0}{\sum}}}
\left({\overset{\mu_1-j}{\underset{a_j=\max\{0,\mu_1-j-\lambda_2\}}{\sum}}}
\ch_{\mathfrak{sl_3(\mathbb C)}}\tau^\ast_{\mu_1-j}  V(\lambda+w_0\mu+j \theta+a_j\alpha_2)\right) q^{\mu_1-j}. $$
Since $\theta$ and $\alpha_2$ are linearly independent roots, for pairs of integers $(j, a_j)$ and $(k, a_k)$,
$$\lambda+w_0\mu+j \theta+a_j\alpha_2 = \lambda+w_0\mu+k \theta+a_k\alpha_2,$$ only if $j=k$ and $a_j=a_k$. Hence, we have
\begin{equation}\label{a-p.def}\begin{array}{ll}
	[\mathcal F_{\lambda, \mu}: V(\nu)]_q  =\left\{\small{\begin{array}{ll} q^{\mu_1-j}, & \text{if} \,\, \nu=\lambda+w_0\mu+j\theta+a_j\alpha_2  \text{ with } 0\leq j\leq \mu_1, \max\{0,\mu_1-j-\lambda_2\}\leq a_j\leq \mu_1-j \\
			0, & \text{ otherwise}. \end{array}}\right. \end{array}\end{equation}
			Putting $q=1$ in \eqref{a-p.def} we obtain \thmref{littlewood coeff}(i).
			
			\label{case1.ch}
			
			\subsection{Proof of \thmref{littlewood coeff}}\label{first.kind.2} Let $(\lambda, \mu)\in P^+(\lambda+\mu,2)$ be a partition of $\lambda+\mu$ of first kind with $\mu_1,\mu_2 \neq 0$. Then, from \propref{onto.filter.1}(i) and  \propref{l*m.dim}(ii),  it follows that,
			$${\mathcal{F}}_{\lambda, \mu}={\mathcal{F}}_{\lambda+\omega_2,\mu-\omega_2}
			\oplus \ker \phi(\lambda, \mu),$$ 
			$$\ker \phi(\lambda, \mu) \cong_{\mathfrak{sl}_3[t]} 
			\tau^\ast_{\mu_2} {\mathcal{F}}_{\lambda+\mu_2(\omega_1-\omega_2), \mu_1\omega_1} \oplus 
			{{\underset{(a_1,a_2)\in {\mathbb S}^{\lambda, \mu}}{\bigoplus}}} \tau^\ast_{|\mu|} V(\lambda+w_0\mu+a_2\alpha_1+a_1\alpha_2)$$
			Using the $\mathfrak{sl}_3[t]$-module decomposition for the successive quotients $\mathcal F_{\lambda+j\omega_2,\mu-j\omega_2}$, $1\leq j\leq \mu_2$, 
			we have, 
			$$\begin{array}{rl}
{\mathcal{F}}(\lambda, \mu)
&= \overset{\mu_2-1}{\underset{j=0}{\bigoplus }}  
\ker \phi (\lambda+j\omega_2,\mu-j\omega_2) \, \oplus \,  {\mathcal{F}}
_{\lambda+\mu_2\omega_2,\mu_1\omega_1}\\

&\cong_{\lie{sl}_3[t]} \overset{\mu_2-1}{\underset{j=0}{\bigoplus }}  
\tau^\ast_{\mu_2-j} \cal F_{\lambda+j\omega_2+(\mu_2-j)(\omega_1-\omega_2), \mu_1\omega_1}\oplus {\mathcal{F}}_{\lambda+\mu_2\omega_2,\mu_1\omega_1}\\
& \quad \qquad	\oplus  \bigoplus\limits_{j=0}^{\mu_2-1} \left( 
\bigoplus\limits_{(a_j, b_j)\in \mathbb S^{\lambda+j\omega_2, \mu-j\omega_2}}
\tau^\ast_{|\mu|-j} 
V(\lambda+w_0\mu +j\theta+b_j\alpha_1+a_j\alpha_2) \right)
\end{array}$$

For $(\ell, j)\in \mathbb Z^2$ , let $M^\lambda_\mu(\ell,j)=\max\{0,|\mu|-\ell-\lambda_2-2j\}.$
Using \ref{case1.ch} we have :
$$\begin{array}{ll}
\ch_{gr} \, \cal F_{\lambda,\mu} &= \sum\limits_{j=0}^{\mu_2}
\ch_{gr} \cal F_{\lambda+j\omega_2+(\mu_2-j)(\omega_1-\omega_2), \mu_1\omega_1} q^{\mu_2-j}\\ 
&+ \sum\limits_{j=0}^{\mu_2-1} \left( \sum\limits_{(a_j, b_j)\in {\mathbb S}^{\lambda+j\omega_2,\mu-j\omega_2}} \ch_{\lie g} \tau^\ast_{|\mu|-j} V(\lambda+w_0\mu+j\theta+b_j\alpha_1+a_j\alpha_2)\right)q^{|\mu|-j}\\
&= \sum\limits_{j=0}^{\mu_2} \left(\sum\limits_{\ell=0}^{\mu_1}
\left(\sum\limits_{a_\ell^j=M^\lambda_\mu(\ell,j)}^{\mu_1-\ell}
\ch_{\mathfrak{g}} \tau^\ast_{|\mu|-\ell-j} V(\lambda+w_0\mu+(\ell+j)\theta+(\mu_2-j)\alpha_1+a^j_\ell\alpha_2)	\right)q^{|\mu|-\ell-j}\right) \\ 
&+ \sum\limits_{j=0}^{\mu_2-1}  
\left( \sum\limits_{(a_j, b_j)\in {\mathbb S}^{\lambda+j\omega_2,\mu-j\omega_2}}
\ch_{\lie g} \tau^\ast_{|\mu|-j} V(\lambda+w_0\mu+j \theta+b_j\alpha_1+a_j\alpha_2)\right)q^{|\mu|-j}
\end{array}$$

Comparing the coefficients of $\alpha_1$ and $\alpha_2$ in the highest weights of the $\mathfrak g$-irreducible components of $\cal F_{\lambda,\mu}$, we now determine the polynomials $[\cal F_{\lambda,\mu}:V(\nu)]_q$. \begin{enumerate}
\item[(i).] Given 
$(a_1,a_2)\in \mathbb S^{\lambda, \mu}$, and integers $(\ell, j)$ such that  $0\leq j\leq \mu_2$ and $0\leq \ell \leq \mu_1$, 
$$\lambda+w_0\mu+(\ell+j)\theta+(\mu_2-j)\alpha_1+a^j_\ell\alpha_2 = 
\lambda+w_0\mu+a_2\alpha_1+a_1\alpha_2$$
only if  $(a_1,a_2)=(j+\ell+a_\ell^j, \mu_2+\ell)$. However, by definition of $\mathbb S^{\lambda, \mu}$, $a_2<\mu_2$, therefore this case cannot occur.  
\item[(ii).] For pair of integers $(j,\ell)$ such that $0\leq j \leq \mu_2$ and $0\leq \ell \leq \mu_1$ and $(a_{j_1}, a_{j_2})\in 
\mathbb S^{\lambda+j\omega_2,\mu-j\omega_2}$,  
$$\lambda+w_0\mu +\ell\theta+\mu_2\alpha_1+a_\ell\alpha_2 = 
\lambda+w_0\mu+j \theta +a_{j_2}\alpha_1+a_{j_1}\alpha_2$$
only if $(a_{j_1}, a_{j_2})=(\ell-j+a_\ell, \mu_2+\ell-j)$. But by definition of $\mathbb S^{\lambda+j\omega_2,\mu-j\omega_2}$, $a_{j_2}<\mu_2-j$, therefore  this case cannot occur.
\item[(iii).] Given triples of integers, $(j, \ell, a_\ell^j)$ and $(j+s, r, a_r^{j+s})$, with $0\leq j\leq j+s \leq \mu_2$,  $\ell, r \in [0,\mu_1]$ and 
$\max\{0, |\mu|-k-\lambda_2-2s\} \leq a_k^s\leq \mu_1-k$ for $s\in \{\ell,r\}$ and $k\in\{j, j+s\}$, 
$$ \begin{array}{l} 
	\lambda+w_0\mu+(\ell+j)\theta+(\mu_2-j)\alpha_1+a^j_\ell\alpha_2
	=\lambda+w_0\mu+(r+j+s)\theta+(\mu_2-j-s)\alpha_1+a^{j+s}_r\alpha_2
\end{array}$$ 
only if $\ell=r,$ and $a_\ell^j=a_r^{j+s}+s=a_\ell^{j+s}+s.$
Note that if $a_\ell^j$ is a non-negative integer in the interval $[|\mu|-\lambda_2-\ell-2j,\mu_1-\ell]$, then $a_\ell^j+1$ also lies in the interval $[|\mu|-\lambda_2-\ell-2j,\mu_1-\ell]$ unless $a_\ell^j\in \{\mu_1-\ell, |\mu|-\lambda_2-\ell-2j\}$. Thus,
\begin{equation} \label{l.P}\begin{array}{lr}
		[\mathcal F_{\lambda,\mu}: V (\lambda+w_0\mu+(\ell+j)\theta+(\mu_2-j)\alpha_1+a \alpha_2)]_q = \sum\limits_{s=0}^{a} q^{|\mu|-\ell-s-j}
\end{array}\end{equation}
whenever $(j,\ell, a)$ is a triple of non-negative integers contained in the set $A^{\lambda,\mu}$, which is defined in \thmref{littlewood coeff} (ii). 
	
	\item[(iv).]  Given triples of integers, $(j, a_j, b_j)$ and $(j+s, a_{j+s}, b_{j+s})$, with $0\leq j\leq  j+s\leq \mu_2-1$, $(a_j, b_j) \in \mathbb{S}^{\lambda+j \omega_2,\mu-j\omega_2}$ and $(a_{j+s}, b_{j+s}) \in \mathbb{S}^{\lambda+(j+s) \omega_2,\mu-(j+s)\omega_2}$,
	$$\begin{array}{l}\lambda+w_0\mu+j \theta+b_j\alpha_1+a_j\alpha_2
		=\lambda+w_0\mu+(j+s)\theta+b_{j+s}\alpha_1+a_{j+s}\alpha_2
	\end{array}$$ only if $(a_{j}, b_{j})=(a_{j+s}+s, b_{j+s}+s)$. From the
	definition of the set $\mathbb S^{\lambda+j\omega_2,\mu-j\omega_2}$, we see that if  
	$(a, b)\in \mathbb S^{\lambda+(j+1)\omega_2,\mu-(j+1)\omega_2},$ for some $0\leq j\leq \mu_2$,
	then $(a+1,b+1)\in \mathbb S^{\lambda+j\omega_2,\mu-j\omega_2}$ unless $b-a=\lambda_2+j-\mu_1$ or $b-a=\mu_2-\lambda_1-j$. Thus, 
	\begin{equation}\label{ninv.i}
		\begin{array}{lr} [\mathcal F_{\lambda, \mu}: V(\lambda+w_0\mu+j \theta+b\alpha_1+a\alpha_2)]_q = \sum\limits_{s=0}^{\min\{a, b\}}q^{|\mu|-s-j}, 
	\end{array}\end{equation}
	whenever $(j,a,b)$ is a triple of non-negative integers contained in the set $B^{\lambda,\mu}$, which is defined in \thmref{littlewood coeff} (ii).
	
	\noindent Putting $q=1$ in \eqref{l.P} and \eqref{ninv.i} and observing that in this case 
	$$[\mathcal F_{\lambda, \mu}: V(\nu)]_{q=1} =0, \qquad \text{when } \nu\notin A^{\lambda,\mu}\cup B^{\lambda,\mu}$$ we obtain \thmref{littlewood coeff}(ii). 
\end{enumerate}
\subsection{Proof of \thmref{littlewood coeff}(iii).} \label{second.kind} Let $(\lambda, \mu)\in P^+(\lambda+\mu,2)$ be a partition  of second kind. Then, from  \propref{onto.filter} and  \propref{l*m.dim}(iii) it follows that,
$${\mathcal{F}}_{\lambda, \mu}={\mathcal{F}}_{\rho^\lambda_\mu,\rho^\mu_\lambda} 
\oplus \ker \phi(\lambda, \mu),$$
$$\ker \phi(\lambda, \mu) \cong_{\mathfrak{sl}_3[t]} 
\bigoplus\limits_{\ell=1}^{\mu_2-\lambda_2}
\bigoplus\limits_{(a, b)\in {_\ell}{\mathbb S}_{\lambda, \mu}} 
\tau^\ast_{\mu_1+\lambda_2+\ell}V(\lambda+w_0\mu+(\mu_2-\lambda_2-\ell)\theta+b\alpha_1+a\alpha_2),$$ where  $\rho^\lambda_\mu = \lambda_1\omega_1+\mu_2\omega_2,$ and 
$\rho^\mu_\lambda = \mu_1\omega_1+\lambda_2\omega_2$.   
In this case, as $(\rho^\lambda_\mu,\rho^\mu_\lambda)$ is a partition of 
$\lambda+\mu$ of first kind, using \ref{first.kind.2},  we get,
\small{$$\begin{aligned}
		\ch_q \cal F_{\lambda, \mu}=& \sum\limits_{\ell=1}^{\mu_2-\lambda_2}q^{|\rho^\mu_\lambda|+\ell} \left( \sum\limits_{(a, b) \in {}_\ell{\mathbb S}_{\lambda, \mu}} \ch_{\mathfrak{g}} \tau^\ast_{\mu_1+\lambda_2+\ell}V(\rho^\lambda_\mu+ w_0 \rho^\mu_\lambda- \ell \theta +b \alpha_2+ a \alpha_1) \right)\\
		&+\sum\limits_{j=0}^{\lambda_2}\left(\sum\limits_{\ell=0}^{\mu_1}\left(\sum\limits_{a_\ell^j={L}^{\lambda, \mu}(\ell, j)}^{\mu_1-\ell} \ch_{\mathfrak{g}}\tau^\ast_{\mu_1+\lambda_2-\ell-j}V( \rho^\lambda_\mu+ w_0 \rho^\mu_\lambda+ (\ell+j) \theta+(\lambda_2-j)\alpha_1+ a^j_\ell \alpha_2) \right)\right)q^{|\rho^\mu_\lambda|-\ell-j}\\
		&+\sum\limits_{j=0}^{\lambda_2-1}\left(\sum\limits_{(a, b) \in {\mathbb S}^{\rho^\lambda_\mu+j \omega_2,\rho^\mu_\lambda-j\omega_2}} \ch_{\mathfrak{g}}\tau^\ast_{\mu_1+\lambda_2-j}V(\rho^\lambda_\mu+ w_0 \rho^\mu_\lambda+ j \theta +b \alpha_1+ a \alpha_2) \right)
		q^{|\rho^\mu_\lambda|-j}
	\end{aligned}$$}
\normalsize
where, for $(\ell,j)\in \mathbb Z^2$, $L^{\lambda,\mu}(\ell,j)=M^{\rho^\lambda_\mu}_{\rho^\mu_\lambda}(\ell,j)=\max\{0,|\rho^\mu_\lambda|-\mu_2-\ell-2j\}$. 
Now, by comparing the coefficients of $\omega_1$ and $\omega_2$ in the highest weights of the $\mathfrak g$-irreducible 
components of $\cal F_{\lambda,\mu}$ we determine the polynomials $[\mathcal F_{\lambda,\mu};V(\nu)]_q$. 
\noindent To begin with we observe that given a positive integer $\ell\in [1,\mu_2-\lambda_2]$, there exists an integer $r\in [0,\mu_2-\lambda_2-1]$ such that $\ell=\mu_2-\lambda_2-r$. Hence it follows from the definition of the set ${_\ell}{\mathbb S}_{\lambda,\mu}$ (ref. \propref{Case.3},  \propref{l*m.dim}(iii)), that if  $(a_1,a_2)\in {_\ell}{\mathbb S}_{\lambda,\mu}$ for $\ell=\mu_2-\lambda_2-r$, then $a_1,a_2$ are integers such that 
$0\leq a_1\leq \mu_1,\, 0\leq a_2\leq \lambda_2$ and 
$\mu_2-\lambda_1-r\leq a_2-a_1\leq \lambda_2-\mu_1+r$.
Therefore,
	\small{\begin{equation}\begin{array}{ll} {}_{\mu_2-\lambda_2-r-1}{\mathbb S}_{\lambda,\mu}= {}_{\mu_2-\lambda_2-r}{\mathbb S}_{\lambda,\mu}\cup \left\{(a_1,a_2)\in \mathbb Z^2: \begin{array}{l}0\leq a_1\leq \mu_1,\ 0\leq a_1\leq \lambda_2,\\
					a_2-a_1\in \{ \mu_2-\lambda_1-r-1,\lambda_2-\mu_1+r+1\}\end{array} \right\}.
			\end{array}\label{inv.S}\end{equation}}
	\normalsize
	
	\item[(i).] Given triplet of integers $(r, a_1,a_2)$ with  $1 \leq r \leq \mu_2-\lambda_2, \  (a_1,a_2) \in {}_r {\mathbb{S}_{\lambda, \mu}}$
	and  $(\ell, j, a_\ell^j)$ with $0 \leq \ell  \leq \mu_1,\ 
	, \ 0 \leq j \leq \lambda_{2}, \ a_\ell^j\in [|\rho^\mu_\lambda|-\lambda_2-\ell-2j,\mu_1-\ell]\cap\mathbb{Z}_+$, 
	$$\rho^\lambda_\mu+ w_0 \rho^\mu_\lambda- r \theta +a_1 \alpha_2+ a_2 \alpha_1= \rho^\lambda_\mu+ w_0 \rho^\mu_\lambda+ (\ell+j) \theta+(\lambda_2-j)\alpha_1+ a^j_\ell \alpha_2 $$
	only if	$(a_1,a_2)=(a_\ell^j+j+r+\ell, \lambda_2+\ell+r)$. This cannot happen since $a_2 \leq \lambda_2$ and $r\leq 1$. Hence it follows from \eqref{l.P} that 
	for $\nu\in A^{\rho^\lambda_\mu, \rho^\mu_\lambda}$,
	\begin{equation} \label{l.P.inv}\begin{array}{lr}
			[\mathcal F_{\lambda,\mu}: V (\rho^\lambda_\mu+w_0\rho^\mu_\lambda+(\ell+j)\theta+(\lambda_2-j)\alpha_1+a \alpha_2)]_q = \sum\limits_{s=0}^{a} q^{|\rho^\mu_\lambda|-\ell-s-j}
	\end{array}\end{equation} 
	\item[(ii)] Given triplets of integers $(\mu_2-\lambda_2-\ell, a_1,a_2)$ and $(\mu_2-\lambda_2-\ell-s,b_1,b_2)$ with 
	$0\leq  \ell\leq  \ell+s \leq \mu_2-\lambda_2-1$, 
	$(a_1,a_2) \in {}_{\mu_2-\lambda_2-\ell}\mathbb{S}_{\lambda,\mu}$ and $(b_1,b_2) \in {}_{\mu_2-\lambda_2-\ell-s}\mathbb{S}_{\lambda,\mu}$, 
	$$\rho^\lambda_\mu+ w_0 \rho^\mu_\lambda- (\mu_2-\lambda_2-\ell) \theta +a_1 \alpha_2+ a_2 \alpha_1 = \rho^\lambda_\mu+ w_0 \rho^\mu_\lambda- (\mu_2-\lambda_2-\ell-s) \theta +b_1 \alpha_2+ b_2 \alpha_1$$ 
	only if $(a_1,a_2)=(b_1+s,b_2+s).$ 
	Thus  along with \eqref{inv.S} it follows that the distinct irreducible components of $\cal F_{\lambda,\mu}$ that are parametrized by  elements of $\bigcup\limits_{\ell=1}^{\mu_2-\lambda_2}{}_\ell{\mathbb S}_{\lambda,\mu}$ are $V(\rho^\lambda_\mu+ w_0 \rho^\mu_\lambda- \ell \theta +a \alpha_2+ b \alpha_1)$. Consequently,
	\begin{equation}\label{q.C.l.m}[\mathcal F_{\lambda, \mu}: V(\rho^\lambda_\mu+ w_0 \rho^\mu_\lambda- \ell \theta +a \alpha_2+ b \alpha_1)]_q =\sum\limits_{s=0}^{\min\{a, b\}}q^{|\mu|-\ell-s} \end{equation}
	where $(\ell,a,b)$ is a triple of non-negative integers contained in the set $C_{\lambda,\mu}$, defined in \thmref{littlewood coeff}(iii).
	
	\item[(iii).] Given triplet of integers $(\ell, a_1,a_2)$ with $1 \leq \ell \leq \mu_2-\lambda_2, \  (a_1,a_2) \in {}_\ell {\mathbb{S}_{\lambda,\mu}}$ and 
	$(j,b_1,b_2)$ with $0 \leq j \leq \lambda_{2}-1, \ (b_1,b_2) \in  {\mathbb{S}^{\rho^{\lambda}_\mu+j \omega_2,\rho^{\mu}_\lambda-j\omega_2}}$,
	$$\rho^\lambda_\mu+ w_0 \rho^\mu_\lambda- \ell \theta +a_1 \alpha_2+ a_2 \alpha_1= \rho^\lambda_\mu+ w_0 \rho^\mu_\lambda+ j \theta +b_2 \alpha_1+ b_1 \alpha_2 $$
	only if $(b_1,b_2)=(a_1-\ell-j,a_2-\ell-j)$ and in this case $a_1-a_2= b_1-b_2$.
	
	We claim that, if  $(a_1,a_2)\in {_\ell}{\mathbb S}_{\lambda,\mu}$, is such that
	$(a_1-\ell-j,a_2-\ell-j) \in \mathbb{S}^{\rho^\lambda+j \omega_2, \rho^\mu_\lambda-j\omega_2} $ for some $0\leq j\leq \lambda_2-1$ and $1\leq \ell\leq \mu_2-\lambda_2$, then
	$\lambda_2-\lambda_1\leq a_2-a_1\leq \mu_2-\mu_1$ and consequently,
	\begin{equation}\label{last_char}
		[\mathcal F_{\lambda, \mu}: V(\rho^\lambda_\mu+ w_0 \rho^\mu_\lambda-\ell \theta +a_2 \alpha_1+ a_1 \alpha_2 )]_q =  \sum\limits_{s=0}^{\min\{a_{1}, a_{2}\}}q^{\mu_1+\lambda_2+\ell-s}.
	\end{equation}
	
	By definition, $$\begin{aligned}
		\mathbb{S}^{\rho^\lambda+j \omega_2, \rho^\mu_\lambda-j\omega_2} &=\{ (a, b) \in \mathbb Z_+^2: 0 \leq a \leq \mu_1, 0 \leq b \leq \lambda_2-j, \lambda_2-\lambda_1-j \leq b-a \leq \mu_2-\mu_1+j\}.\end{aligned}$$
			Suppose contrary to our claim there exists $(a_1,a_2)\in {_\ell}{\mathbb S}_{\lambda,\mu}$ such that $(a_1-\ell-j,a_2-\ell-j) \in \mathbb{S}^{\rho^\lambda+j \omega_2, \rho^\mu_\lambda-j\omega_2} $ and  
			$$a_2-a_1\in \{ \mu_2-\mu_1+k,\, \lambda_2-\lambda_1-k: 1\leq k \leq j\}.$$ 
			Since by definition of ${}_{_\ell}\mathbb{S}_{\lambda,\mu}$, for $(a_1,a_2)\in {}_{_\ell}\mathbb{S}_{\lambda,\mu}$,  
			$\mu_1-\mu_2+\ell \leq a_1-a_2 \leq \lambda_1-\lambda_2-\ell$ it follows that  
			$$ \text{if }\, a_2-a_1=\mu_2-\mu_1+k, \text{ then } \mu_1-\mu_2+\ell \leq \mu_1-\mu_2-k,$$ implying $\ell \leq -k$. This contradicts the fact $\ell, k$ are both positive integers. Similarly, if $a_2-a_1=\lambda_2-\lambda_1-k$,  then  
			$\lambda_1-\lambda_2+k \leq \lambda_1-\lambda_2-\ell$ implying $k \leq -\ell$ which contradicts that positivity of both $\ell$ and $k$. This proves our claim.
			
			
			\noindent Since $1\leq \ell \leq \mu_2-\lambda_2$, putting $\lambda= \mu_2-\lambda_2-r$ where $0 \leq r < \mu_2-\lambda_2$, \eqref{last_char} can be written as 
			$$[\mathcal F_{\lambda, \mu}: V(\rho^\lambda_\mu+ w_0 \rho^\mu_\lambda-(\mu_2-\lambda_2-r) \theta +a_2 \alpha_1+ a_1 \alpha_2 )]_q =  \sum\limits_{s=0}^{\min\{a_{1}, a_{2} \}}q^{|\mu|-r-s}.$$
			Summing up we get, 
			\begin{equation}\label{q.B.l.m}[\mathcal F_{\lambda, \mu}: V(\lambda+ w_0 \mu+r \theta +a_2 \alpha_1+ a_1 \alpha_2 )]_q =  \sum\limits_{s=0}^{\min\{a_{1}, a_{2} \}}q^{|\mu|-r-s} \quad 
				\text{for } (r,b_1,b_2) \in B_{\lambda,\mu}\end{equation} 
			where $B_{\lambda,\mu}$ is defined in \thmref{littlewood coeff}(iii).
			
			\noindent Putting $q=1$ in \eqref{l.P.inv},  \eqref{q.B.l.m}, \eqref{q.C.l.m}  and observing that in this case 
			$$[\mathcal F_{\lambda, \mu}: V(\nu)]_{q=1} =0, \qquad \text{when } \nu\notin A^{\rho^\lambda_\mu,\rho^\mu_\lambda}\cup B_{\lambda,\mu}\cup C_{\lambda,\mu}$$ we obtain \thmref{littlewood coeff}(iii). This completes the proof of \thmref{littlewood coeff}.

			\section{A proof of Saturated Tensor Product Theorem in type {$A_2$} using \thmref{littlewood coeff}}\label{section 8}
			
			\subsection{} Given $n,k\in \mathbb N$, let $p(n,r+1)$ be the set of all partitions of $n$ having $r+1$ parts. One can attach to a partition $\bon=(n_1\geq n_2\geq \cdots\geq n_{r+1})\in p(n,r+1)$ a dominant integral weight of $\lie{sl}_{r}$ given by 
			$\lambda_\bon = \sum\limits_{i=1}^{r}(n_i-n_{i+1})\omega_i$. It is known that two partitions $\bon,\bon'$ having $r+1$ parts correspond to the same dominant integral weight of $\lie{sl}_r$ provided, $n_1-n_1'=\cdots=n_{r+1}-n_{r+1}'$ or equivalently if $\sum\limits_{i=1}^{r+1} n_i = \sum\limits_{i=1}^{r+1}n_i'\mod r+1.$ 
			Conversely, one can associate with $\lambda=\sum\limits_{i=1}^r\lambda_i\omega_i\in P^+$ a partition $\bop_\lambda= (\sum\limits_{i=1}^r\lambda_i\geq\sum\limits_{i=2}^r\lambda_i\geq \cdots \geq \lambda_r )$ of $|\bop_\lambda|=\sum\limits_{j=1}^r j\lambda_j$.
			Using the ``honeycomb model,'' the saturation problem for tensor products was solved  in \cite{KT}, which can be rephrased as follows :
			\begin{saturation} \label{saturation} 
				Let  $\bon\in p(n,k+1)$, $\bom\in p(m,k+1)$ and $\bor\in p(n+m,k+1)$. Then 
				$c^{\lambda_\bor}_{\lambda_\bon,\lambda_\bom}\neq 0$ if and only if $c^{N\lambda_\bor}_{N\lambda_\bon, N\lambda_\bom}\neq 0$ for some $N>1$. 
			\end{saturation}  
		
		In this section, we show that our characterization of $c^\nu_{\lambda,\mu}$  can be used to prove the saturation theorem in the case when $r=2$. From \thmref{littlewood coeff} it is easy to see that if $c^\nu_{\lambda,\mu}\neq 0$, then for all $N>1$, $c^{N\nu}_{N\lambda,N\mu}\neq 0$. In this section, we give a case-by-case proof of the converse statement of \thmref{saturation} for $r=2$. The following is an important observation that we use in the proof: 
		if $\lambda=\lambda_\bon$, $\mu=\lambda_\bom$ and $\nu=\lambda_\bor$ are such that $\sum\limits_{i=1}^{k+1} n_i+m_i =\sum\limits_{i=1}^{k+1}r_i$, then $|\bop_\lambda|+|\bop_\mu|\equiv |\bop_\nu|\mod 3$. 
		
		With this in mind, we prove the following result.
		\begin{thm} Given a triple $(\lambda,\mu,\nu)$ of dominant integral weights of $\mathfrak{sl}_3(\mathbb C)$ satisfying the condition $|\bop_\lambda|+|\bop_\mu|\equiv|\bop_\nu|\mod 3$, then $c^{\mu}_{\lambda,\mu}\neq 0$ if and only if 
			$c^{N\nu}_{N\lambda,N\mu}\neq 0$ for some $N>1$.
		\end{thm}
		
		\subsubsection{Proof of \thmref{saturation}: Case 1.} Suppose $(\lambda,\mu_1\omega_1)\in P^+(\lambda+\mu_1\omega_1,2)$ is a partition of first kind and 
		$\nu\in P^+$ is such that $c^{N\nu}_{N\lambda,N\mu_1\omega_1}\neq 0$ for some $N>1$.
		Then by \thmref{littlewood coeff}, there exists a pair of non-negative integers $(j,a_j)$ satisfying  
		$0\leq j\leq N\mu_1, \quad \max\{0, N\mu_1-j-N\lambda_2\} \leq a_j \leq N\mu_1-j $ such that
		$$N \nu= N\lambda+w_0(N\mu)+j\theta+a_j\alpha_2$$ Comparing the coefficients of $\omega_1$ and $\omega_2$ we get, \begin{equation}\label{a-p.1} N\nu_1 = N\lambda_1+j-a_j, \quad 
			N\nu_2 =N\lambda_2-N\mu_1+j+2a_j, \end{equation} which implies that,
		\begin{equation} \label{a-p} 
			j  = \frac{N}{3}(2 \nu_1-2\lambda_1+\nu_2-\lambda_2+\mu_1), \quad 
			a_j= \frac{N}{3}(\nu_1-\lambda_1-\mu_1-\nu_2+\lambda_2).\end{equation}
		Using the condition $\bop_{N\lambda}+\bop_{N\mu}\equiv\bop_{N\nu}\mod 3$, we see that  since  $N(\lambda_1+2\lambda_2+\mu_1-\nu_1-2\nu_2)\equiv 0 \mod 3$,
		we have,  $N(\lambda_1+\mu_1-\nu_1)=2N(\nu_2-\lambda_2)+3r$ for some $r\in \mathbb Z$. Substituting this value of $N(\lambda_1+\mu_1-\nu_1)$ in \eqref{a-p}, we see that $a_j= -\frac{N}{3}(2\lambda_2-2\nu_2+3r-\nu_2+\lambda_2)$ for some $r\in \mathbb Z$. Hence it follows that $N$ divides $a_j$, and consequently, it follows from \eqref{a-p.1} that $N$ divides $j$. As there exits $a_j',j'$ such that $a_j=Na_j', j=Nj'$. As $0\leq Nj'\leq N\mu_1$ and $Nj'-N\lambda_2\leq Na_j'\leq Nj'$,  $0\leq j'\leq \mu_1$ and $j'-\lambda_2\leq a_j'\leq j'$, therefore by \thmref{littlewood coeff}, 
		$c^{\nu}_{\lambda,\mu}\neq 0,$ which proves the saturation theorem in this case.
		\subsubsection{Proof of \thmref{saturation}: Case  2.} Suppose $(\lambda,\mu)$ is a partition of $\lambda+\mu$ of first kind with $\mu_i>0$ for $i=1,2$ and there exists $N>1$ such that $c^{N\nu}_{N\lambda, N\mu}\neq 0$. By \thmref{littlewood coeff}, two cases can occur:\\
		i. $N\nu=N\lambda+w_0N\mu+(\ell+j)\theta+(\mu_2-j)\alpha_1+a\alpha_2,\,  \text{with}\, (j,\ell,a)\in A^{N\lambda,N\mu}$ \\
		ii. $N\nu=	N\lambda+w_0N\mu+j\theta+b\alpha_1+a\alpha_2,\, \text{with}\, (j,a,b)\in B^{N\lambda,N\mu}$ \\
		\noindent {\it{Case 2 }}(i). If $N\nu =N\lambda+w_0N\mu+(\ell+j)\theta+(\mu_2-j)\alpha_1+a\alpha_2=N\lambda+(N\mu_2+\ell-(a+j))\omega_1-(N|\mu|-\ell-2(a+j)) \omega_2,$\\ comparing the coefficients of the fundamental weights we get,
		\begin{equation} \label{first.2.l.a}N\nu_1-N\lambda_1-N\mu_2= \ell- (a+j), \quad \quad N\nu_2-N\lambda_2+N|\mu|= \ell+2(a+j). \end{equation} 
		Hence we have 
		\begin{equation}\label{first.2}
			\ell=\frac{N}{3}(2 \nu_1+\nu_2-2\lambda_1-\lambda_2-2\mu_2+|\mu|), \quad 
			a+j= \frac{N}{3}(\nu_2-\nu_1+\lambda_1-\lambda_2+\mu_2+|\mu|).
		\end{equation}
		Using $|\bop_{N\lambda}|+|\bop_{N\mu}|\equiv |\bop_{N\nu}|\mod 3$, we see that since $N(|\mu|+\mu_2+|\lambda|+\lambda_2-|\nu|-\nu_2) \equiv 0\mod 3$,
		$$\begin{aligned}
			2 \nu_1+\nu_2-2\lambda_1-\lambda_2-2\mu_2+|\mu|= 3(|\nu|-|\lambda|-\mu_2)+|\mu|+\mu_2+|\lambda|+\lambda_2-|\nu|-\nu_2 &\equiv 0\mod 3,\\
		\end{aligned}$$
		From \eqref{first.2} it thus follows that $N$ divides $\ell$ and hence by \eqref{first.2.l.a}, $N$ divides $a+j.$ \\
		For $j=0$, this implies that $N$ divides both $\ell$ and $a$, hence using the same arguments as in case 1, we see that the theorem holds in this case.\\
		For $j\neq 0$, $a \in \{N\mu_1-\ell, N|\mu|-N\lambda_2-\ell-2j\} \cap \mathbb Z_+$ for $0 \leq \ell \leq N\mu_1$. Since $N$ divides $\ell$, we see that if $a=N\mu_1-\ell$, then $N$ divides $a$ and since $N$ divides 
		$a +j$, it follows that in this case $N$ divides $j$. On the other hand if 
		$a=N|\mu|-N\lambda_2-\ell-2j,$  using the fact that $N$ divides $\ell$ and $a +j$, we see that $N$ divides  $N|\mu|-N\lambda_2-\ell-2j+j=N|\mu|-N\lambda_2-\ell-j$ which implies that in this case $N$ divides $j$. Now using the same arguments as in case 1, we see that  the theorem also holds for $j\neq 0$.\\
		
		\noindent {\it{Case 2 }}(ii). If $N\nu=	N\lambda+w_0N\mu+j\theta+b\alpha_1+a\alpha_2= (N\lambda_1-N\mu_2+j+2b-a)\omega_1+(N\lambda_{2}-N\mu_1+j-b+2a)\omega_2$, 
		for $0 \leq j \leq N\mu_2-1$, $(a,b) \in \mathbb S_{ninv}(N\lambda+ j \omega_2, N\mu-j \omega_2), \quad  b-a \in \{N(\lambda_2-\mu_1)+j, N(\mu_2-\lambda_1)-j\},$ comparing the coefficients of the fundamental weights we get, 
		\begin{equation} \label{first.2(ii)}
			N\nu_1-N\lambda_1+N\mu_2 =j-a+2b, \quad N\nu_2-N\lambda_2+N\mu_1=j+2a-b.
		\end{equation} Hence we have
		\begin{equation}\label{first}\begin{array}{c} a+j=\frac{N}{3}(2 \nu_2+2\mu_1-2\lambda_2+\nu_1-\lambda_1+\mu_2), \quad b+j=\frac{N}{3}(\nu_2+\mu_1+2\nu_1-2\lambda_1+2\mu_2-\lambda_2),\\
				\text{and}\quad \quad
				b-a= \frac{N}{3}(-\nu_2-\mu_1+\nu_1-\lambda_1+\mu_2+\lambda_2).\end{array}
		\end{equation}
		Since $   |\nu|+\nu_2-|\lambda|-\lambda_2-|\mu|-\mu_2 \equiv 0 \mod 3$,
		$$ |\nu|+\nu_2-|\lambda|-\lambda_2-|\mu|-\mu_2+3(\mu_1+\mu_2)=
		\nu_1+2\nu_2-\lambda_1-2\lambda_2+2\mu_1+\mu_2\equiv 0\mod 3.$$ 
		It now follows from \eqref{first} that $N$ divides $a+j$. Thus using \eqref{first.2(ii)}, we see that $N$ also divides $b+j$. \\
		This implies that for $j=0$, $N$ divides both $a$ and $b$. Therefore  similar arguments as in Case 1, show that the theorem holds in this case.\\
		When $j\neq 0$, $b-a \in \{N(\lambda_2-\mu_1)+j, N(\mu_2-\lambda_1)-j\}. $\\
		If $b-a= N(\lambda_2-\mu_1)+j$, then by \eqref{first}, $N(\lambda_2-\mu_1)+j= \frac{N}{3}(-\nu_2-\mu_1+\nu_1-\lambda_1+\mu_2+\lambda_2)$, which gives, 
		$$3j= {N}(\nu_1-\nu_2+2\mu_1+\mu_2-\lambda_1-2\lambda_2)= 3N(-\nu_2+\mu_1+\mu_2)+N(|\nu|+\nu_2-|\mu|-\mu_2-|\lambda|-\lambda_2).$$ Using the condition that 
		$|N\lambda|+N\lambda_2+N|\mu|+N\mu_2-N|\nu|-N\nu\equiv 3,$  it follows  
		that $j\equiv 0\mod N$ and hence $a\equiv 0\mod N$ and $b\equiv 0\mod N$.\\
		Similarly if $b-a=N(\mu_2-\lambda_1)-j$, then by \eqref{first}, 
		$N(\mu_2-\lambda_1)-j= \frac{N}{3}(-\nu_2-\mu_1+\nu_1-\lambda_1+\mu_2+\lambda_2)$, implying
		$$3j=N(2\mu_2-2\lambda_1+\nu_2-\nu_1+\mu_1-\lambda_2)=N(|\lambda|+\lambda_2+|\mu|+\mu_2-|\nu|-\nu_2+3(\lambda_1+\lambda_2-\nu_2)).$$ Using the condition that 
		$|N\lambda|+N\lambda_2+N|\mu|+N\mu_2-N|\nu|-N\nu\equiv 3,$  it follows that in this case also $j\equiv 0\mod N$ and hence $a\equiv 0\mod N$ and $b\equiv 0\mod N$. Now similar arguments as used in case 1, show that the theorem also holds when $(\lambda,\mu)$ is a partition of first kind with $\mu_i>0$ for $i=1,2.$ 
		\label{Case.2.sat}
		\subsubsection{Proof of \thmref{saturation}: Case 3.} Suppose $(\lambda,\mu)$ is a partition of $\lambda+\mu$ of second kind and there exists some $N>1$ such that $c^{N\nu}_{N\lambda,N\mu} \neq 0$. By \thmref{littlewood coeff}, 
		three cases can occur:-\\
		i. $N\nu= \rho^{N\lambda}_{N\mu}+ w_0 \rho^{N\mu}_{N\lambda}- \ell \theta +a \alpha_2+ b \alpha_1, (\ell,a,b) \in C_{N\lambda,N\mu}$\\
		ii. $N\nu= N\lambda+ w_0 N\mu+j \theta +b_1 \alpha_2+ b_2 \alpha_1, (j,b_1,b_2) \in {B}_{N\lambda,N\mu}$ \\
		iii. $ N\nu =\rho^{N\lambda}_{N\mu}+w_0\rho^{N\mu}_{N\lambda}+(\ell+j)\theta+(\lambda_2-j)\alpha_1+a \alpha_2, (j,\ell,a) \in {A}_{N\lambda,N\mu}$\\
		As Case 3(ii) and (iii) follow from \secref{Case.2.sat}, 
		here we consider \textit{Case 3}(i). \\
		\noindent 
		If $N\nu= \rho^{N\lambda}_{N\mu}+ w_0 \rho^{N\mu}_{N\lambda}- \ell \theta +a \alpha_2+ b \alpha_1,$ for $(\ell,a,b) \in C_{N\lambda,N\mu}$, then
		$$N \nu = (N\lambda_1-N\lambda_2-\ell-a+2b)\omega_1+(N\mu_2-N\mu_1-\ell+2a-b)\omega_2.$$
			Comparing the coefficients of the fundamental weights we get,
			$$
			N \nu_1= N\lambda_1-N\lambda_2-a+2b-\ell, \quad N\nu_2= N\mu_2-N\mu_1-b+2a-\ell.$$ 
			Hence we have,
			\begin{equation}\label{second_inv}\begin{array}{c} a-\ell = \frac{N}{3}(\nu_1+2\nu_2-\lambda_1+\lambda_2+2\mu_1-2\mu_2), \quad b-\ell =\frac{N}{3}(2 \nu_1+\nu_2-2\lambda_1+2\lambda_2-\mu_2+\mu_1) \\ 
					\text{ and } b-a_= \frac{N}{3}(\nu_1-\nu_2-\lambda_1+\lambda_2-\mu_1+\mu_2). \end{array}
			\end{equation}
			Using $|\bop_{N\lambda}|+|\bop_{N\mu}|\equiv |\bop_{N\nu}|\mod 3$, we see that since $N(|\mu|+\mu_2+|\lambda|+\lambda_2-|\nu|-\nu_2) \equiv 0\mod 3$,
			$$\begin{aligned}
				&	2 \nu_1+\nu_2-2\lambda_1+2\lambda_2-\mu_2+\mu_1 = 3(\nu_1+\nu_2-\lambda_1-\mu_2)-|\nu|-\nu_2+|\lambda|+\lambda_2+|\mu|+\mu_2
				\equiv 0 \mod 3,\\
				& \text{ and } \nu_1+2\nu_2-\lambda_1+\lambda_2+2\mu_1-2\mu_2 =
				|\nu|+\nu_2-|\lambda|-\lambda_2-|\mu|-\mu_2+3 (\mu_1+\lambda_2) \equiv 0\mod 3. 
			\end{aligned}$$ Therefore, it follows from \eqref{second_inv} that $N$ divides both $b-\ell$ and $a-\ell$.
			\noindent Thus when $\ell =0$, $N$ divides $a$ and $b$ and from similar arguments as in case 1 it follows that  $\nu= \rho^\lambda_\mu+ w_0 \rho^\mu_\lambda+ \frac{a}{N} \alpha_2+\frac{b}{N} \alpha_1 $ for $(0,\frac{a}{N},\frac{b}{N}) \in C_{\lambda,\mu} $. Hence the theorem holds in this case.\\
			When $\ell \neq 0$, $b-a \in \{ N \mu_2-N\lambda_1-\ell, N \lambda_2-N\mu_1+\ell\}.$\\ If $b-a = N \mu_2-N\lambda_1-\ell $, using \eqref{second_inv}, we get
			\begin{equation}\label{r+1_inv} \begin{aligned}
					\ell=	N \mu_2-N\lambda_1-\frac{N}{3} (\nu_1-\nu_2-\lambda_1+\lambda_2-\mu_1+\mu_2) = \frac{N}{3}(2 \mu_2-2 \lambda_1-\nu_1+\nu_2-\lambda_2+\mu_1).    
			\end{aligned}\end{equation}
			Since $ 2 \mu_2-2 \lambda_1-\nu_1+\nu_2-\lambda_2+\mu_1 = |\mu|+\mu_2+|\lambda|+\lambda_2-|\nu|-\nu_2-3(\lambda_1+\lambda_2-\nu_2)\equiv 0\mod 3,
			$ it follows from \eqref{r+1_inv} that $N$ divides $\ell $. Hence, using \eqref{second_inv} we see that $N$ divides $a$ and $b$. Taking $N a'= a, \ N b'= b,$ and $N\ell'=\ell$, we get $\nu=  \rho^\lambda_\mu+ w_0 \rho^\mu_\lambda- \ell' \theta + a' \alpha_2+b' \alpha_1 $ for $(\ell',a',b') \in C_{\lambda,\mu}$, which proves the theorem in this case. Now, similar arguments can be used for $b-a= N \lambda_2-N \mu_1+ \ell$. This completes the proof of the theorem.\qed
			
			
			Using \thmref{littlewood coeff}, the following result follows from the proof of \thmref{saturation}.
			\begin{cor} Given a triple $(\lambda,\mu,\nu)$ of dominant integral weights of $\mathfrak{sl}_3(\mathbb C)$ satisfying the condition $|\bop_\lambda|+|\bop_\mu|\equiv|\bop_\nu|\mod 3$, $c^{N\nu}_{N\lambda,N\mu} = N(c^\nu_{\lambda,\mu}-1)+1$, for all $N\in \mathbb N.$ 
			\end{cor}

			\end{document}